\documentclass[a4paper,11pt,reqno]{amsart}

\usepackage{amsmath}
\usepackage{amssymb}
\usepackage{amsthm}
\usepackage{amsfonts}
\usepackage{graphicx}
\usepackage{psfrag}
\usepackage{color}
\usepackage{bbm}
\usepackage{url}
\usepackage{ifthen}
\usepackage[labelfont={footnotesize,bf},textfont=footnotesize]{caption}
\usepackage[labelfont=footnotesize,textfont=footnotesize]{subfig}
\usepackage[T1]{fontenc}

\DeclareMathOperator{\convop}{conv}

\DeclareMathOperator{\affop}{aff}
\DeclareMathOperator{\linop}{lin}

\newcommand{\ints}[1]{[{#1}]}
\newcommand{\setdef}[2]{\{{#1}\,:\,{#2}\}}
\newcommand{\suchthat}{\;:\;}
\newcommand{\symgr}[1]{{\mathfrak{S}}_{#1}}
\newcommand{\cyclgr}[1]{{\mathfrak{C}}_{#1}}
\newcommand{\R}{\mathbbm{R}}
\newcommand{\Q}{\mathbbm{Q}}

\newcommand{\Z}{\mathbbm{Z}}
\newcommand{\zerovec}{\mathbf{0}}
\newcommand{\bigo}[1]{O({#1})}
\newcommand{\aff}[1]{\affop({#1})}
\newcommand{\lin}[1]{\linop({#1})}
\newcommand{\dispsum}{\displaystyle \sum}
\newcommand{\card}[1]{\lvert {#1} \rvert}
\newcommand{\basis}{\mathcal{B}}
\newcommand{\scal}[2]{\langle{#1},{#2}\rangle}

\newtheorem{definition}{Definition}
\newtheorem{theorem}{Theorem}
\newtheorem{proposition}[theorem]{Proposition}
\newtheorem{lemma}[theorem]{Lemma}
\newtheorem{corollary}[theorem]{Corollary}
\newtheorem{claim}{Claim}
\newtheorem{observation}{Observation}

\newtheoremstyle{plainNoItalics}{}{}{\normalfont}{}{\bfseries}{.}{ }{}
\theoremstyle{plainNoItalics}

\newtheorem*{remark*}{Remark}

\makeatletter
\newenvironment{myenumerate}{%
  \advance\@enumdepth \@ne%
  \edef\@enumctr{enum\romannumeral\the\@enumdepth}%
  \begin{list}{\csname label\@enumctr\endcsname}%
    {\def\makelabel##1{\hss\llap{\upshape##1}}%
    \usecounter{\@enumctr}%
      \setlength{\topsep}{1mm}%
      \setlength{\partopsep}{0mm}%
      \setlength{\parskip}{0mm}%
      \setlength{\parsep}{0mm}%
      \setlength{\itemsep}{0mm}%
      \settowidth{\labelwidth}{\csname label\@enumctr\endcsname}%
      \setlength{\labelsep}{2mm}%
      \setlength{\leftmargin}{0mm}%
      \addtolength{\leftmargin}{\labelwidth}%
      \addtolength{\leftmargin}{\labelsep}%
      \setlength{\itemindent}{0mm}}}%
  {\end{list}
}
\makeatother

\newenvironment{myitemize}{%
\begin{list}{$\circ$}%
{\setlength{\topsep}{1ex}%
\setlength{\partopsep}{0mm}%
\setlength{\parskip}{0ex}%
\setlength{\parsep}{0mm}%
\setlength{\itemsep}{0ex}%
\setlength{\labelwidth}{4mm}%
\setlength{\leftmargin}{0ex}%
\addtolength{\leftmargin}{\labelwidth}%
\addtolength{\leftmargin}{\labelsep}%
\setlength{\itemindent}{0mm}}}%
{\end{list}}

\newcounter{partslistcounter}
\newenvironment{partslist}{%
\setcounter{partslistcounter}{0}
\begingroup
\newcommand{\partsitem}{%
\stepcounter{partslistcounter}%
\ifthenelse{\thepartslistcounter>1}{\medskip\par}{}
\noindent\textit{Part~(\arabic{partslistcounter}):} }
}%
{\endgroup}

\graphicspath{{Pictures/}}

\makeatletter
\renewcommand{\subsection}{\@startsection{subsection}{2}%
  {\z@}%
  {.7\linespacing\@plus\linespacing}%
  {.5\linespacing}%
  {\normalfont\scshape\centering}}
\makeatother

\DeclareMathOperator{\orbiop}{O}
\DeclareMathOperator{\rowop}{row}
\DeclareMathOperator{\colop}{col}
\DeclareMathOperator{\diagop}{diag}
\DeclareMathOperator{\weightop}{weight}

\newcommand{\mat}[2]{\mathcal{M}_{{#1},{#2}}}
\newcommand{\matpack}[2]{\mathcal{M}^{\le}_{{#1},{#2}}}
\newcommand{\matcov}[2]{\mathcal{M}^{\ge}_{{#1},{#2}}}
\newcommand{\matpart}[2]{\mathcal{M}^{=}_{{#1},{#2}}}
\newcommand{\matmaxG}[3]{\mat{#1}{#2}^{\max}({#3})}

\newcommand{\row}[1]{\rowop_{#1}}
\newcommand{\col}[1]{\colop_{#1}}

\newcommand{\orbifullG}[3]{\orbiop_{#1,#2}(#3)}
\newcommand{\orbipackG}[3]{\orbiop^{\le}_{#1,#2}(#3)}
\newcommand{\orbipartG}[3]{\orbiop^{=}_{#1,#2}(#3)}
\newcommand{\orbicoverG}[3]{\orbiop^{\ge}_{#1,#2}(#3)}

\newcommand{\orbipack}[2]{\orbiop^{\leq}_{#1,#2}(\symgr{#2})}
\newcommand{\orbipart}[2]{\orbiop^{=}_{#1,#2}(\symgr{#2})}

\newcommand{\orbipartinds}[2]{\mathcal{I}_{{#1},{#2}}}

\newcommand{\diagcol}[2]{\langle{#1},{#2}\rangle}
\newcommand{\diagleq}[2]{\diagop^{\leq}({#1},{#2})}
\newcommand{\diaggeq}[2]{\diagop^{\geq}({#1},{#2})}

\newcommand{\setof}[1]{I^{#1}}
\newcommand{\pointof}[1]{\chi^{#1}}

\title{Packing and Partitioning Orbitopes}

\author{}
\date{11/22/2006}
\author[Kaibel and Pfetsch]{Volker Kaibel \and Marc E.~Pfetsch}
\address{Zuse Institute Berlin, Takustr.~7, 14195 Berlin, Germany
}
\email{[kaibel,pfetsch]@zib.de}
\thanks{Supported by the DFG Research Center \textsc{Matheon} in Berlin}

\begin{document}

\begin{abstract}
  We introduce \emph{orbitopes} as the convex hulls of $0/1$-matrices that
  are lexicographically maximal subject to a group acting on the columns.
  Special cases are packing and partitioning orbitopes, which arise from
  restrictions to matrices with at most or exactly one $1$-entry in each
  row, respectively. The goal of investigating these polytopes is to gain
  insight into ways of breaking certain symmetries in integer programs by
  adding constraints, e.g., for a well-known formulation of the graph
  coloring problem.

  We provide a thorough polyhedral investigation of packing and
  partitioning orbitopes for the cases in which the group acting on the
  columns is the cyclic group or the symmetric group.  Our main results are
  complete linear inequality descriptions of these polytopes by
  facet-defining inequalities. For the cyclic group case, the descriptions
  turn out to be totally unimodular, while for the symmetric group case,
  both the description and the proof are more involved. The
  associated separation problems can be solved in linear time.
\end{abstract}

\maketitle

% +++++++++++++++++++++++++++++++++++++++++++
% Introduction
% +++++++++++++++++++++++++++++++++++++++++++

\section{Introduction}

Symmetries are ubiquitous in discrete mathematics and geometry. They are
often responsible for the tractability of algorithmic problems and for the
beauty of both the investigated structures and the developed methods. It is
common knowledge, however, that the presence of symmetries in integer
programs may severely harm the ability to solve them. The reasons for this
are twofold. First, the use of branch-and-bound methods usually leads to an
unnecessarily large search tree, because equivalent solutions are found
again and again. Second, the quality of LP relaxations of such programs
typically is extremely poor.

A classical approach to ``break'' such symmetries is to add constraints
that cut off equivalent copies of solutions, in hope to resolve these
problems.  There are numerous examples of this in the literature; we will
give a few references for the special case of graph coloring below. Another
approach was developed by Margot~\cite{Mar02,Mar03}. He studies a
branch-and-cut method that ensures to investigate only one representative
of each class of equivalent solutions by employing methods from
computational group theory. Furthermore, the symmetries are also used to
devise cutting planes. Methods for symmetry breaking in the context of
constraint programming have been developed, for instance, by Fahle,
Schamberger, and Sellmann~\cite{FahSS01b} and Puget~\cite{Pug05}.

The main goal of this paper is to start an investigation of the polytopes
that are associated with certain symmetry breaking inequalities.  In order
to clarify the background, we first discuss the example of a well-known
integer programming (IP) formulation for the graph coloring
problem.\smallskip

Let $G = (V, E)$ be a loopless undirected graph without isolated nodes. A
\emph{(vertex) coloring} of $G$ using at most~$C$ colors is an assignment
of colors $\{1, \dots, C\}$ to the nodes such that no two adjacent nodes
receive the same color. The \emph{graph coloring} problem is to find a
vertex coloring with as few colors as possible. This is one of the
classical NP-hard problems~\cite{GarJ79}. It is widely believed to be among
the hardest problems in combinatorial optimization. In the following
classical IP formulation, $V = \{1, \dots, n\}$ are the nodes of~$G$
and~$C$ is some upper bound on the number of colors needed.
\begin{equation}
\label{eq:intro:model}
\begin{array}{lr@{\;}ll@{\qquad}l}
  \min & \dispsum_{j=1}^{C} y_j &&&\\[3ex]
       & x_{ij} + x_{kj}         & \leq y_j & \{i,k\} \in E,\; j \in \{1,\dots,C\} & \text{(i)}\\
       & \dispsum_{j=1}^{C} x_{ij} & =  1     & i \in V & \text{(ii)}\\
       & \multicolumn{2}{r}{x_{ij} \in \{0,1\}} & i \in V,\; j \in \{1,\dots,C\} & \text{(iii)}\\
       & \multicolumn{2}{r}{y_j \in \{0,1\}} & j \in \{1, \dots, C\} & \text{(iv)}
\end{array}
\end{equation}
In this model, variable $x_{ij}$ is~$1$ if and only if color~$j$ is
assigned to node~$i$ and variable~$y_j$ is~$1$ if color~$j$ is used.
Constraints~(i) ensure that color~$j$ is assigned to at most one of the two
adjacent nodes~$i$ and~$k$; it also enforces that~$y_j$ is~$1$ if color~$j$
is used, because there are no isolated nodes. Constraints~(ii) guarantee
that each node receives exactly one color.

It is well known that this formulation exhibits symmetry: Given a solution
$(x, y)$, any permutation of the colors, i.e., the columns of~$x$ (viewed
as an $n \times C$-matrix) and the components of~$y$, results in a
valid solution with the same objective function value. Viewed abstractly,
the symmetric group of order~$C$ acts on the solutions $(x,y)$ (by
permuting the columns of~$x$ and the components of~$y$) in such a way that
the objective function is constant along every orbit of the group action.
Each orbit corresponds to a symmetry class of feasible colorings of the
graph. Note that ``symmetry'' here always refers to the symmetry of
permuting colors, not to symmetries of the graph.

The weakness of the LP-bound mentioned above is due to the fact that the
point $(x^\star,y^\star)$ with $x^\star_{ij} = 1/C$ and $y^\star_j = 2/C$
is feasible for the LP relaxation with objective function value~$2$.  The
symmetry is responsible for the feasibility of $(x^\star,y^\star)$,
since~$x^\star$ is the barycenter of the orbit of an arbitrary~$x \in
\{0,1\}^{n \times C}$ satisfying (ii) in~\eqref{eq:intro:model}.

% Furthermore, the value of the LP relaxation is~$2$:
% a feasible solution is given by setting $x_{ij} = 1/C$ and $y_j = 2/C$ for
% all~$i$ and~$j$. A lower bound is obtained by adding constraints~(i) for
% all~$C$ and an arbitrary edge~$\{i,k\}$, which yields
% \[
% \sum_{j=1}^C y_j \geq \sum_{j=1}^C x_{ij} + \sum_{j=1}^C x_{kj} = 2,
% \]
% by using equations~(ii) for~$i$ and~$j$.

It turned out that the symmetries make the above IP-formulation for the
graph coloring problem difficult to solve. One solution is to
develop different formulations for the graph coloring problem. This line
has been pursued, e.g., by Mehrotra and Trick~\cite{MehT96}, who devised a
column generation approach. See Figueiredo, Barbosa, Maculan, and de
Souza~\cite{FigBMS02} and Cornaz~\cite{Cor06} for alternative models.

Another solution is to enhance the IP-model by additional inequalities that
cut off as large parts of the orbits as possible, keeping at least one
element of each orbit in the feasible region. M\'{e}ndez-D\'{i}az and
Zabala~\cite{DiaZ06} showed that a branch-and-cut algorithm using this kind
of symmetry breaking inequalities performs well in practice. The polytope
corresponding to~\eqref{eq:intro:model} was investigated by Camp\^{e}lo,
Corr\^{e}a, and Frota~\cite{CamCF04} and Coll, Marenco,
M\'{e}ndez-D\'{i}az, and Zabala~\cite{CollMDZ02}.  Ramani, Aloul, Markov,
and Sakallah~\cite{RamAMS04} studied symmetry breaking in connection with
SAT-solving techniques to solve the graph coloring problem.

The strongest symmetry breaking constraints that M\'{e}ndez-D\'{i}az and
Zabala~\cite{DiaZ01,DiaZ06} introduced are the inequalities
\begin{equation}\label{eq:SymmetryBreak}
x_{ij} - \sum_{k = 1}^{i-1} x_{k,j-1} \leq 0,\quad\text{ for all }i \text{
  and }j \geq 2.
\end{equation}
From each orbit, they cut off all points except for one representative that
is the maximal point in the orbit with respect to a lexicographic ordering.
A solution $(x,y)$ of the above IP-model is such a representative if and
only if the columns of~$x$ are in decreasing lexicographic order. We
introduce a generalization and strengthening of
Inequalities~\eqref{eq:SymmetryBreak} in Section~\ref{sec:IPFormulations}.
\medskip

\begin{figure}
  \centering
  \includegraphics[height=4cm]{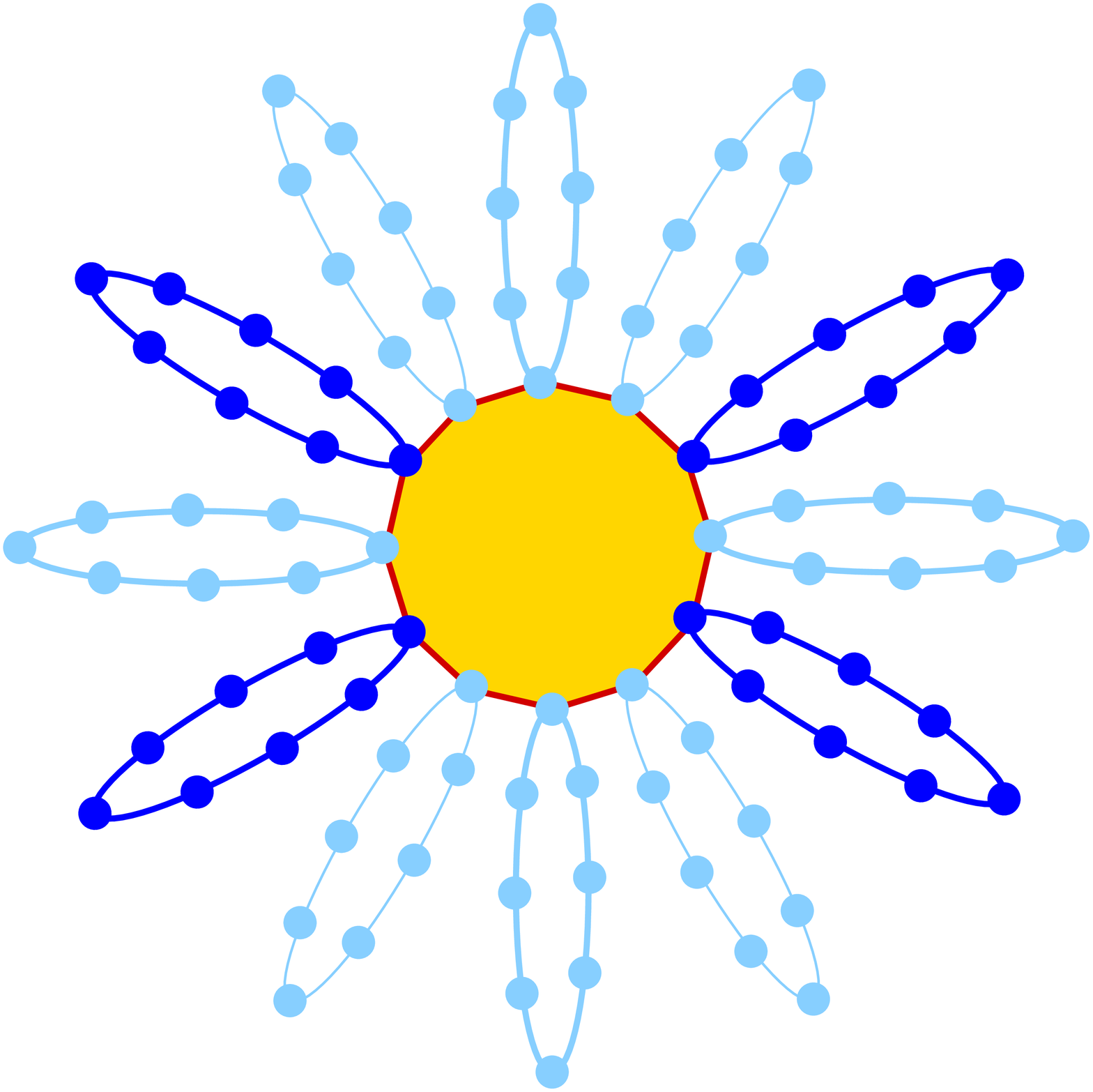}\hspace{2cm}
  \includegraphics[height=4cm]{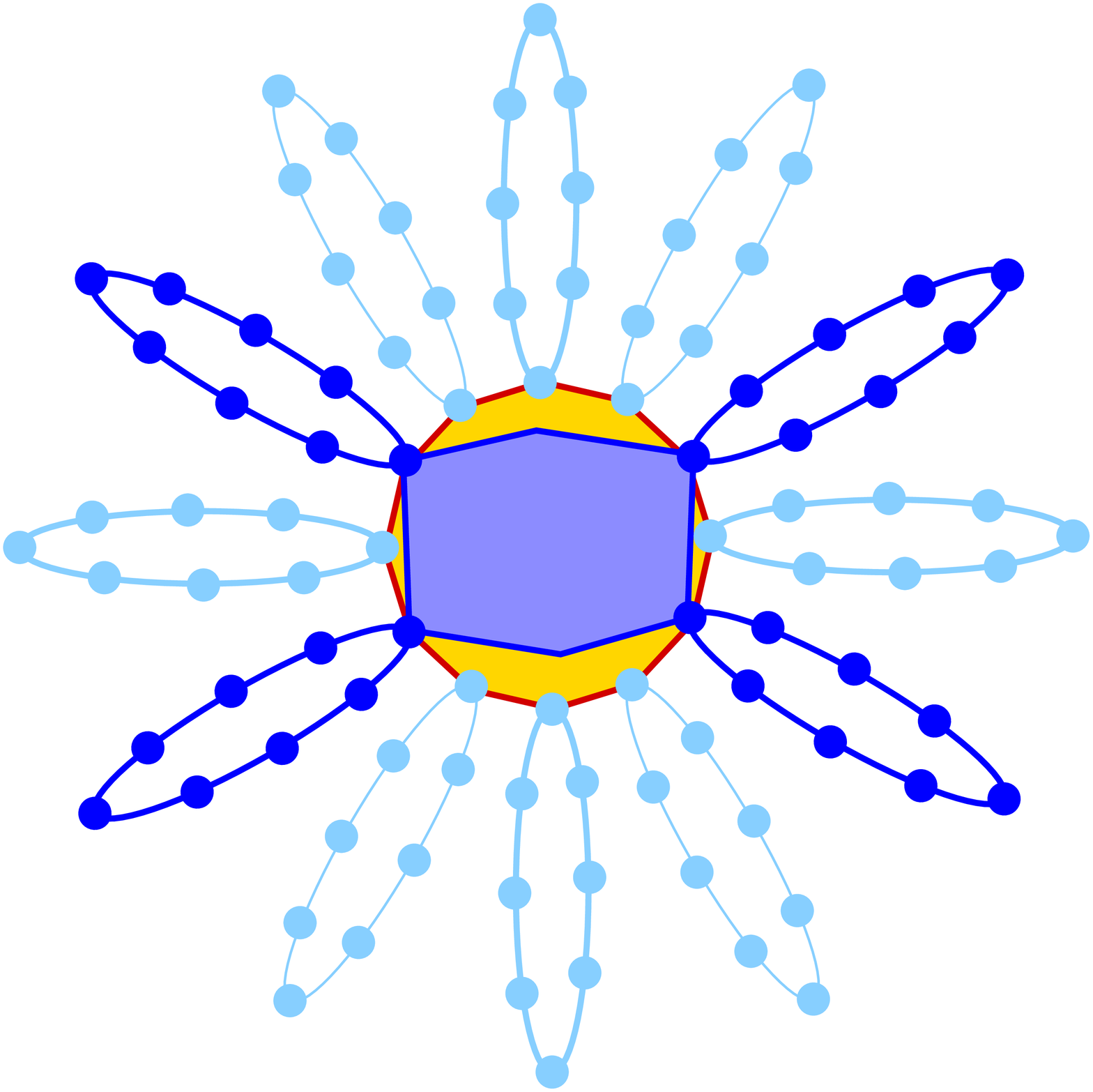}
  \caption{Breaking symmetries by orbitopes. The left figure illustrates an
    orbitope, i.e., the convex hull of the representatives of a large
    system of orbits. For a concrete problem, like graph coloring, only a
    subset of the orbits are feasible (the dark orbits). Combining a (symmetric)
    IP-formulation for the concrete problem with the orbitope removes the
    symmetry from the formulation (right figure).}
  \label{fig:GenericIllustration}
\end{figure}

Breaking symmetries by adding inequalities like~\eqref{eq:SymmetryBreak}
does not depend on the special structure of the graph coloring problem.
These inequalities single out the lexicographic maximal representative from
each orbit (with respect to the symmetric group acting on the columns) of
the whole set of all 0/1-matrices with exactly one $1$-entry per row.  The
goal of this paper is to investigate the structure of general ``symmetry
breaking polytopes'' like the convex hull of these representatives. We call
these polytopes \emph{orbitopes}.  The idea is that general knowledge on
orbitopes (i.e., valid inequalities) can be utilized for different
symmetric IPs in order to address both the difficulties arising from the
many equivalent solutions and from the poor LP-bounds.  In particular with
respect to the second goal, for concrete applications it will be desirable
to combine the general knowledge on orbitopes with concrete polyhedral
knowledge on the problem under investigation in oder to derive strengthened
inequalities.  For the example of graph coloring, we indicate that (and
how) this can be done in Section~\ref{sec:ClosingRemarks}.
Figure~\ref{fig:GenericIllustration} illustrates the geometric situation.

The case of a symmetric group acting on the columns is quite important. It
does not only appear in IP-formulations for the graph coloring problem, but
also in many other contexts like, e.g., block partitioning of
matrices~\cite{BorFM98}, $k$-partitioning in the context of frequency
assignment~\cite{Eis01}, or line-planning in public
transport~\cite{BorGP05}.  However, other groups are interesting as well.
For instance, in the context of timetabling in public transport
systems~\cite{SerU89}, cyclic groups play an important role.

We thus propose to study different types of orbitopes, depending on the
group acting on the columns of the variable-matrix and on further
restrictions like the number of $1$-entries per row being exactly
one (\emph{partitioning}), at most one (\emph{packing}), at least one
(\emph{covering}), or arbitrary (\emph{full}).\smallskip

The main results of this paper are complete and irredundant linear
descriptions of packing and partitioning orbitopes for both the symmetric
group and for the cyclic group acting on the columns of the
variable-matrix. We also provide (linear time) separation algorithms for
the corresponding sets of inequalities. While this work lays the
theoretical foundations on orbitopes, a thorough computational
investigation of the practical usefulness of the results will be the
subject of further studies (see also the remarks in
Section~\ref{sec:ClosingRemarks}).

The outline of the paper is as follows. In
Section~\ref{sec:BasicOrbitopes}, we introduce some basic notations and
define orbitopes. In Section~\ref{sec:OptimizingOverOrbitopes} we show that
optimization over packing and partitioning orbitopes for symmetric and
cyclic groups can be done in polynomial time. In
Section~\ref{sec:PackingAndPartitioningOrbitopesForCyclicGroups} we give
complete (totally unimodular) linear descriptions of packing and
partitioning orbitopes for cyclic groups.
Section~\ref{PackingAndPartitioningOrbitopesForSymmetricGroups} deals with
packing and partitioning orbitopes for symmetric groups, which turn out to
be more complicated than their counterparts for cyclic groups. Here,
besides (strengthenings of) Inequalities~\eqref{eq:SymmetryBreak}, one
needs exponentially many additional inequalities, the ``shifted column
inequalities'', which are introduced in
Section~\ref{sec:ShiftedColumnInequalities}. We show that the corresponding
separation problem can be solved in linear time, see
Section~\ref{sec:SeparationAlgorithmForSCIs}.
Section~\ref{sec:CompleteInequalityDescriptions} gives a complete linear
description, and Section~\ref{sec:Facets} investigates the facets of the
polytopes. We summarize the results for symmetric groups in
Section~\ref{SummaryOfResultsSymmetricGroup} for easier reference. Finally,
we close with some remarks in Section~\ref{sec:ClosingRemarks}.

% +++++++++++++++++++++++++++++++++++++++++++
% Definitions and terminology
% +++++++++++++++++++++++++++++++++++++++++++

\section{Orbitopes: General Definitions and Basic Facts}
\label{sec:BasicOrbitopes}

We first introduce some basic notation.  For a positive integer~$n$, we
define $\ints{n} := \{1,2,\dots,n\}$. We denote by $\zerovec$ the 0-matrix
or 0-vector of appropriate sizes. Throughout the paper let~$p$ and~$q$ be
positive integers. For $x \in \R^{\ints{p} \times \ints{q}}$ and $S
\subseteq \ints{p} \times \ints{q}$, we write
\[
x(S) := \sum_{(i,j) \in S} x_{ij}.
\]
For convenience, we use $S - (i,j)$ for $S \setminus \{(i,j)\}$ and $S +
(i,j)$ for $S \cup \{(i,j)\}$, where $S \subseteq \ints{p} \times \ints{q}$
and $(i,j) \in \ints{p} \times \ints{q}$. If~$p$ and~$q$ are clear from the
context, then $\row{i} := \{(i,1), (i,2), \dots, (i,q)\}$ are the entries
of the $i$th row.

Let $\mat{p}{q} := \{0,1\}^{\ints{p} \times \ints{q}}$ be the set of
$0/1$-matrices of size $p \times q$. We define
\begin{myitemize}
  \item $\matpack{p}{q} := \setdef{x \in \mat{p}{q}}{x(\row{i}) \leq 1 \text{ for all }i}$
  \item $\matpart{p}{q} := \setdef{x \in \mat{p}{q}}{x(\row{i}) = 1 \text{ for all }i}$
  \item $\matcov{p}{q} := \setdef{x \in \mat{p}{q}}{x(\row{i}) \geq 1 \text{ for all }i}$.
\end{myitemize}
Let $\prec$ be the lexicographic ordering of~$\mat{p}{q}$ with respect to
the ordering
\[
(1,1) < (1,2) < \dots < (1,q) < (2,1) < (2,2) < \dots < (2,q) < \dots < (p,q)
\]
of matrix positions, i.e., $A \prec B$ with $A = (a_{ij}), B = (b_{ij}) \in
\mat{p}{q}$ if and only if $a_{k\ell} < b_{k\ell}$, where $(k,\ell)$ is the
first position (with respect to the ordering above) where~$A$ and~$B$
differ.

Let $\symgr{n}$ be the group of all permutations
of~$\ints{n}$ (\emph{symmetric group}) and let~$G$ be a subgroup
of~$\symgr{q}$, acting on~$\mat{p}{q}$ by permuting columns. Let
$\matmaxG{p}{q}{G}$ be the set of matrices of~$\mat{p}{q}$ that are
$\prec$-maximal within their orbits under the group action~$G$.

We can now define the basic objects of this paper.

\begin{definition}[Orbitopes]\
  \begin{myenumerate}
  \item The \emph{full orbitope} associated with the group~$G$ is
    \[
    \orbifullG{p}{q}{G} := \convop\,\matmaxG{p}{q}{G}.
    \]
  \item We associate with the group~$G$ the following restricted orbitopes:
    \begin{align*}
      & \orbipackG{p}{q}{G} :=
      \convop(\matmaxG{p}{q}{G}\cap\matpack{p}{q})
      \quad\text{\emph{(packing orbitope)}}
      \\
      & \orbipartG{p}{q}{G} :=
      \convop(\matmaxG{p}{q}{G}\cap\matpart{p}{q})
      \quad\text{\emph{(partitioning orbitope)}}
      \\
      & \orbicoverG{p}{q}{G} :=
      \convop(\matmaxG{p}{q}{G}\cap\matcov{p}{q})
      \quad\text{\emph{(covering orbitope)}}
    \end{align*}
  \end{myenumerate}
\end{definition}

\begin{remark*}
  By definition, $\orbipartG{p}{q}{G}$ is a face of both
  $\orbipackG{p}{q}{G}$ and $\orbicoverG{p}{q}{G}$.
\end{remark*}

In this paper, we will be only concerned with the cases of~$G$ being the
\emph{cyclic group} $\cyclgr{q}$ containing all~$q$ cyclic permutations of~$\ints{q}$
(Section~\ref{sec:PackingAndPartitioningOrbitopesForCyclicGroups}) or the
symmetric group $\symgr{q}$
(Section~\ref{PackingAndPartitioningOrbitopesForSymmetricGroups}).
Furthermore, we will restrict attention to packing and partitioning
orbitopes. For these, we have the following convenient characterizations of
vertices:

\begin{observation}\label{obs:charvert}\
  \begin{myenumerate}
  \item A matrix of $\mat{p}{q}$ is contained in
    $\matmaxG{p}{q}{\symgr{q}}$ if and only if its columns are in non-increasing
    lexicographic order (with respect to the order $\prec$ defined
    above).%
    \label{obs:charvert:sym}
  \item A matrix of $\matpack{p}{q}$ is contained in
    $\matmaxG{p}{q}{\cyclgr{q}}$ if and only if its first column is
    lexicographically not smaller than the remaining ones (with respect to
    the order $\prec$).%
    \label{obs:charvert:cycPack}
  \item In particular, a matrix of $\matpart{p}{q}$ is contained in
    $\matmaxG{p}{q}{\cyclgr{q}}$ if and only if it has a $1$-entry at
    position $(1,1)$.%
    \label{obs:charvert:cycPart}
  \end{myenumerate}
\end{observation}

% ---------------------------------------------------------------
% Optimizing over Orbitopes
% ---------------------------------------------------------------

\subsection{Optimizing over Orbitopes}
\label{sec:OptimizingOverOrbitopes}

The main aim of this paper is to provide complete descriptions of
$\orbipartG{p}{q}{\symgr{q}}$, $\orbipackG{p}{q}{\symgr{q}}$,
$\orbipartG{p}{q}{\cyclgr{q}}$, and $\orbipackG{p}{q}{\cyclgr{q}}$ by
systems of linear equations and linear inequalities. If these orbitopes
admit ``useful'' linear descriptions then the corresponding linear
optimization problems should be solvable efficiently, due to the
equivalence of optimization and separation, see Gr\"otschel, Lov{\'a}sz,
and Schrijver~\cite{GroLS93}.

We start with the cyclic group operation, since the optimization problem is
particularly easy in this case.

\begin{theorem}\label{thm:cyclpartopt}
  Both the linear optimization problem over
  $\matmaxG{p}{q}{\cyclgr{q}}\cap\matpack{p}{q}$ and over
  $\matmaxG{p}{q}{\cyclgr{q}}\cap\matpart{p}{q}$ can be solved in time
  $\bigo{pq}$.
\end{theorem}

\begin{proof}
  We first give the proof for the packing case.

  For a vector $c \in \Q^{\ints{p}\times\ints{q}}$, we consider the linear
  objective function
  \[
  \scal{c}{x} := \sum_{i=1}^p \sum_{j=1}^q c_{ij}\, x_{ij}.
  \]
  The goal is to find a matrix~$A^\star \in \matmaxG{p}{q}{\cyclgr{q}} \cap
  \matpack{p}{q}$ such that $\scal{c}{A^\star}$ is maximal. Let~$A^{\star}$
  be such a $c$-maximal matrix, and let $a^{\star} \in \{0,1\}^p$ be its
  first column. If $a^{\star} = \zerovec$, then $A^{\star} = \zerovec$ by
  Part~\eqref{obs:charvert:cycPack} of Observation~\ref{obs:charvert}. By
  the same observation it follows that if $a^{\star} \neq \zerovec$ and
  $i^{\star} \in \ints{p}$ is the minimum row-index~$i$ with
  $a^{\star}_{i}=1$, then $A^{\star}$ has only zero entries in its
  first~$i^{\star}$ rows, except for the $1$-entry at position
  $(i^{\star},1)$ (there is at most one $1$-entry in each row).
  Furthermore, each row $i > i^{\star}$ of $A^{\star}$ either has no
  $1$-entry or it has its (unique) $1$-entry at some position where~$c$ is
  maximal in row~$i$.

  Thus, we can compute an optimal solution as follows: (1) For each $i \in
  \ints{p}$ determine a vector $b^i \in \{0,1\}^q$ that is the zero vector
  if~$c$ does not have any positive entries in row~$i$ and otherwise is the
  $j$-th standard unit vector, where~$j \in \ints{q}$ is chosen such that
  $c_{ij} = \max \setdef{c_{i\ell}}{\ell \in \ints{q}}$; set $\sigma_i :=
  0$ in the first case and $\sigma_i := c_{ij}$ in the second.  (2) Compute
  the values $s_p := \sigma_p$ and $s_i := \sigma_i+s_{i+1}$ for all $i =
  p-1, p-2, \dots, 1$.  (3) Determine $i^{\star}$ such that
  $c_{i^{\star},1} + s_{i^{\star}+1}$ is maximal among $\setdef{c_{i,1} +
    s_{i+1}}{i \in \ints{p}}$. (4) If $c_{i^{\star},1} + s_{i^{\star}+1}
  \leq 0$, then~$\zerovec$ is an optimal solution. Otherwise, the matrix
  whose $i$-th row equals~$b^i$ for $i \in \{i^{\star}+1, \dots, p\}$ and
  which is all-zero in the first $i^{\star}$ rows, except for a $1$-entry
  at position $(i^{\star},1)$, is optimal.

  From the description of the algorithm it is easy to see that its running
  time is bounded by $\bigo{pq}$ (in the unit-cost model).

  The partitioning case is then straightforward and even becomes easier due
  to Part~\eqref{obs:charvert:cycPart} of Observation~\ref{obs:charvert}.
\end{proof}

\begin{theorem}\label{thm:sympartopt}
  Both the linear optimization problem over
  $\matmaxG{p}{q}{\symgr{q}}\cap\matpack{p}{q}$ and over
  $\matmaxG{p}{q}{\symgr{q}}\cap\matpart{p}{q}$ can be solved in time
  $\bigo{p^2q}$.
\end{theorem}

\begin{proof}
  We give the proof for the partitioning case, indicating the necessary
  modifications for the packing case at the relevant points.

  As in the proof of Theorem~\ref{thm:cyclpartopt}, we maximize the linear
  objective function given by~$\scal{c}{x}$ for $c \in
  \Q^{\ints{p}\times\ints{q}}$. We describe a two-step approach.

  In the first step, for $i_1, i_2 \in \ints{p}$ with $i_1 \leq i_2$ and $j
  \in \ints{q}$, we let $M(i_1,i_2,j)$ be $c$-maximal among the matrices in
  $\{0,1\}^{\{i_1, i_1+1, \dots, i_2\} \times \ints{j}}$ with exactly (in
  the packing case:\ at most) one $1$-entry in every row. Denote by
  $\mu(i_1,i_2,j)$ the $c$-value of $M(i_1,i_2,j)$, i.e.,
  \[
  \mu(i_1,i_2,j) = \sum_{k=i_1}^{i_2} \sum_{\ell=1}^j c_{k\ell}\,
  M(i_1,i_2,j)_{k\ell}\,.
  \]
  The values~$\mu(i_1,i_2,j)$ can be computed in time $\bigo{p^2q}$ as
  follows. First, we compute all numbers $\lambda(i,j) = \max
  \setdef{c_{i\ell}}{\ell \in \ints{j}}$ (in the packing case:
  $\lambda(i,j) = \max(0, \setdef{c_{i\ell}}{\ell \in \ints{j}})$) for all
  $i \in \ints{p}$ and $j \in \ints{q}$. This can clearly be done in
  $\bigo{pq}$ steps by using the recursions
  $\lambda(i,j)=\max\{\lambda(i,j-1),c_{ij}\}$ for $j\ge 2$. Then, after
  initializing $\mu(i,i,j) = \lambda(i,j)$ for all $i \in \ints{p}$ and $j
  \in \ints{q}$, one computes $\mu(i_1,i_2,j) = \mu(i_1,i_2-1,j) +
  \lambda(i_2,j)$ for all $j \in \ints{q}$, $i_1 = 1, 2, \dots, p$, and
  $i_2 = i_1+1, i_1+2, \dots, q$; see Figure~\ref{fig:sympartopt}.

  In the second step, for $i \in \ints{p}$ and $j \in \ints{q}$, let
  $T(i,j)$ be $c$-maximal among the matrices in $\{0,1\}^{\{i,i+1,\dots,p\}
    \times \ints{q}}$ with exactly (in the packing case:\ at most) one
  $1$-entry in every row and with columns $j, j+1, \dots, q$ being
  in non-increasing
  lexicographic order. Thus, by Part~\eqref{obs:charvert:sym}
  of Observation~\ref{obs:charvert}, $T(1,1)$ is an optimal solution to our
  linear optimization problem.  Denote by $\tau(i,j)$ the $c$-value of
  $T(i,j)$, i.e.,
  \[
  \tau(i,j) = \sum_{k=i}^p \sum_{\ell=1}^q c_{k\ell}\, T(i,j)_{k\ell}.
  \]

  \begin{figure}
    \centering
    \newcommand{\mystyle}[1]{\footnotesize {#1}}
    \psfrag{mus}{\textcolor{white}{\mystyle{$\mu(i_1,i_2-1,j)$}}}
    \psfrag{lambda}{\mystyle{$\lambda(i_2,j)$}}
    \psfrag{i1}{\mystyle{$i_1$}}
    \psfrag{i2}{\mystyle{$i_2$}}
    \psfrag{i}{\mystyle{$i$}}
    \psfrag{j}{\mystyle{$j$}}
    \psfrag{k}{\mystyle{$k$}}
    \psfrag{mu}{\textcolor{white}{\mystyle{$M$}}}
    \psfrag{tau}{\textcolor{white}{\mystyle{$T$}}}
    \includegraphics[height=.17\textheight]{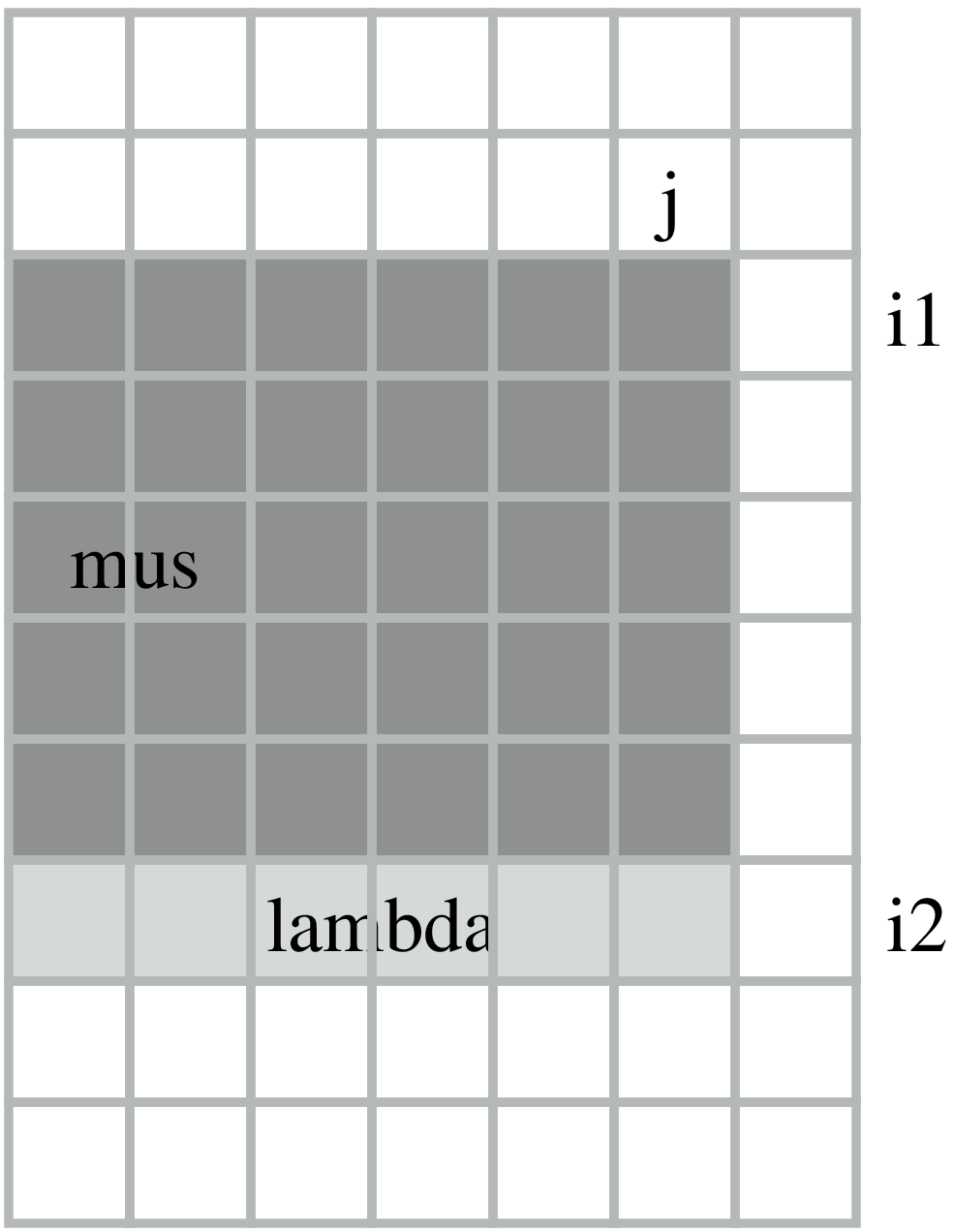}\hspace{10ex}
    \includegraphics[height=.17\textheight]{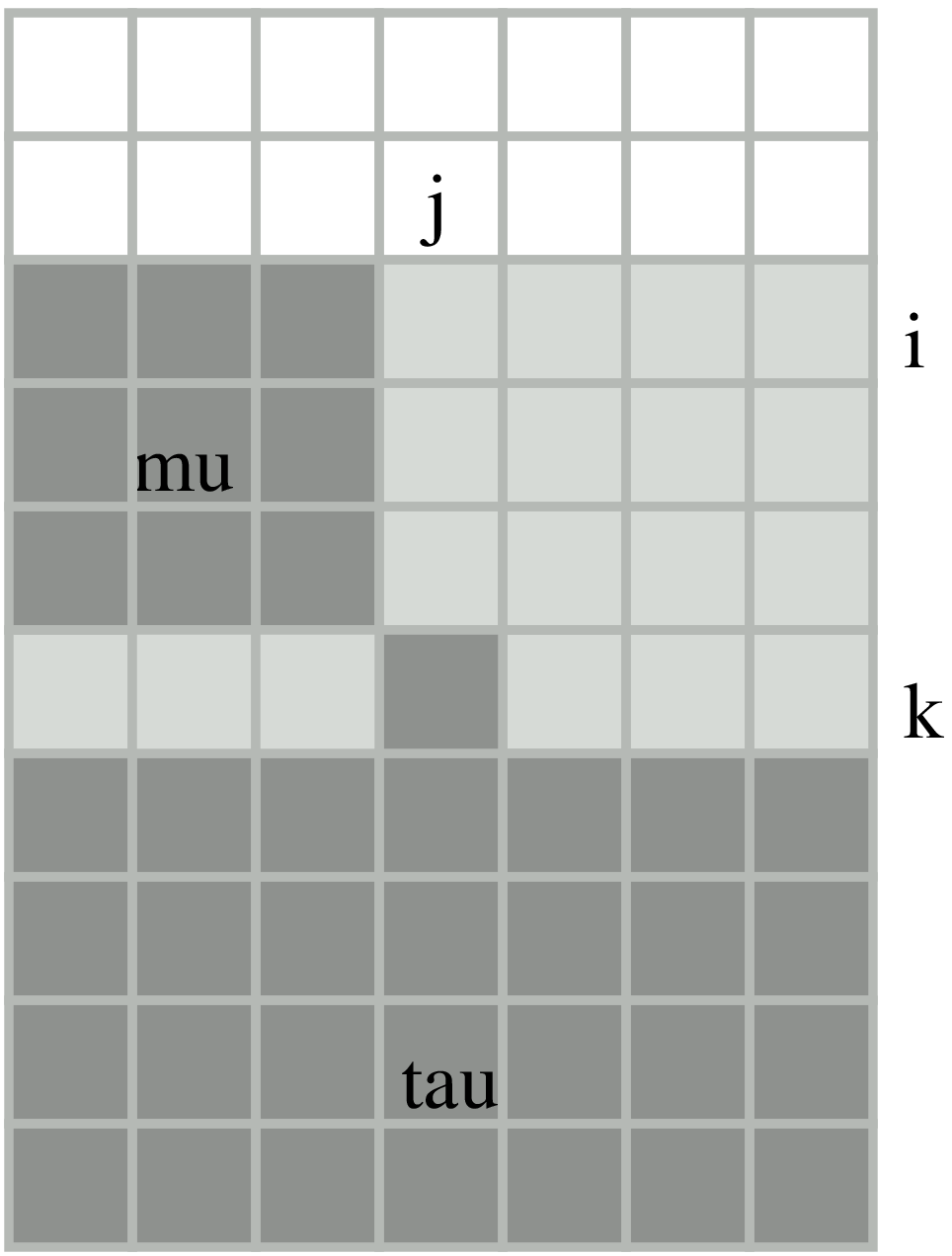}
    \caption{Illustration of the proof of Theorem~\ref{thm:sympartopt}.
      \emph{Left:} Computation of $\mu(i_1,i_2,j)$.  \emph{Right:}
      Computation of~$\tau(i,j)$ via the dynamic programming
      relation~\eqref{eq:sympartopt}. Indicated are the
      matrix~$M(i,k-1,j-1)$ and corresponding term $\mu(i,k-1,j-1)$ and
      matrix $T(k+1,j+1)$ with corresponding term $\tau(k+1,j+1)$.}
    \label{fig:sympartopt}
  \end{figure}

  Let $k\in \{i,i+1,\dots,p+1\}$ be the index of the first row,
  where~$T(i,j)$ has a $1$-entry in column~$j$ (with $k=p+1$ if there is no
  such $1$-entry); see Figure~\ref{fig:sympartopt}. Then $T(i,j)$ has a
  $c$-maximal matrix~$T$ in rows $k+1,\dots,p$ with exactly (in the packing
  case:\ at most) one $1$-entry per row and lexicographically sorted
  columns $j+1,\dots,q$ (contributing $\tau(k+1,j+1)$). In row~$k$, there
  is a single $1$-entry at position~$(k,j)$ (contributing~$c_{kj}$). And in
  rows $i,\dots,k-1$, we have a $c$-maximal matrix~$M$ with exactly (in the
  packing case:\ at most) one $1$-entry per row in the first $j-1$ columns
  (contributing $\mu(i,k-1,j-1)$) and zeroes in the remaining columns.
  Therefore, we obtain
  \[
  \tau(i,j) = \mu(i,k-1,j-1) + c_{kj} + \tau(k+1,j+1).
  \]
  Hence, considering all possibilities for~$k$, we have
  \begin{align}\label{eq:sympartopt}
    \tau(i,j) = \max \, \{ \; & \mu(i,k-1,j-1) + c_{kj} + \tau(k+1,j+1) \; :\\
    & k \in \{i,i+1,\dots,p+1\}\},\notag
  \end{align}
  for all $i \in \ints{p}$ and $j \in \ints{q}$.
  For convenience we define
  $\mu(k_1,k_2,0) = 0$ for $k_1, k_2 \in \ints{p}$ with $k_1 \leq k_2$ and
  $\mu(k,k-1,\ell) = 0$ for all $k \in \ints{p}$ and $\ell \in \{0,1,
  \dots, q\}$.
  Furthermore, we set $c_{p+1,\ell} = 0$ for all $\ell \in \ints{q}$.
  Finally, we define $\tau(p+2,\ell)=\tau(p+1,\ell) = \tau(k,q+1) = 0$ for all $k \in
  \ints{p}$ and $\ell \in \ints{q+1}$.

  Thus, by dynamic programming, we can compute the table $\tau(i,j)$ via
  Equation~\eqref{eq:sympartopt} in the order $i = p, p-1, \dots, 1$, $j =
  q, q-1, \dots, 1$. For each pair $(i,j)$ the evaluation
  of~\eqref{eq:sympartopt} requires no more than $\bigo{p}$ steps, yielding
  a total running time bound of $\bigo{p^2q}$.

  Furthermore, if during these computations for each $(i,j)$ we store a
  maximizer~$k(i,j)$ for~$k$ in~\eqref{eq:sympartopt}, then we can easily
  reconstruct the optimal solution $T(1,1)$ from the $k$-table without
  increasing the running time asymptotically: For $i \in \ints{p}$, $j \in
  \ints{q}$ the matrix $T(i,j)$ is composed of $M(i,k(i,j)-1,j-1)$ (if
  $k(i,j) \geq i+1$ and $j \geq 2$), $T(k(i,j)+1,j+1)$ (if $k(i,j) \leq
  p-1$ and $j \leq q-1$), and having $0$-entries everywhere else, except
  for a $1$-entry at position $(k(i,j),j)$ (if $k(i,j) \leq p$).  Each
  single matrix $M(i_1,i_2,j)$ can be computed in $\bigo{(i_2-i_1)j}$
  steps.  Furthermore, for the matrices $M(i_1,i_2,j)$ needed during the
  recursive reconstruction of $T(1,1)$, the sets $\{i_1,\dots,i_2\} \times
  \ints{j}$ are pairwise disjoint (see Figure~\ref{fig:sympartopt}). Thus,
  these matrices all together can be computed in time $\bigo{pq}$.  At the
  end there might be a single $T(k,q+1)$ to be constructed, which trivially
  can be done in $\bigo{pq}$ steps.
\end{proof}

Thus, with respect to complexity theory there are no ``obstructions'' to
finding complete linear descriptions of packing and partitioning orbitopes
for both the cyclic and the symmetric group action.  In fact, for cyclic
group actions we will provide such a description in
Theorem~\ref{thm:orbipartcycl} and Theorem~\ref{thm:orbipackcycl} for the
partitioning and packing case, respectively. For symmetric group actions we
will provide such a description for partitioning orbitopes in
Theorems~\ref{thm:resultsPart} and for packing orbitopes in
Theorem~\ref{thm:resultsPack}. The algorithm used in the proof of
Theorem~\ref{thm:cyclpartopt} (for cyclic groups) is trivial, while the one
described in the proof of Theorem~\ref{thm:sympartopt} (for symmetric
groups) is a bit more complicated. This is due to the simpler
characterization of the cyclic case in Observation~\ref{obs:charvert} and
is reflected by the fact that the proofs of Theorems~\ref{thm:resultsPart}
and~\ref{thm:resultsPack} (for symmetric groups) need much more work than
the ones of Theorems~\ref{thm:orbipartcycl} and~\ref{thm:orbipackcycl} (for
cyclic groups).

The algorithms described in the above two proofs heavily rely on the fact
that we are considering only matrices with at most one $1$-entry per row.
For cyclic group operations, the case of matrices with more ones per row
becomes more involved, because we do not have a simple characterization
(like the one given in parts~\ref{obs:charvert:cycPack}
and~\ref{obs:charvert:cycPart} of Observation~\ref{obs:charvert}) of the
matrices in $\matmaxG{p}{q}{\cyclgr{q}}$ anymore. For the action of the
symmetric group, though we still have the characterization provided by
Part~\eqref{obs:charvert:sym} of Observation~\ref{obs:charvert}, the
dynamic programming approach used in the proof of
Theorem~\ref{thm:sympartopt} cannot be adapted straight-forwardly without
resulting in an exponentially large dynamic programming table (unless $q$
is fixed). These difficulties apparently are reflected in the structures of
the corresponding orbitopes (see the remarks in
Section~\ref{sec:ClosingRemarks}).

% +++++++++++++++++++++++++++++++++++++++++++++++++++++
% Packing and partitioning orbitopes for cyclic groups
% +++++++++++++++++++++++++++++++++++++++++++++++++++++

\section{Packing and Partitioning Orbitopes for Cyclic Groups}
\label{sec:PackingAndPartitioningOrbitopesForCyclicGroups}

From the characterization of the vertices in
parts~\eqref{obs:charvert:cycPack} and~\eqref{obs:charvert:cycPart} of
Observation~\ref{obs:charvert} one can easily derive IP-formulations of
both the partitioning orbitope $\orbipartG{p}{q}{\cyclgr{q}}$ and the
packing orbitope $\orbipackG{p}{q}{\cyclgr{q}}$ for the cyclic
group~$\cyclgr{q}$. In fact, it turns out that these formulations do
already provide linear descriptions of the two polytopes, i.e., they are
totally unimodular. We refer the reader to Schrijver~\cite[Chap.~19]{Sch86}
for more information on total unimodularity.

It is easy to see that for the descriptions given in
Theorems~\ref{thm:orbipartcycl} and~\ref{thm:orbipackcycl} below, the
separation problem can be solved in time $\bigo{pq}$.

\begin{theorem}\label{thm:orbipartcycl}
  The partitioning orbitope $\orbipartG{p}{q}{\cyclgr{q}}$ for the cyclic
  group~$\cyclgr{q}$ equals the set of all $x \in \R^{\ints{p} \times
    \ints{q}}$ that satisfy the following linear constraints:

  \begin{myitemize}
  \item the equations $x_{11} = 1$ and $x_{1j} = 0$ for all $2 \leq j \leq q$,
  \item the nonnegativity constraints $x_{ij} \geq 0$ for all
    $2 \leq i \leq p$ and $j \in \ints{q}$,
  \item the row-sum equations $x(\row{i}) = 1$ for all $2 \leq i \leq p$.
  \end{myitemize}
  This system of constraints is non-redundant.
\end{theorem}

\begin{proof}
  The constraints $x(\row{i}) = 1$ for $i \in \ints{p}$ and $x_{ij}\ge 0$ for $i
  \in \ints{p},j\in\ints{q}$ define an integral polyhedron, since they
  describe a transshipment problem (and thus, the coefficient matrix is
  totally unimodular).  Hence, the constraint system given in the statement
  of the theorem describes an integer polyhedron, because it defines a face
  of the corresponding transshipment polytope.
  % Let~$A$ be the coefficient matrix of the constraint system in the
  % theorem. Thus, $A$ has $pq$ columns (and $(p+1)q-1$ rows), which we
  % assume to be ordered according to the lexicographic order~$\prec$ of
  % $\ints{p} \times \ints{q}$. Let~$\tilde{A}$ be the submatrix of~$A$
  % belonging to the~$p$ row-sum equations. Because in order~$\prec$ the
  % entries of each row appear consecutively, $\tilde{A}$ is a
  % 0/1-interval-matrix (i.e., if there are two $1$-entries in the same
  % row, then all entries between the two $1$-entries are $1$'s as well),
  % and hence $\tilde{A}$ is totally unimodular.  Since~$A$ arises
  % from~$\tilde{A}$ by adding an identity matrix, $A$ is totally
  % unimodular, too. In fact, the constraint system in the statement of the
  % Theorem defines a face of a transshipment polytope.

  By Part~\eqref{obs:charvert:cycPart} of Observation~\ref{obs:charvert},
  the set of integer points satisfying this constraint system is
  $\matpart{p}{q}\cap\matmaxG{p}{q}{\cyclgr{q}}$.  Hence the given
  constraints completely describe $\orbipartG{p}{q}{\cyclgr{q}}$. The
  non-redundancy follows from the fact that dropping any of the constraints
  enlarges the set of feasible integer solutions.
\end{proof}

\begin{figure}
  \centering
  \newcommand{\mystyle}[1]{\footnotesize {#1}}
  \psfrag{i}{\mystyle{$i$}}
  \includegraphics[height=.17\textheight]{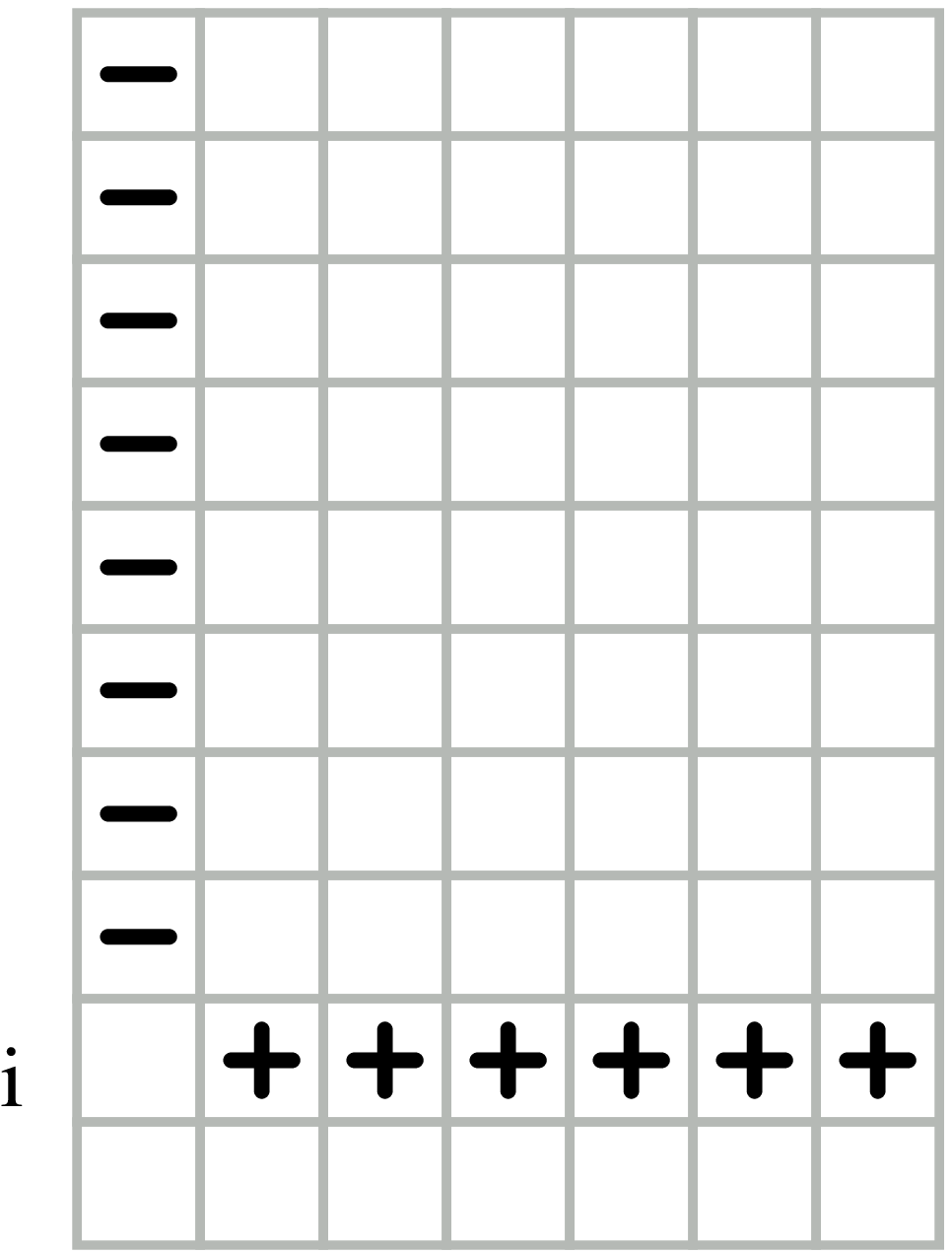}\label{sfl:packcyclcoline}
  \caption[]{Example of the coefficient
    vector for an inequality of type~\eqref{eq:orbipackcycl:1}; ``$-$'' stands for a
    $-1$, ``$+$'' for a $+1$.}
  \label{fig:orbipackcycl}
\end{figure}

\begin{theorem}\label{thm:orbipackcycl}
  The packing orbitope $\orbipackG{p}{q}{\cyclgr{q}}$ for the cyclic
  group~$\cyclgr{q}$ equals the set of all $x \in \R^{\ints{p} \times
    \ints{q}}$ that satisfy the following linear constraints:
  \begin{myitemize}
  \item the constraints $0\le x_{11} \leq 1$ and $x_{1j} = 0$ for all $2 \leq
    j \leq q$,
  \item the nonnegativity constraints $x_{ij} \geq 0$ for all $2 \leq i
    \leq p$ and $j \in \ints{q}$,
  \item the row-sum inequalities $x(\row{i}) \leq 1$ for all $2 \leq i \leq
    p$,
  \item the inequalities
    \begin{equation}\label{eq:orbipackcycl:1}
      \sum_{j=2}^q x_{ij} - \sum_{k=1}^{i-1} x_{k1} \leq 0
    \end{equation}
    for all $2 \leq i \leq p$ (see
    Figure~\ref{fig:orbipackcycl} for an example).
  \end{myitemize}
  This system of constraints is non-redundant.
\end{theorem}

\begin{figure}
  \centering
  \newcommand{\mystyle}[1]{\tiny {#1}}
  % paths
  \psfrag{P1}{\mystyle{$P_1$}}
  \psfrag{P2}{\mystyle{$P_2$}}
  \psfrag{Pi-1}{\mystyle{$P_{i-1}$}}
  \psfrag{Pi}{\mystyle{$P_i$}}
  \psfrag{Pi+1}{\mystyle{$P_{i+1}$}}
  \psfrag{Pp}{\mystyle{$P_p$}}
  % first row:
  \psfrag{v11}{\mystyle{$v_{11}$}}
  \psfrag{v21}{\mystyle{$v_{21}$}}
  \psfrag{vi-11}{\mystyle{$v_{i-1,1}$}}
  \psfrag{vi1}{\mystyle{$v_{i1}$}}
  \psfrag{vi+11}{\mystyle{$v_{i+1,1}$}}
  \psfrag{vp1}{\mystyle{$v_{p1}$}}
  \psfrag{vp+11}{\mystyle{$v_{p+1,1}$}}
  % vertical middle
  \psfrag{vij-1}{\mystyle{$v_{i,j-1}$}}
  \psfrag{vij}{\mystyle{$v_{ij}$}}
  \psfrag{viq}{\mystyle{$v_{iq}$}}
  % horizontal arcs
  \psfrag{a11}{\mystyle{$\alpha_{11}$}}
  \psfrag{ai-11}{\mystyle{$\alpha_{i-1,1}$}}
  \psfrag{ai1}{\mystyle{$\alpha_{i1}$}}
  \psfrag{ap1}{\mystyle{$\alpha_{p1}$}}
  % vertical arcs: 2nd column
  \psfrag{a22}{\mystyle{$\alpha_{22}$}}
  \psfrag{a2j}{\mystyle{$\alpha_{2j}$}}
  \psfrag{a2q}{\mystyle{$\alpha_{2q}$}}
  % vertical arcs: ith column
  \psfrag{ai2}{\mystyle{$\alpha_{i2}$}}
  \psfrag{aij}{\mystyle{$\alpha_{ij}$}}
  \psfrag{aiq}{\mystyle{$\alpha_{iq}$}}
  % vertical arcs: pth column
  \psfrag{ap2}{\mystyle{$\alpha_{p2}$}}
  \psfrag{apj}{\mystyle{$\alpha_{pj}$}}
  \psfrag{apq}{\mystyle{$\alpha_{pq}$}}
  \includegraphics[height=.3\textheight]{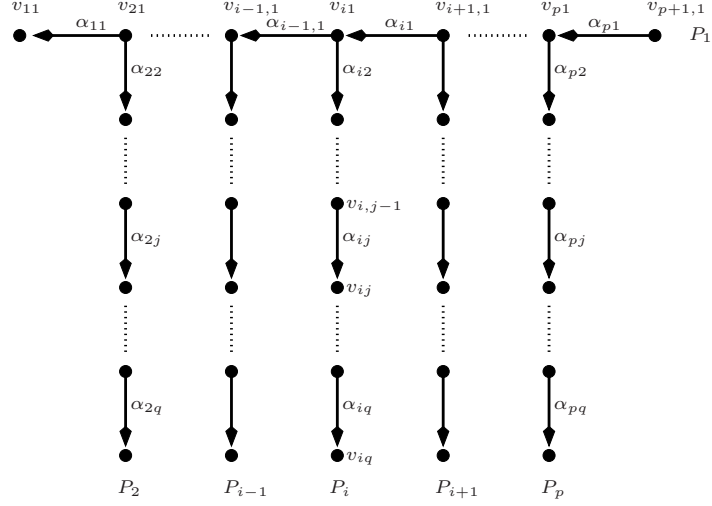}\label{sfl:networkmatrix}
  \caption[]{The network matrix constructed in the proof of Theorem~\ref{thm:orbipackcycl}.}
  \label{fig:networkmatrix}
\end{figure}

\begin{proof}
  From Part~\eqref{obs:charvert:cycPack} of Observation~\ref{obs:charvert}
  it follows that an integer point is contained in
  $\orbipackG{p}{q}{\cyclgr{q}}$ if and only if it satisfies the
  constraints described in the statement, where
  Inequalities~\eqref{eq:orbipackcycl:1} ensure that the first column
  of~$x$ is lexicographically not smaller than the other ones (note that we
  have at most one $1$-entry in each row of~$x$). Dropping any of the
  constraints enlarges the set of integer solutions, which proves the
  statement on non-redundancy.  Thus, as in the proof of the previous
  theorem, it remains to show that the polyhedron defined by the
  constraints is integral. We prove this by showing that the coefficient
  matrix~$A$ of the row-sum inequalities $x(\row{i}) \leq 1$ (for $2 \leq i
  \leq p$) and Inequalities~\eqref{eq:orbipackcycl:1} (for all $2 \leq i
  \leq p$) is a network matrix (and thus, totally unimodular). Adding the
  nonnegativity constraints amounts to adding an identity matrix and preserves
  total unimodularity, which also holds for the inclusion of $x_{11}\le 1$ into the system.

  In order to establish the claim on the network structure of~$A$, we will
  identify a directed tree~$T$, whose arcs are in bijection with $\ints{p}
  \times \ints{q}$ (the set of indices of the columns of~$A$), such that
  there are pairs of nodes $(v_r, w_r)$ of~$T$ in bijection with the row
  indices $r \in \ints{2(p-1)}$ of~$A$ with the following property. The
  matrix~$A$ has a $(+1)$-entry in row~$r$ and column~$(i,j)$, if the
  unique path~$\pi_r$ from node~$v_r$ to node~$w_r$ in the tree~$T$ uses
  arc~$(i,j)$ in its direction from~$i$ to~$j$, a $(-1)$-entry, if~$\pi_r$
  uses $(i,j)$ in its reverse direction, and a $0$-entry, if~$\pi_r$ does
  not use~$(i,j)$.

  For the construction of the tree~$T$, we take a directed path~$P_1$ of
  length~$p$ on nodes $\{v_{11}, v_{21}, \dots, v_{p+1,1}\}$ with arcs
  $\alpha_{i1} := (v_{i+1,1}, v_{i1})$ for $i \in \ints{p}$; see
  Figure~\ref{fig:networkmatrix}.  For each $2 \leq i \leq p$, we append a
  directed path~$P_i$ of length~$q-1$ to node~$v_{i1}$, where~$P_i$ has
  node set $\{v_{i1}, v_{i2}, \dots, v_{iq}\}$ and arcs $\alpha_{ij} :=
  (v_{i,j-1}, v_{ij})$ for $2 \leq j \leq q$. Choosing the pair
  $(v_{i+1,1},v_{iq})$ for the $i$-th row sum-inequality and the pair
  $(v_{11},v_{iq})$ for the $i$-th Inequality~\eqref{eq:orbipackcycl:1},
  finishes the proof (using the bijection between the arcs of~$T$ and the
  columns of~$A$ indicated by the notation $\alpha_{ij}$).
\end{proof}

% +++++++++++++++++++++++++++++++++++++++++++++++++++++++++++++++
% Packing and partitioning orbitopes for symmetric groups
% +++++++++++++++++++++++++++++++++++++++++++++++++++++++++++++++

\section{Packing and Partitioning Orbitopes for Symmetric Groups}
\label{PackingAndPartitioningOrbitopesForSymmetricGroups}

For packing orbitopes~$\orbipack{p}{q}$ and partitioning
orbitopes~$\orbipart{p}{q}$ with respect to the symmetric group it follows
readily from the characterization in Part~\eqref{obs:charvert:sym} of
Observation~\ref{obs:charvert} that the equations
\begin{equation}\label{eq:fixtriang}
  x_{ij} = 0\qquad\text{ for all } i < j
\end{equation}
are valid. Thus, we may drop all variables corresponding to components in
the upper right triangle from the formulation and consider
\[
\orbipack{p}{q},\; \orbipart{p}{q} \subset \R^{\orbipartinds{p}{q}}
\qquad\text{with}\quad \orbipartinds{p}{q} := \setdef{(i,j) \in \ints{p}
  \times \ints{q}}{i \geq j}.
\]
We also adjust the definition of
\[
\row{i} := \{(i,1), (i,2), \dots, (i,\min\{i,q\})\} \qquad \text{ for } i
\in \ints{p}
\]
and define the $j$th column for $j \in \ints{q}$ as
\[
\col{j} := \{(j,j), (j+1,j), \dots, (p,j)\}.
\]
Furthermore, we restrict ourselves to the case
\[
p \geq q \geq 2
\]
in this context. Because of~\eqref{eq:fixtriang}, the case of $q > p$ can
be reduced to the case $p = q$ and the case of $q = 1$ is of no interest.

The next result shows a very close relationship between packing and
partitioning orbitopes for the case of symmetric group actions.

\begin{proposition}\label{prop:projpartpack}
  The polytopes $\orbipart{p}{q}$ and $\orbipack{p-1}{q-1}$ are affinely
  isomorphic via orthogonal projection of $\orbipart{p}{q}$ onto the
  space
  \[
  \mathcal{L} :=
  \setdef{x \in \R^{\orbipartinds{p}{q}}}{x_{i1} = 0 \text{ for all }i \in \ints{p}}
  \]
  $($and the canonical identification of this space with
  $\R^{\orbipartinds{p-1}{q-1}})$.
\end{proposition}

\begin{proof}
  The affine subspace
  \[
  \mathcal{A} :=
  \setdef{x \in \R^{\orbipartinds{p}{q}}}{ x(\row{i}) = 1 \text{ for all }i}
  \]
  of $\R^{\orbipartinds{p}{q}}$ clearly contains $\orbipart{p}{q}$. Let
  $\pi : \mathcal{A} \rightarrow \R^{\orbipartinds{p-1}{q-1}}$ be the
  orthogonal projection mentioned in the statement (identifying
  $\mathcal{L}$ in the canonical way with $\R^{\orbipartinds{p-1}{q-1}}$);
  note that the first row is removed since it only contains the
  element~$(1,1)$.  Consider the linear map $\phi:
  \R^{\orbipartinds{p-1}{q-1}} \rightarrow \R^{\orbipartinds{p}{q}}$
  defined by
  \[
  \phi(y)_{ij} =
  \begin{cases}
    1 - y(\row{i-1}) & \text{if }j=1\\
    y_{i-1,j-1}    & \text{otherwise}
  \end{cases}
  \qquad\text{for }(i,j) \in \orbipartinds{p}{q}
  \]
  (where $\row{0} = \varnothing$ and $y(\varnothing) = 0$). This is the
  inverse of~$\pi$, showing that $\pi$ is an affine isomorphism. As we have
  $\pi(\orbipart{p}{q})=\orbipack{p-1}{q-1}$, this finishes the proof.
\end{proof}

It will be convenient to address the elements in $\orbipartinds{p}{q}$
via a different ``system of coordinates'':
\[
\diagcol{\eta}{j} := (j+\eta-1,j) \qquad\text{for }j \in \ints{q},\;
1 \leq \eta \leq p-j+1.
\]
Thus (as before) $i$ and~$j$ denote the row and the columns, respectively,
while~$\eta$ is the index of the diagonal (counted from above) containing
the respective element; see Figure~\ref{fig:scispart}~\subref{sfl:coord}
for an example. For $(k,j) = \diagcol{\eta}{j}$ and $x \in
\R^{\orbipartinds{p}{q}}$, we write $x_{\diagcol{\eta}{j}} := x_{(k,j)} :=
x_{kj}$.

For $x \in \{0,1\}^{\orbipartinds{p}{q}}$ we denote by $\setof{x} :=
\setdef{(i,j) \in \orbipartinds{p}{q}}{x_{ij} = 1}$ the set of all
coordinates (positions in the matrix), where~$x$ has a $1$-entry.
Conversely, for $I \subseteq \orbipartinds{p}{q}$, we use $\pointof{I} \in
\{0,1\}^{\orbipartinds{p}{q}}$ for the 0/1-point with $\pointof{I}_{ij} =
1$ if and only if $(i,j) \in I$.

For $(i,j)\in\orbipartinds{p}{q}$, we define the \emph{column}
\[
\colop(i,j) = \{(j,j), (j+1,j), \dots, (i-1,j), (i,j)\} \subseteq \orbipartinds{p}{q},
\]
and for $(i,j) = \diagcol{\eta}{j}$ we write $\colop\diagcol{\eta}{j} :=
\colop(i,j)$. Of course, we have $\colop\diagcol{\eta}{j} =
\{\diagcol{1}{j}, \diagcol{2}{j}, \dots, \diagcol{\eta}{j}\}$.

The rest of this section is organized as follows. First, in
Section~\ref{sec:IPFormulations}, we deal with basic facts about integer
points in packing and partitioning orbitopes for the symmetric group.  To
derive a linear description of $\orbipack{p}{q}$ and $\orbipart{p}{q}$ that
only contains integer vertices, we need additional inequalities, the
\emph{shifted column inequalities}, which are introduced in
Section~\ref{sec:ShiftedColumnInequalities}. We then show that the
corresponding separation problem can be solved in linear time
(Section~\ref{sec:SeparationAlgorithmForSCIs}).
Section~\ref{sec:CompleteInequalityDescriptions} proves the completeness of
the linear description and Section~\ref{sec:Facets} investigates the facets
of the polytopes.

% ---------------------------------------------------------------
% Characterization of Integer Vertices
% ---------------------------------------------------------------

\subsection{Characterization of Integer Points}
\label{sec:IPFormulations}

We first derive a crucial property of the vertices of $\orbipack{p}{q}$.

\begin{lemma}\label{lem:xcapcol}
  Let~$x$ be a vertex of $\orbipack{p}{q}$ with $\diagcol{\eta}{j} \in
  \setof{x}$ $(j \geq 2)$. Then we have $\setof{x} \cap
  \colop\diagcol{\eta}{j-1} \neq \varnothing$.
\end{lemma}

\begin{proof}
  With $\diagcol{\eta}{j} = (i,j)$ we have $x_{ij} = 1$, which implies
  $x_{i,j-1} = 0$ (since~$x$ has at most one $1$-entry in row~$i$). Thus,
  $\setof{x} \cap \colop\diagcol{\eta}{j-1} = \varnothing$ would yield
  $x_{k,j-1} = 0$ for all $k \leq i$, contradicting the lexicographic order
  of the columns of~$x$ (see Part~\eqref{obs:charvert:sym} of
  Observation~\ref{obs:charvert}).
\end{proof}

\begin{definition}[Column inequality] For $(i,j) \in \orbipartinds{p}{q}$
  and the set $B = \{(i,j)$, $(i,j+1), \dots, (i,\min\{i,q\})\}$, we call
  \[
  x(B)-x(\colop(i-1,j-1)) \leq 0
  \]
  a \emph{column inequality}; see
  Figure~\ref{fig:scispart}~\subref{sfl:scipart1} for an example with
  $(i,j) = (9,5)$.
\end{definition}

The column inequalities are strengthenings of the symmetry breaking inequalities
\begin{equation}\label{eq:shortColInequality}
  x_{ij} - x(\colop(i-1,j-1)) \leq 0,
\end{equation}
introduced by~M\'{e}ndez-D\'{i}az and Zabala~\cite{DiaZ01} in the context
of vertex-coloring (see \eqref{eq:SymmetryBreak} in the introduction).

\begin{proposition}\label{prop:IPsym}
  A point $x \in \{0,1\}^{\orbipartinds{p}{q}}$ is contained in
  $\orbipack{p}{q}$ $(\orbipart{p}{q})$ if and only if~$x$ satisfies the
  row-sum constraints $x(\rowop(i)) \leq 1$ $(x(\rowop(i)) = 1)$ for all $i
  \in \ints{p}$ and all column inequalities.
\end{proposition}

\begin{proof}
  By Lemma~\ref{lem:xcapcol}, Inequalities~\eqref{eq:shortColInequality}
  are valid for $\orbipack{p}{q}$ (and thus, for its face $\orbipart{p}{q}$
  as well).  Because of the row-sum constraints, all column inequalities
  are valid as well. Therefore, it suffices to show that a point $x \in
  \{0,1\}^{\orbipartinds{p}{q}}$ that satisfies the row-sum constraints
  $x(\rowop(i)) \leq1$ and all column inequalities is contained in
  $\matmaxG{p}{q}{\symgr{q}}$.

  Suppose, this was not the case. Then, by Part~\eqref{obs:charvert:sym} of
  Observation~\ref{obs:charvert}, there must be some~$j \in \ints{q}$ such
  that the $(j-1)$-st column of~$x$ is lexicographically smaller than the
  $j$th column. Let~$i$ be minimal with $x_{ij} = 1$ (note that column~$j$
  cannot be all-zero). Thus, $x_{k,j-1} = 0$ for all $k < i$. This implies
  $x(\colop(i-1,j-1)) = 0 < 1 = x_{ij}$, showing that the column inequality
  $x(B) - x(\colop(i-1,j-1)) \leq 0$ is violated by the point~$x$ for the
  bar $B = \{(i,j),(i,j+1), \dots, (i,\min\{i,q\})\}$.
\end{proof}

% ---------------------------------------------------------------
% Shifted Column Inequalities
% ---------------------------------------------------------------

\subsection{Shifted Column Inequalities}
\label{sec:ShiftedColumnInequalities}

Proposition~\ref{prop:IPsym} provides a characterization of the vertices of
the packing- and partitioning orbitopes for symmetric groups among the
integer points. Different from the situation for cyclic groups (see
Theorems~\ref{thm:orbipartcycl} and~\ref{thm:orbipackcycl}), however, the
inequalities in this characterization do not yield complete descriptions of
these orbitopes. In fact, we need to generalize the concept of a column
inequality in order to arrive at complete descriptions. This will yield
exponentially many additional facets (see
Proposition~\ref{prop:orbipack:facets}).

\begin{definition}[Shifted columns]
  A set $S = \{\diagcol{1}{c_1}, \diagcol{2}{c_2}, \dots,
  \diagcol{\eta}{c_{\eta}}\} \subset \orbipartinds{p}{q}$ with $\eta \geq
  1$ and $c_1 \leq c_2 \leq \dots \leq c_{\eta}$ is called a
  \emph{shifted column}. It is a \emph{shifting} of each of the columns
  \[
  \colop\diagcol{\eta}{c_{\eta}}, \colop\diagcol{\eta}{c_{\eta}+1}, \dots,\colop\diagcol{\eta}{q}.
  \]
\end{definition}

\begin{remark*}\

  \vspace{-1ex}% ugly, but seems to be necessary (cannot set topsep to 0ex here)
  \begin{myitemize}
  \item As a special case we have column $\colop(i,j)$, which is the
    shifted column $\{\diagcol{1}{j}, \diagcol{2}{j}, \dots,
    \diagcol{\eta}{j}\}$ for $\diagcol{\eta}{j} = (i,j)$.
  \item By definition, if $S = \{\diagcol{1}{c_1}, \diagcol{2}{c_2}, \dots,
    \diagcol{\eta}{c_{\eta}}\} \subset \orbipartinds{p}{q}$ is a shifted
    column, then so is $\{\diagcol{1}{c_1}, \diagcol{2}{c_2}, \dots,
    \diagcol{\eta'}{c_{\eta'}}\}$ for every $1 \leq \eta' \leq \eta$.
  \end{myitemize}
\end{remark*}

\begin{lemma}\label{lem:xcapsc}
  Let $x$ be a vertex of $\orbipack{p}{q}$ with $\diagcol{\eta}{j} \in
  \setof{x}$ $(j \geq 2)$. Then we have $\setof{x} \cap S \neq \varnothing$
  for all shiftings~$S$ of $\colop\diagcol{\eta}{j-1}$.
\end{lemma}

\begin{proof}
  We proceed by induction on~$j$. The case $j = 2$ follows from
  Lemma~\ref{lem:xcapcol}, because the only shifting of
  $\colop\diagcol{\eta}{1}$ is $\colop\diagcol{\eta}{1}$ itself. Therefore,
  let $j \geq 3$, and let $S = \{\diagcol{1}{c_1}, \diagcol{2}{c_2}, \dots,
  \diagcol{\eta}{c_{\eta}}\}$ be a shifting of $\colop\diagcol{\eta}{j-1}$
  (hence, $c_1 \leq c_2 \leq \dots \leq c_{\eta} \leq j-1$). Since by
  assumption $\diagcol{\eta}{j} \in \setof{x}$, Lemma~\ref{lem:xcapcol}
  yields that there is some $\eta' \leq \eta$ with $\diagcol{\eta'}{j-1}
  \in \setof{x}$. If $\diagcol{\eta'}{j-1}\in S$, then we are done.
  Otherwise, $c_{\eta'} < j-1$ holds. Hence, $\{\diagcol{1}{c_1},
  \diagcol{2}{c_2}, \dots, \diagcol{\eta'}{c_{\eta'}}\}$ is a shifting of
  ($\colop\diagcol{\eta'}{c_{\eta'}}$ and hence of)
  $\colop\diagcol{\eta'}{j-2}$, which, by the inductive hypothesis, must
  intersect~$\setof{x}$.
\end{proof}

\begin{figure}
  \centering
  \newcommand{\mystyle}[1]{\footnotesize {#1}}
  \psfrag{i}{\mystyle{$i$}}
  \psfrag{j}{\mystyle{$j$}}
  \psfrag{eta}{\mystyle{$\eta$}}
  \subfloat[][]{\includegraphics[height=.17\textheight]{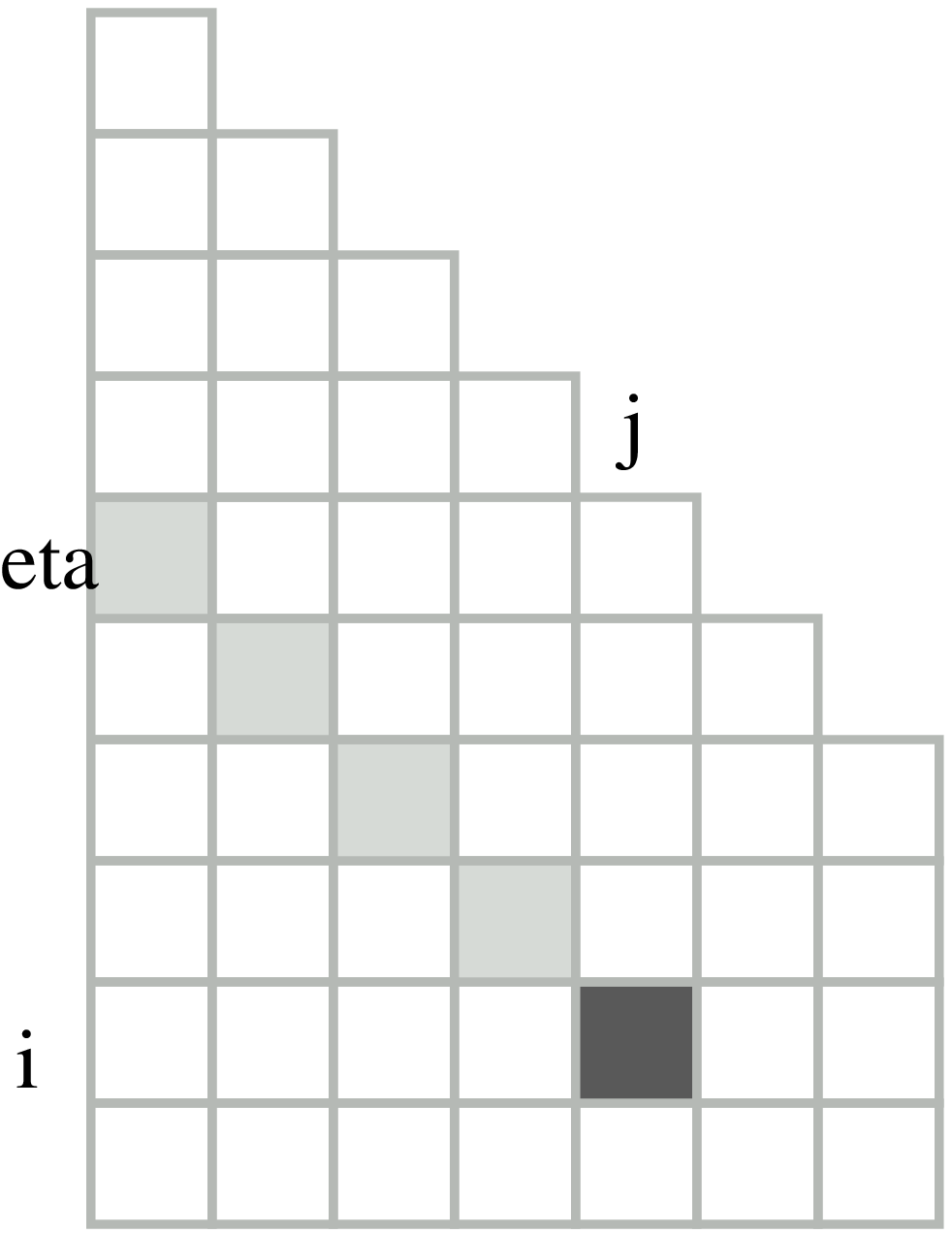}\label{sfl:coord}}\hfill
  \subfloat[][]{\includegraphics[height=.17\textheight]{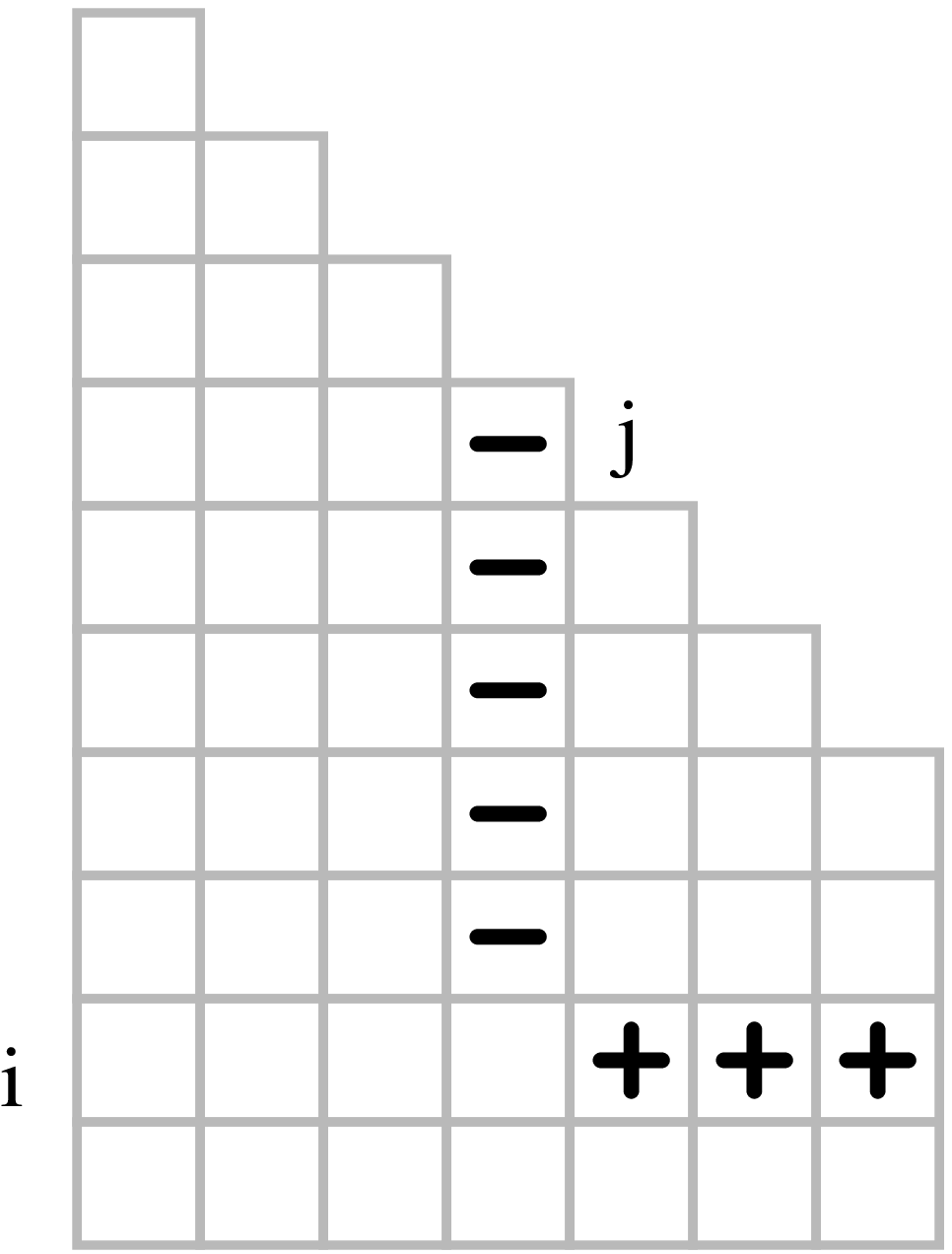}\label{sfl:scipart1}}\hfill
  \subfloat[][]{\includegraphics[height=.17\textheight]{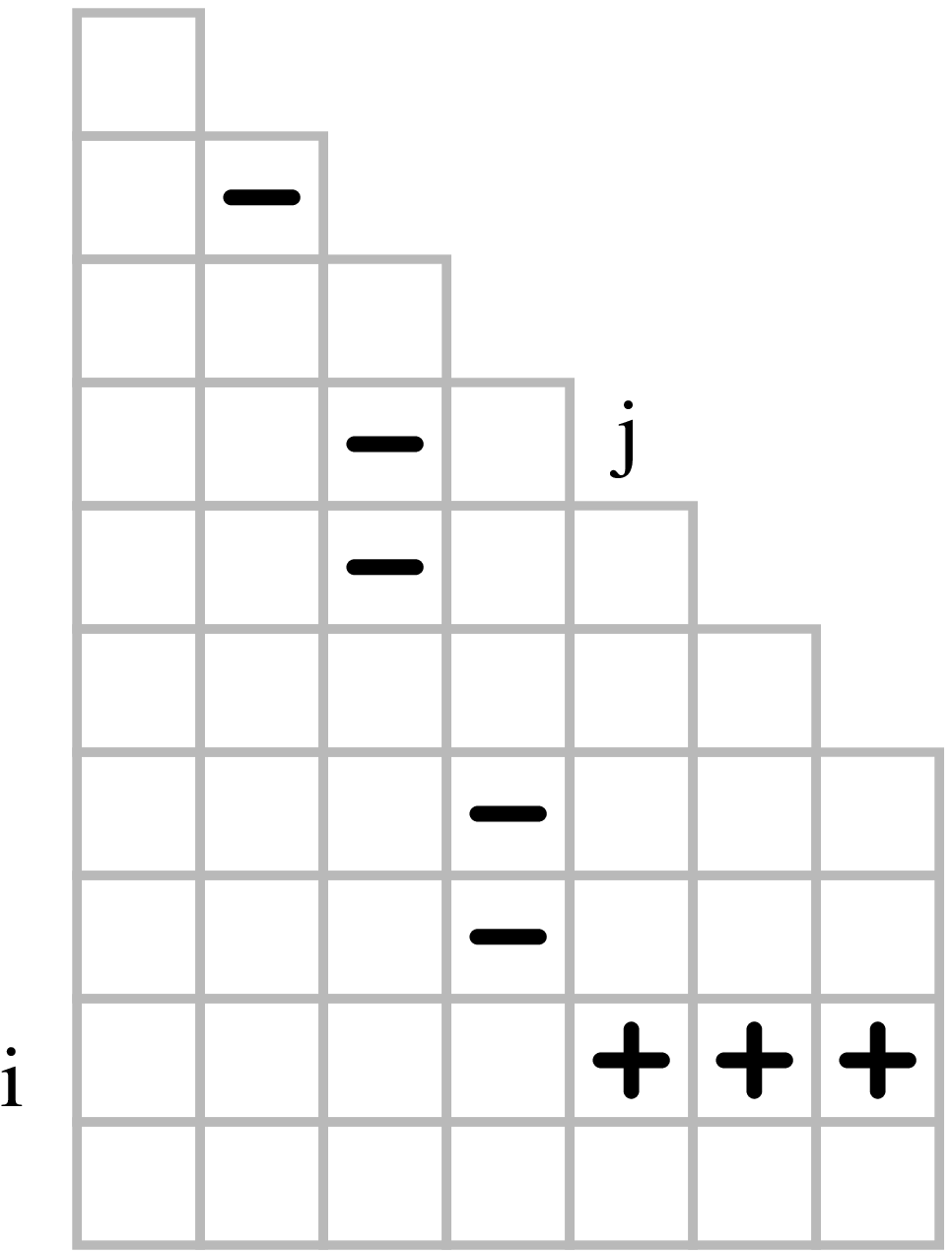}\label{sfl:scipart2}}\hfill
  \subfloat[][]{\includegraphics[height=.17\textheight]{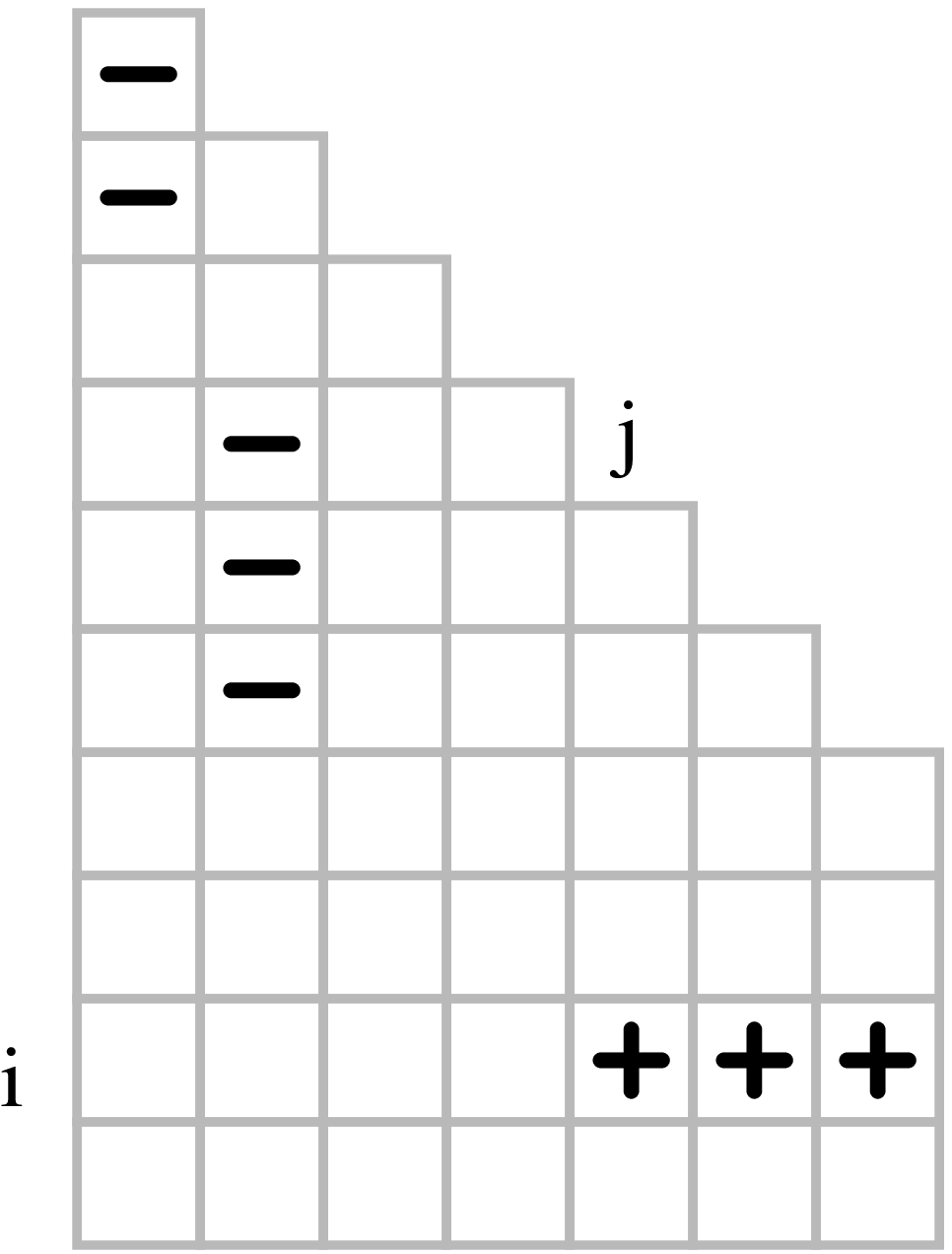}\label{sfl:scipart3}}
  \caption[]{\subref{sfl:coord} Example for coordinates $(9,5) =
    \diagcol{5}{5}$. \subref{sfl:scipart1}--\subref{sfl:scipart3} Shifted
    column inequalities with leader $\diagcol{5}{5}$, see
    Definition~\ref{def:SCIs}.  All SCI inequalities are
    $\leq$-inequalities with right-hand sides zero and ``$-$'' stands for a
    $(-1)$-coefficient, ``$+$`` for a $(+1)$ coefficient. The shifted
    column of~\subref{sfl:scipart2} is $\{\diagcol{1}{2}, \diagcol{2}{3},
    \diagcol{3}{3}, \diagcol{4}{4}, \diagcol{5}{4}\}$.}
  \label{fig:scispart}
\end{figure}

\begin{definition}[Shifted column inequalities]\label{def:SCIs}
  For $(i,j) = \diagcol{\eta}{j} \in \orbipartinds{p}{q}$, $B = \{(i,j),
  (i,j+1), \dots, (i,\min\{i,q\})\}$, and a shifting~$S$ of
  $\colop\diagcol{\eta}{j-1}$, we call
  \[
  x(B) - x(S) \leq 0
  \]
  a \emph{shifted column inequality} \emph{(SCI)}. The set~$B$ is
  the \emph{bar} of the SCI, and $(i,j)$ is the \emph{leader} of (the bar
  of) the SCI. The set $S$ is the \emph{shifted column (SC)} of the
  SCI. See Figure~\ref{fig:scispart} for examples.
\end{definition}

In particular, all column inequalities are shifted column inequalities.
The class of shifted column inequalities, however, is substantially richer:
It contains exponentially many inequalities (in~$q$).

\begin{proposition}\label{prop:SCIvalid}
  Shifted column inequalities are valid both for the packing orbitopes
  $\orbipack{p}{q}$ and for the partitioning orbitopes $\orbipart{p}{q}$.
\end{proposition}

\begin{proof}
  As $\orbipart{p}{q}$ is a face of $\orbipack{p}{q}$, it is enough to
  prove the proposition for packing orbitopes $\orbipack{p}{q}$. Therefore,
  let $(i,j) = \diagcol{\eta}{j} \in \orbipartinds{p}{q}$, with $j \geq 2$,
  and let $S= \{\diagcol{1}{c_1}, \diagcol{2}{c_2}, \dots,
  \diagcol{\eta}{c_{\eta}}\}$ be a shifting of $\colop\diagcol{\eta}{j-1}$.
  Denote by~$B$ the bar of the corresponding SCI.

  Let $x \in \{0,1\}^{\orbipartinds{p}{q}}$ be a vertex of
  $\orbipack{p}{q}$.  If $B \cap \setof{x} = \varnothing$, then clearly
  $x(B) - x(S) = 0 - x(S) \leq 0$ holds. Otherwise, there is a unique
  element $(i,j') = \diagcol{\eta'}{j'} \in B \cap \setof{x}$. As $j' \geq
  j$, we have $\eta' \leq \eta$.  Therefore $S' = \{\diagcol{1}{c_1},
  \diagcol{2}{c_2}, \dots, \diagcol{\eta'}{c_{\eta'}}\} \subseteq S$ is a
  shifting of $\colop\diagcol{\eta'}{j'-1}$. Thus, by
  Lemma~\ref{lem:xcapsc}, we have $S' \cap \setof{x} \neq \varnothing$.
  This shows $x(S) \geq x(S') \geq 1$, implying $x(B) - x(S) \leq 1 - 1 =
  0$.
\end{proof}

% ---------------------------------------------------------------
% A linear time separation algorithm for SCIs
% ---------------------------------------------------------------

\subsection{A Linear Time Separation Algorithm for SCIs}
\label{sec:SeparationAlgorithmForSCIs}

\begin{figure}
  \centering
  \newcommand{\mystyle}[1]{\footnotesize {#1}}
  \psfrag{i}{\mystyle{$i$}}
  \psfrag{j}{\mystyle{$j$}}
  \psfrag{eta}{\mystyle{$\eta$}}
  \includegraphics[height=.17\textheight]{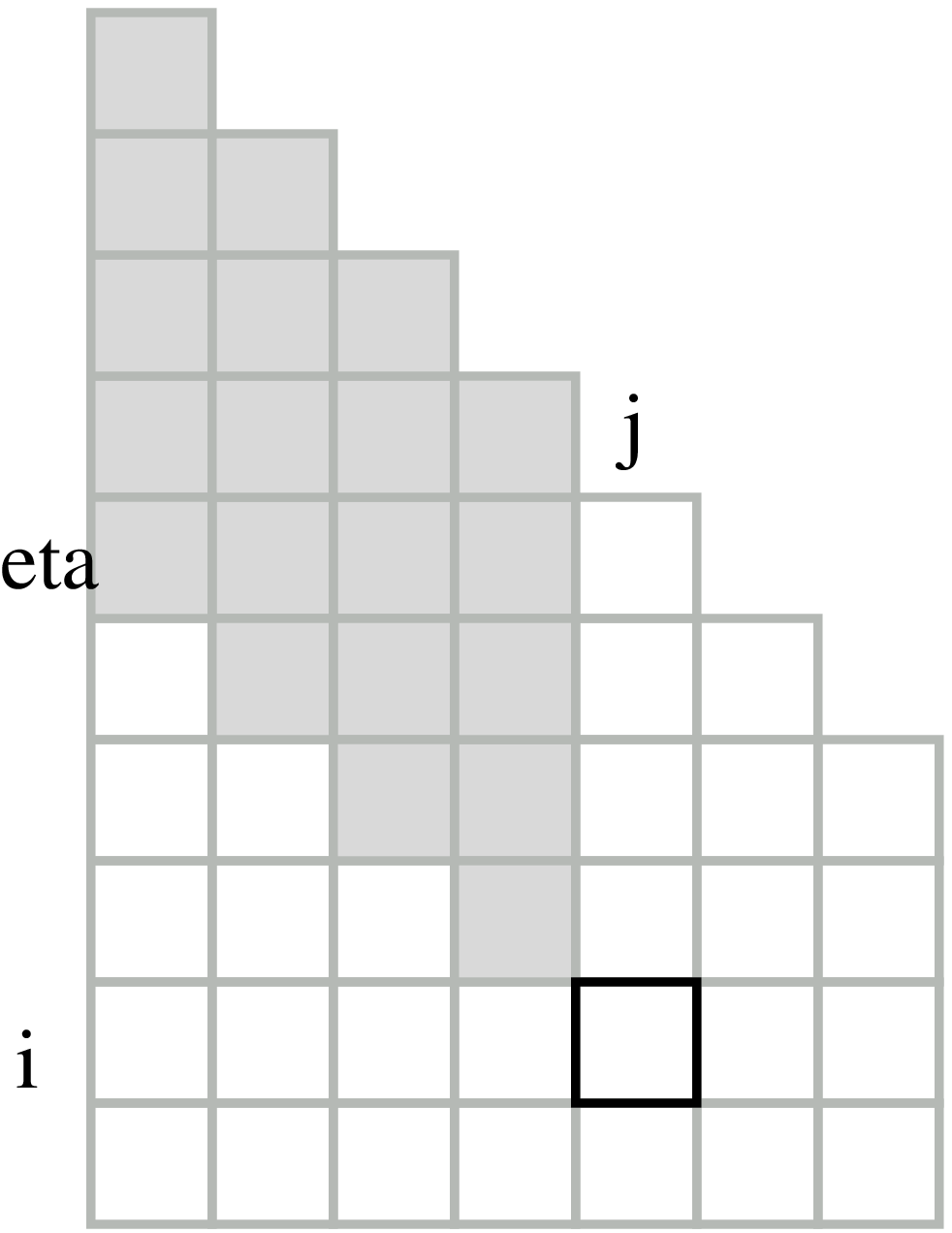}\label{sfl:sep1}\hspace{10ex}
  \includegraphics[height=.17\textheight]{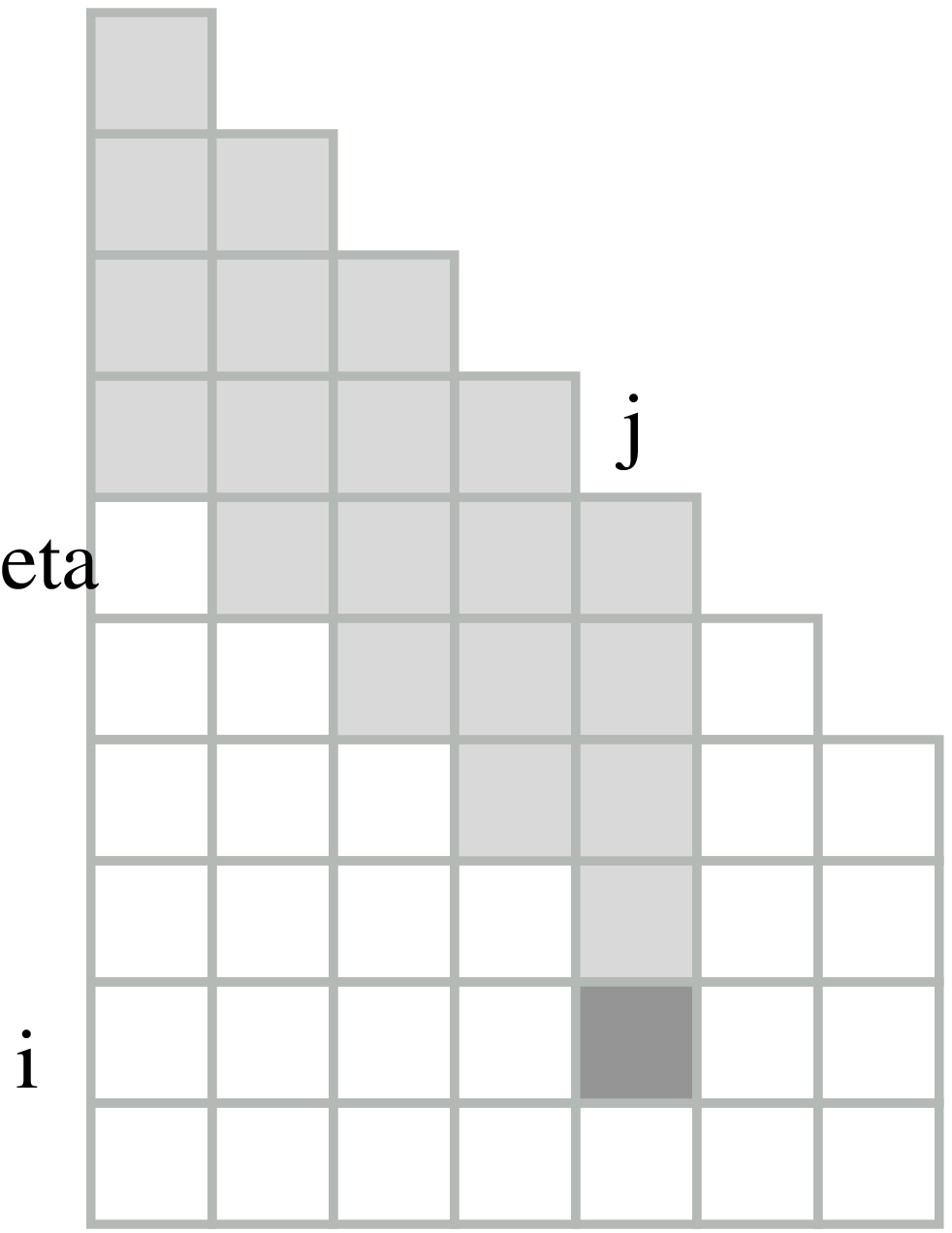}\label{sfl:sep2}
  \caption[]{The two cases arising in the dynamic programming algorithm of
    Section~\ref{sec:SeparationAlgorithmForSCIs}.}
  \label{fig:sep}
\end{figure}

In order to devise an efficient separation algorithm for SCIs, we need a
method to compute minimal shifted columns with respect to a given weight
vector $w \in \Q^{\orbipartinds{p}{q}}$.  The crucial observation is the
following. Let $S = \{\diagcol{1}{c_1}, \diagcol{2}{c_2}, \dots,
\diagcol{\eta}{c_{\eta}}\}$ with $1 \leq c_1 \leq c_2 \leq \dots \leq
c_{\eta} \leq j$ be a shifting of $\colop\diagcol{\eta}{j}$ for
$\diagcol{\eta}{j} \in \orbipartinds{p}{q}$ with $\eta > 1$. If $c_{\eta} <
j$, then~$S$ is a shifting of $\colop\diagcol{\eta}{j-1}$ (\emph{Case 1}).
If $c_{\eta} = j$, then
\[
S - \diagcol{\eta}{j} = \{\diagcol{1}{c_1}, \diagcol{2}{c_2}, \dots, \diagcol{\eta-1}{c_{\eta-1}}\}
\]
is a shifting of $\colop\diagcol{\eta-1}{j}$ (\emph{Case 2}); see
Figure~\ref{fig:sep}.

For all $\diagcol{\eta}{j} \in \orbipartinds{p}{q}$, let
$\omega\diagcol{\eta}{j}$ be the weight of a $w$-minimal shifting of
$\colop\diagcol{\eta}{j}$. The table $(\omega\diagcol{\eta}{j})$ can be
computed by dynamic programming as follows; we also compute a table of
values $\tau\diagcol{\eta}{j} \in \{1,2\}$, for each $\diagcol{\eta}{j}$,
which are needed later to reconstruct the corresponding shifted columns:

\begin{myenumerate}
\item For $j = 1, 2, \dots, q$, initialize $\omega\diagcol{1}{j} := \min
  \setdef{w_{\diagcol{1}{\ell}}}{\ell \in \ints{j}}$.
\item For $\eta = 2, 3, \dots, p$, initialize $\omega\diagcol{\eta}{1} :=
  \omega\diagcol{\eta-1}{1} + w_{\diagcol{\eta}{1}}$.
\item For $\eta = 2, 3, \dots, p$, $j = 2, 3, \dots, q$ (with
  $\diagcol{\eta}{j} \in \orbipartinds{p}{q}$): Compute
  \[
  \omega_1 := \omega\diagcol{\eta}{j-1} \quad\text{and}\quad
  \omega_2 := \omega\diagcol{\eta-1}{j} + w_{\diagcol{\eta}{j}}
  \]
  corresponding to Cases~1 and~2, respectively.
  Then set
  \[
  \omega\diagcol{\eta}{j} = \min \{\omega_1,\; \omega_2\}
  \quad\text{and}\quad
  \tau\diagcol{\eta}{j} =
  \begin{cases}
    1 & \text{if }\omega_1 \leq \omega_2\\
    2 & \text{otherwise}.
  \end{cases}
  \]
\end{myenumerate}

Thus, the tables $(\omega\diagcol{\eta}{j})$ and $(\tau\diagcol{\eta}{j})$
can be computed in time $\bigo{pq}$.  Furthermore, for a given
$\diagcol{\eta}{j} \in \orbipartinds{p}{q}$, we can compute a $w$-minimal
shifting $S\diagcol{\eta}{j}$ of $\colop\diagcol{\eta}{j}$ in time
$\bigo{\eta}$ from the table~$(\tau\diagcol{\eta}{j})$: We have
$S\diagcol{1}{j} = \{\diagcol{1}{j}\}$ for all $j \in \ints{q}$,
$S\diagcol{\eta}{1} = \colop\diagcol{\eta}{1}$ for all $\eta \in \ints{p}$,
and
\[
S\diagcol{\eta}{j} =
\begin{cases}
  S\diagcol{\eta}{j-1} & \text{if }\tau\diagcol{\eta}{j} = 1 \\
  S\diagcol{\eta-1}{j}\cup\{\diagcol{\eta}{j}\} & \text{if }
  \tau\diagcol{\eta}{j} = 2
\end{cases}
\]
for all other $\diagcol{\eta}{j}$. This proves the following result.

\begin{theorem}\label{thm:sepSCI}
  Let $w \in \Q^{\orbipartinds{p}{q}}$ be a given weight vector. There is
  an $\bigo{pq}$ time algorithm that simultaneously computes the weights of
  $w$-minimal shiftings of $\colop\diagcol{\eta}{j}$ for all
  $\diagcol{\eta}{j} \in \orbipartinds{p}{q}$ and a data structure that
  afterwards, for a given $\diagcol{\eta}{j}$, allows to determine a
  corresponding shifted column in~$\bigo{\eta}$ steps.
\end{theorem}

In particular, we obtain the following:

\begin{corollary}\label{cor:sepSCI}
  The separation problem for shifted column inequalities can be solved in
  linear time $\bigo{pq}$.
\end{corollary}

\begin{proof}
  Let a point $x^{\star} \in \Q^{\orbipartinds{p}{q}}$ be given. We can
  compute the $x^{\star}$-values $\beta(i,j) := x^{\star}(B(i,j))$ of all
  bars $B(i,j) = \{(i,j), (i,j+1), \dots, (i,\min\{i,q\})\}$ in linear time
  in the following way: First, we initialize $\beta(i,\ell) =
  x^{\star}_{i\ell}$ for all $i \in \ints{p}$ and $\ell = \min\{i,q\}$.
  Then, for each $i\in\ints{p}$, we calculate the value $\beta(i,j) =
  x^{\star}_{ij} + \beta(i,j+1)$ for $j = \min\{i,q\}-1, \min\{i,q\}-2,
  \dots, 1$.

  Using Theorem~\ref{thm:sepSCI} (and the notations introduced in the
  paragraphs preceeding it), we compute the table
  $(\omega\diagcol{\eta}{j})$ and the mentioned data structure in time
  $\bigo{pq}$. Then in time $\bigo{pq}$ we check whether there exists an
  $(i,j) = \diagcol{\eta}{j}\in\orbipartinds{p}{q}$ with $j \geq 2$ and
  $\omega\diagcol{\eta}{j-1} < \beta(i,j)$. If there exists such an
  $\diagcol{\eta}{j}$, we compute the corresponding shifted column
  $S\diagcol{\eta}{j-1}$ (in additional time $\bigo{\eta} \subseteq
  \bigo{p}$), yielding an SCI that is violated by~$x^{\star}$. Otherwise
  $x^{\star}$ satisfies all SCIs.
\end{proof}

Of course, the procedure described in the proof of the corollary can be
modified to find a maximally violated SCI if $x ^{\star}$ does not satisfy
all SCIs.

% ------------------------------------------------------
% Complete inequality descriptions
% ------------------------------------------------------

\subsection{Complete Inequality Descriptions}
\label{sec:CompleteInequalityDescriptions}

In this section we prove that nonnegativity constraints, row-sum equations,
and SCIs suffice to describe partitioning and packing orbitopes for
symmetric groups. The proof will be somewhat more involved than in the case
of cyclic groups. In particular, the coefficient matrices are not
totally unimodular anymore. In order to see this, consider the three column
inequalities
\begin{align*}
& x_{3,3}-x_{2,2}\leq 0,\quad
x_{4,3}+x_{4,4}-x_{2,2}-x_{3,2}\leq 0,
\quad \text{and} \\
& x_{5,4}+x_{5,5}-x_{3,3}-x_{4,3}\leq 0.
\end{align*}
The submatrix of the coefficient matrix belonging to these three rows and
the columns corresponding to $(2,2)$, $(3,3)$, and $(4,3)$ is the matrix
\[
\left(
\begin{array}{rrr}
-1 & +1 & 0 \\
-1 &  0 & +1 \\
0 & -1 & -1
\end{array}
\right),
\]
whose determinant equals $-2$. Note that the above three inequalities
define facets both of~$\orbipack{p}{q}$ and $\orbipart{p}{q}$ for $p\geq q
\geq 5$ (see Propositions~\ref{prop:orbipack:facets}
and~\ref{prop:orbipart:facets}, respectively).

\begin{proposition}\label{prop:orbipart:descr}
  The partitioning orbitope $\orbipart{p}{q}$ is completely described by
  the nonnegativity constraints, the row-sum equations, and the shifted
  column inequalities:
  \begin{align*}
    \orbipart{p}{q} =
    \{ \, & x \in \R^{\orbipartinds{p}{q}} \suchthat
    x \geq \zerovec,\; x(\row{i}) = 1 \text{ for }i = 1, \dots, p,\\
    & x(B) - x(S) \leq 0 \text{ for all SCIs with SC } S \text{ and bar }B\, \}.
  \end{align*}
\end{proposition}

\begin{proof}
  Let~$P$ be the polyhedron on the right-hand side of the statement above.
  From Propositions~\ref{prop:IPsym} and~\ref{prop:SCIvalid} we know
  already that
  \[
  P \cap \Z^{\orbipartinds{p}{q}} = \orbipart{p}{q} \cap \Z^{\orbipartinds{p}{q}}
  \]
  holds. Thus, it suffices to show that~$P$ is an integral polytope
  (as~$\orbipart{p}{q}$ is by definition). In the following, we first
  describe the strategy of the proof.

  For the rest of the proof, fix an arbitrary vertex~$x^\star$ of~$P$. A
  \emph{basis}~$\mathcal{B}$ of~$x^{\star}$ is a cardinality
  $\card{\orbipartinds{p}{q}}$ subset of the constraints describing~$P$
  that are satisfied with equality by~$x^{\star}$ with the property that
  the $\card{\orbipartinds{p}{q}} \times
  \card{\orbipartinds{p}{q}}$-coefficient matrix of the left-hand sides of
  the constraints in~$\mathcal{B}$ is non-singular.  Thus, the equation
  system obtained from the constraints in~$\mathcal{B}$ has~$x^{\star}$ as
  its unique solution.

  We will show that there exists a basis ~$\basis^{\star}$ of~$x^{\star}$
  that does not contain any SCI. Thus, $\basis^{\star}$ contains a
  subset of the~$p$ row-sum equations and at
  least~$\card{\orbipartinds{p}{q}}-p$ nonnegativity constraints. This
  shows that~$x^{\star}$ has at most~$p$ nonzero entries and,
  since~$x^{\star}$ satisfies the row-sum equations, it has a nonzero entry
  in every row.  Therefore, $\basis^{\star}$ contains all~$p$ row-sum
  equations, and all~$p$ nonzero entries must in fact be~$1$.
  Hence,~$x^{\star}$ is a 0/1-point. So the existence of such a basis
  proves the proposition.

  The \emph{weight} of a shifted column $S = \{\diagcol{1}{c_1},
  \diagcol{2,c_2}, \dots, \diagcol{\eta}{c_{\eta}}\}$ with $1 \leq c_1 \leq
  c_2 \leq \dots \leq c_{\eta} < q$ (we will not need shifted columns with
  $c_{\eta} = q$ here, as they do not appear in SCIs) is
  \[
  \weightop(S) := \sum_{i=1}^{\eta} c_i\, q^i.
  \]
  In particular, if $S_1$ and~$S_2$ are two shifted columns with
  $\card{S_1} < \card{S_2}$, then we have $\weightop(S_1)<\weightop(S_2)$.
  The \emph{weight} of an SCI is the weight of its shifted column, and the
  \emph{weight} of a basis~$\basis$ is the sum of the weights of the
  SCIs contained in~$\basis$ (note that a shifted column can appear in
  several SCIs).

  A basis of~$x^{\star}$ that contains all row-sum equations and all
  nonnegativity constraints corresponding to $0$-entries of~$x^{\star}$ is
  called \emph{reduced}. As the coefficient vectors (of the
  left-hand sides) of these constraints are linearly independent, some
  reduced basis of~$x^{\star}$ exists. Hence, there is also a reduced
  basis~$\basis^{\star}$ of~$x^{\star}$ of minimal weight.

  To prove the proposition, it thus suffices to establish the following claim.

  \begin{claim}\label{claim:minrednoSCI}
    A reduced basis of~$x^{\star}$ of minimal weight does not contain any
    SCI.
  \end{claim}

  The proof of Claim~\ref{claim:minrednoSCI} consists of three parts:
  \begin{myenumerate}
  \item We show that a reduced basis of~$x^{\star}$ does not contain any
    ``trivial SCIs'' (Claim~\ref{claim:TrivialSCIs}).
  \item We prove that a reduced basis of~$x^{\star}$ of minimal weight
    satisfies three structural conditions on its (potential) SCIs
    (Claim~\ref{claim:BasicStructure}).
  \item Finally, assuming that a reduced basis of~$x^{\star}$ with minimal
    weight contains at least one SCI, we will derive a contradiction by
    constructing a different solution~$\Tilde{x}\not=x^{\star}$ of the
    corresponding equation system.
  \end{myenumerate}

  We are now ready to start with Part~1. We call an SCI with shifted
  column~$S$ \emph{trivial} if $x^{\star}(S) = 0$ holds or if we have
  $x^{\star}(S) = 1$ and $x^{\star}_{k\ell} = 0$ for all $(k,\ell) \in S -
  (i,j)$ for some $(i,j) \in S$ (thus satisfying $x^{\star}_{ij}=1$) (see
  Figure~\ref{fig:forbidden}~\subref{sfl:forbidden1}).

  \begin{figure}
    \centering
    \newcommand{\mystyle}[1]{\footnotesize {#1}}
    \psfrag{0}{\mystyle{0}}
    \psfrag{1}{\mystyle{1}}
    \psfrag{?}{\mystyle{?}}
    \psfrag{*}{\mystyle{$\star$}}
    \subfloat[][]{\includegraphics[height=.17\textheight]{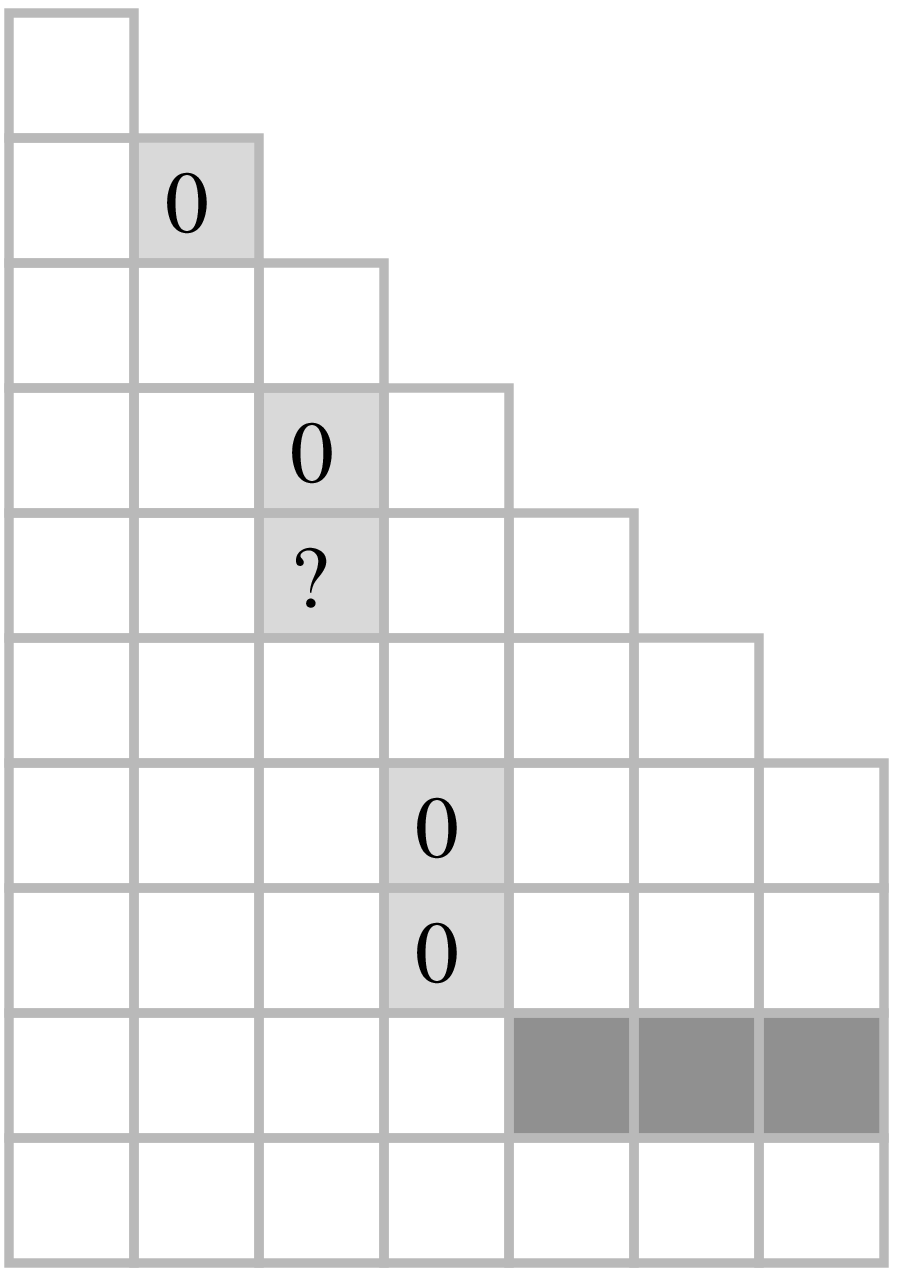}\label{sfl:forbidden1}}\hfill
    \subfloat[][]{\includegraphics[height=.17\textheight]{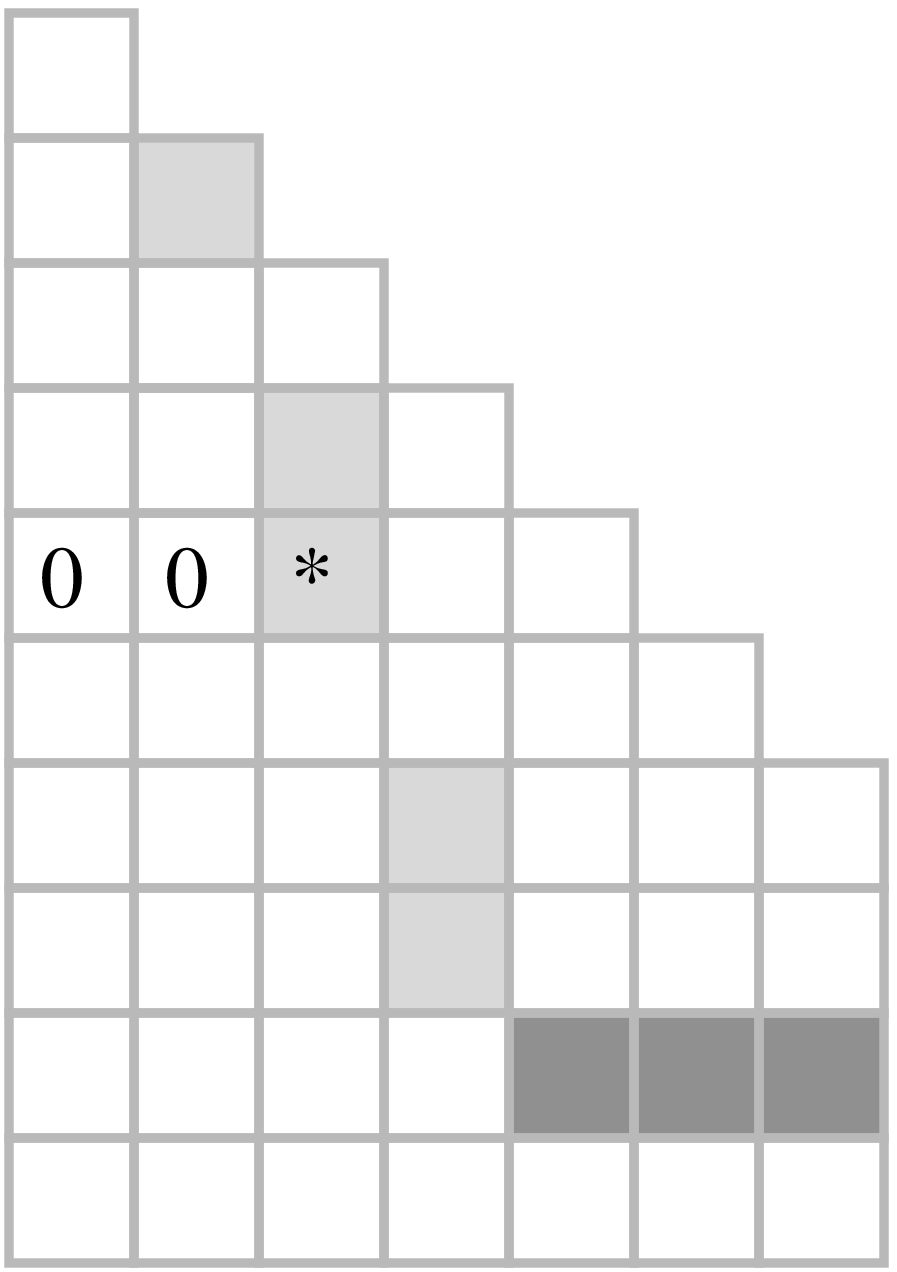}\label{sfl:forbidden2}}\hfill
    \subfloat[][]{\includegraphics[height=.17\textheight]{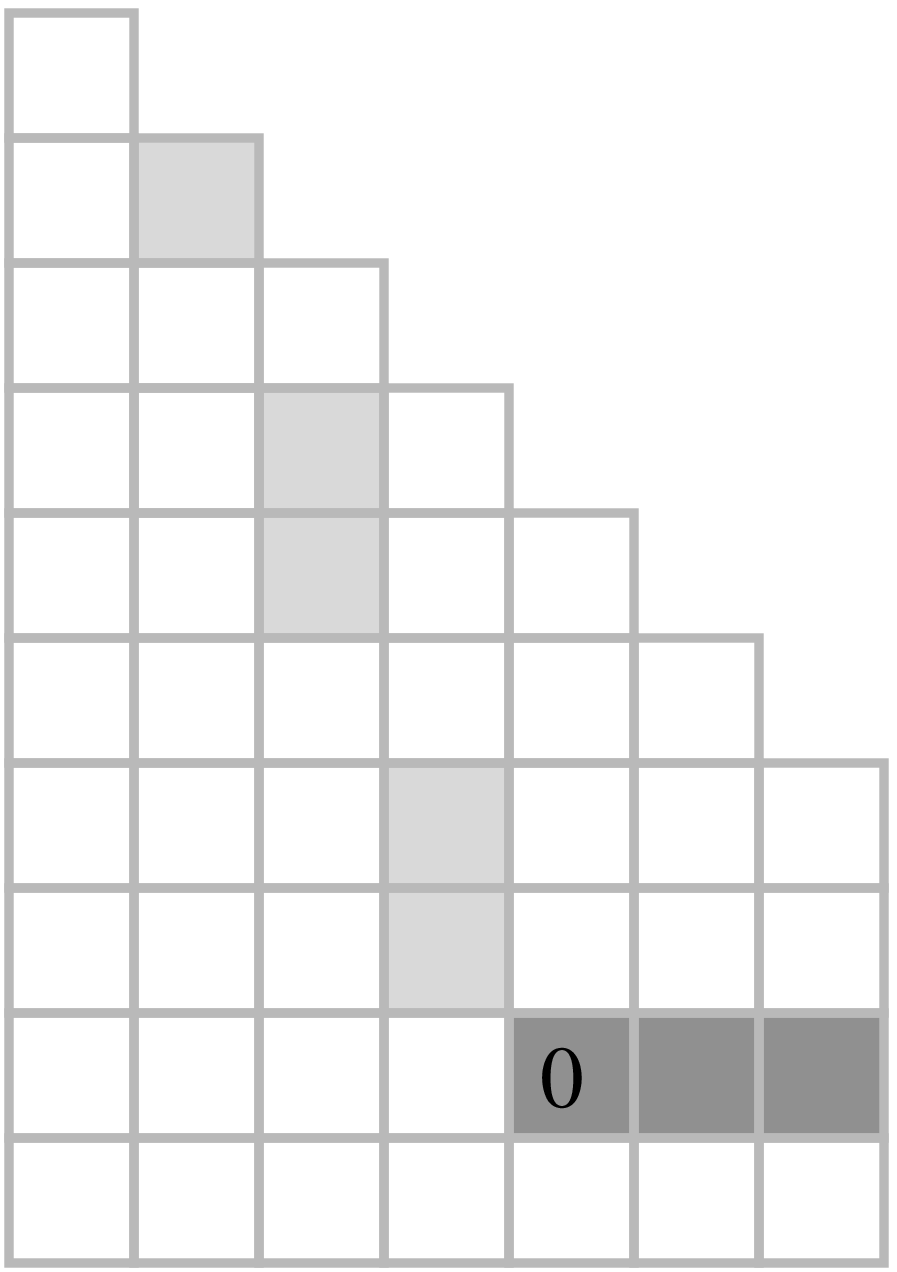}\label{sfl:forbidden3}}\hfill
    \subfloat[][]{\includegraphics[height=.17\textheight]{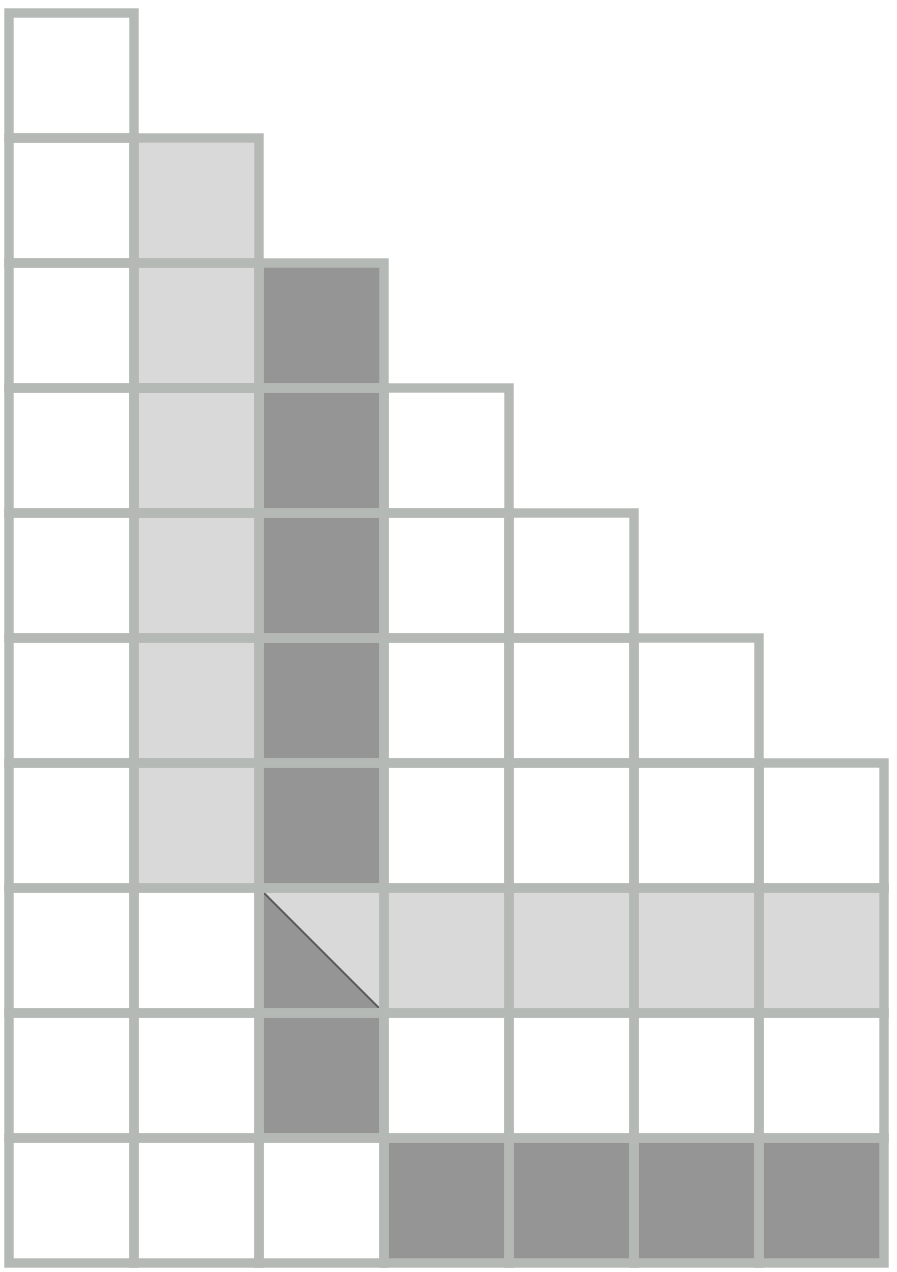}\label{sfl:forbidden4}}
    \caption[]{Illustration of trivial SCIs and of the three types of
      configurations not present in reduced bases of minimal weight, see
      Claim~\ref{claim:BasicStructure}. Bars are shown in dark gray,
      shifted columns in light gray. Figure~\subref{sfl:forbidden1} shows
      trivial SCIs (``?'' refers to a $0$ or $1$).
      Figures~\subref{sfl:forbidden2}, \subref{sfl:forbidden3}, and
      \subref{sfl:forbidden4} refer to parts \eqref{reduced:firstRowEntry},
      \eqref{reduced:nonzeroLeader}, and \eqref{reduced:containLeader} of
      Claim~\ref{claim:BasicStructure}, respectively (``$\star$'' indicates
      any nonzero number).  }
    \label{fig:forbidden}
  \end{figure}

  \begin{claim}\label{claim:TrivialSCIs}
    A reduced basis~$\basis$ of~$x^{\star}$ does not contain any trivial
    SCIs.
  \end{claim}
  \begin{proof}
    Let~$S$ be the shifted column~$S$ and~$B$ be the bar of some SCI that
    is satisfied with equality by~$x^{\star}$.

    If~$x^{\star}(S)=0$, then the coefficient vector of the SCI is a linear
    combination of the coefficient vectors of the inequalities $x_{ij} \geq
    0$ for $(i,j) \in S\cup B$, which all are contained in~$\basis$ (due to
    $x^{\star}(B)=x^{\star}(S)=0$). Since the coefficient vectors of the
    inequalities in~$\basis$ form a non-singular matrix, the SCI can not be
    in~$\basis$. (By ``coefficient vector'' we always mean the vector
    formed by the coefficients of the left-hand side of a constraint.)

    If~$S$ contains exactly one entry $(k,\ell) \in S$ with
    $x^\star_{k\ell} = 1$, then we have $x^{\star}(S) = x^{\star}(B) = 1$.
    Let~$i$ be the index of the row that contains the bar~$B$. The
    nonnegativity constraints $x_{rs} \geq 0$ for $(r,s) \in S - (k,\ell)$,
    $x_{ks} \geq 0$ for $(k,s) \in \row{k} - (k,\ell)$, and $x_{is} \geq 0$
    for $(i,s) \in \row{i} \setminus B$ are contained in~$\basis$.

    Since the coefficient vector of the considered SCI can linearly be
    combined from the coefficient vectors of these nonnegativity
    constraints and of the row-sum equations $x(\row{k}) = 1$ and
    $x(\row{i}) = 1$, this SCI cannot be contained in~$\basis$.
  \end{proof}

  \begin{claim}\label{claim:BasicStructure}
    A minimal weight reduced basis~$\basis$ of~$x^{\star}$ satisfies the
    following three conditions:
    \begin{myenumerate}
    \item\label{reduced:firstRowEntry}%
      If $(k,\ell)$ is contained in the shifted column of some SCI
      in~$\basis$, then there exists some $s < \ell$ with $x^{\star}_{ks} >
      0$.
    \item\label{reduced:nonzeroLeader}%
      If $(i,j)$ is the leader of an SCI in~$\basis$, then $x^\star_{ij} >
      0$ holds.
    \item\label{reduced:containLeader}%
      If $(i,j)$ is the leader of an SCI in~$\basis$, then there is no SCI
      in $\basis$ whose shifted column contains $(i,j)$.
    \end{myenumerate}
    See Figure~\ref{fig:forbidden},
    \subref{sfl:forbidden2}--\subref{sfl:forbidden4} for an illustration of
    the three conditions.
  \end{claim}

  \begin{figure}
    \centering
    \newcommand{\mystyle}[1]{\footnotesize {#1}}
    \psfrag{i}{\mystyle{$i$}}
    \psfrag{j}{\mystyle{$j$}}
    \psfrag{k}{\mystyle{$k$}}
    \psfrag{l}{\mystyle{$\ell$}}
    \psfrag{0}{\mystyle{$0$}}
    \psfrag{C}{\mystyle{$C$}}
    \psfrag{B}{\mystyle{\textcolor{white}{$B$}}}
    \psfrag{B'}{\mystyle{$B'$}}
    \psfrag{S'}{\mystyle{$S'$}}
    \psfrag{S3}{\mystyle{$S_3$}}
    \psfrag{S4}{\mystyle{$S_4$}}
    \psfrag{B2}{\mystyle{$B_2$}}
    \psfrag{B3}{\textcolor{white}{\mystyle{$B_3$}}}
    \begingroup
    \renewcommand{\thesubfigure}{\arabic{subfigure}}
    \subfloat[][]{\includegraphics[height=.17\textheight]{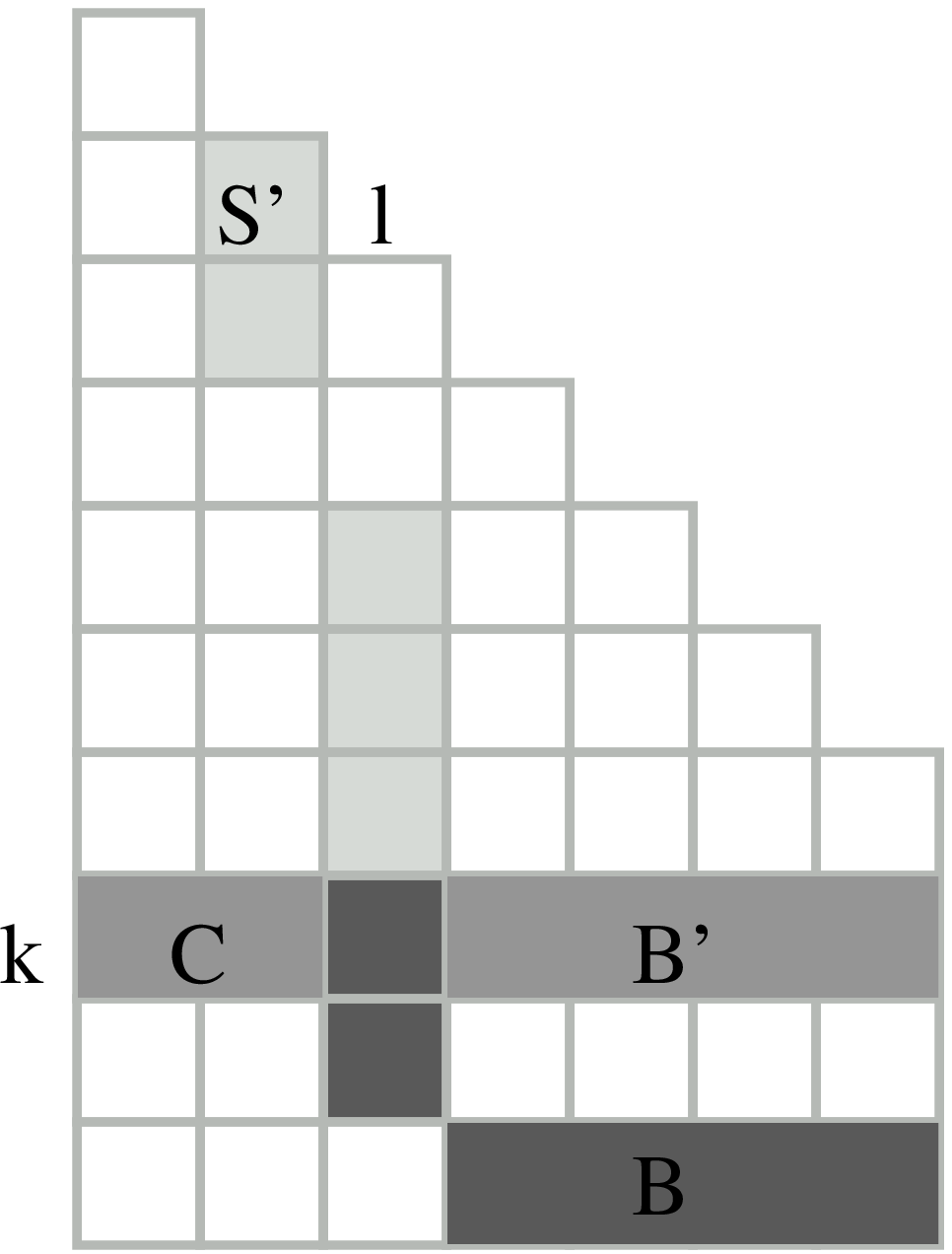}\label{sfl:ensure1}}\qquad
    \subfloat[][]{\includegraphics[height=.17\textheight]{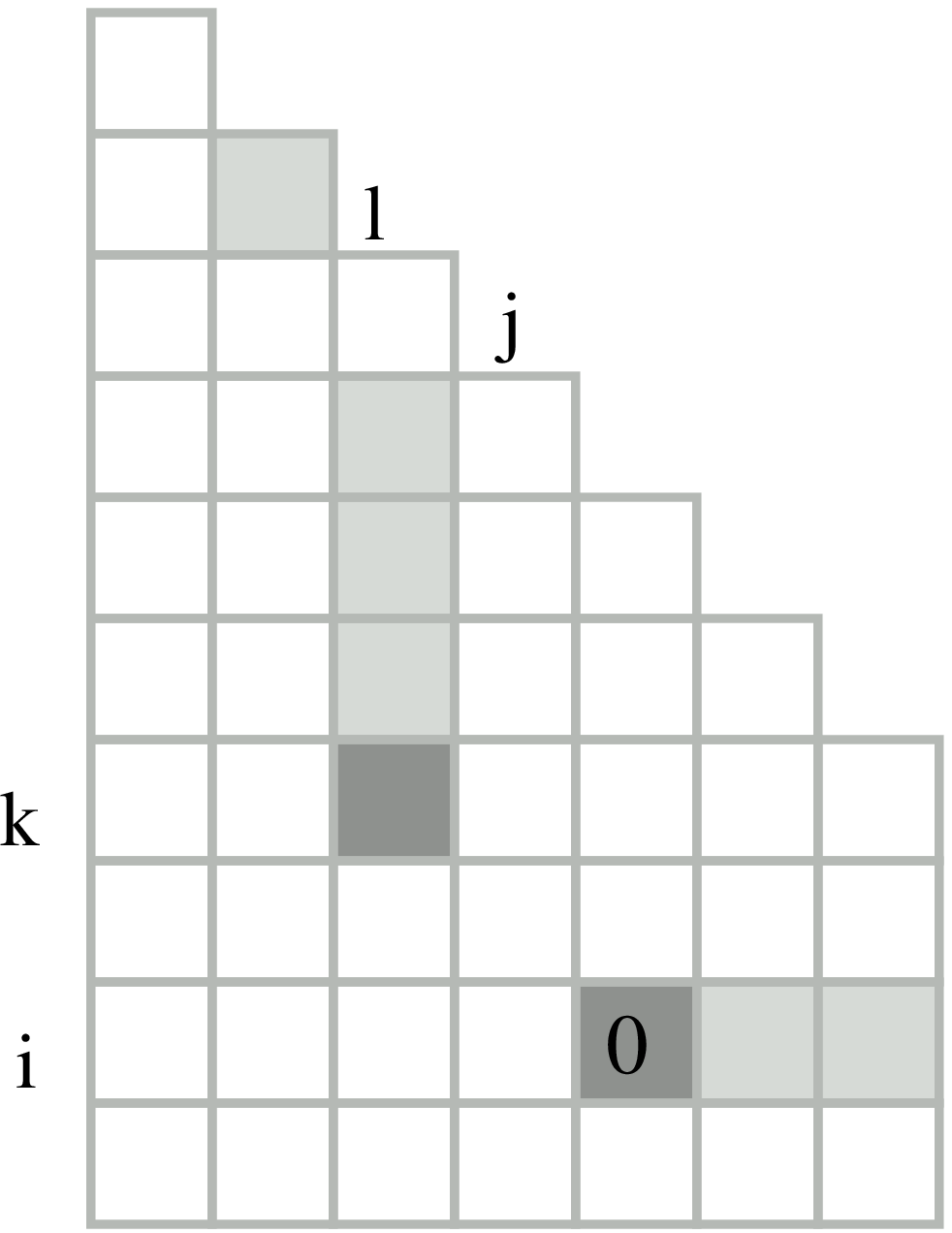}\label{sfl:ensure2}}\qquad
    \subfloat[][]{\includegraphics[height=.17\textheight]{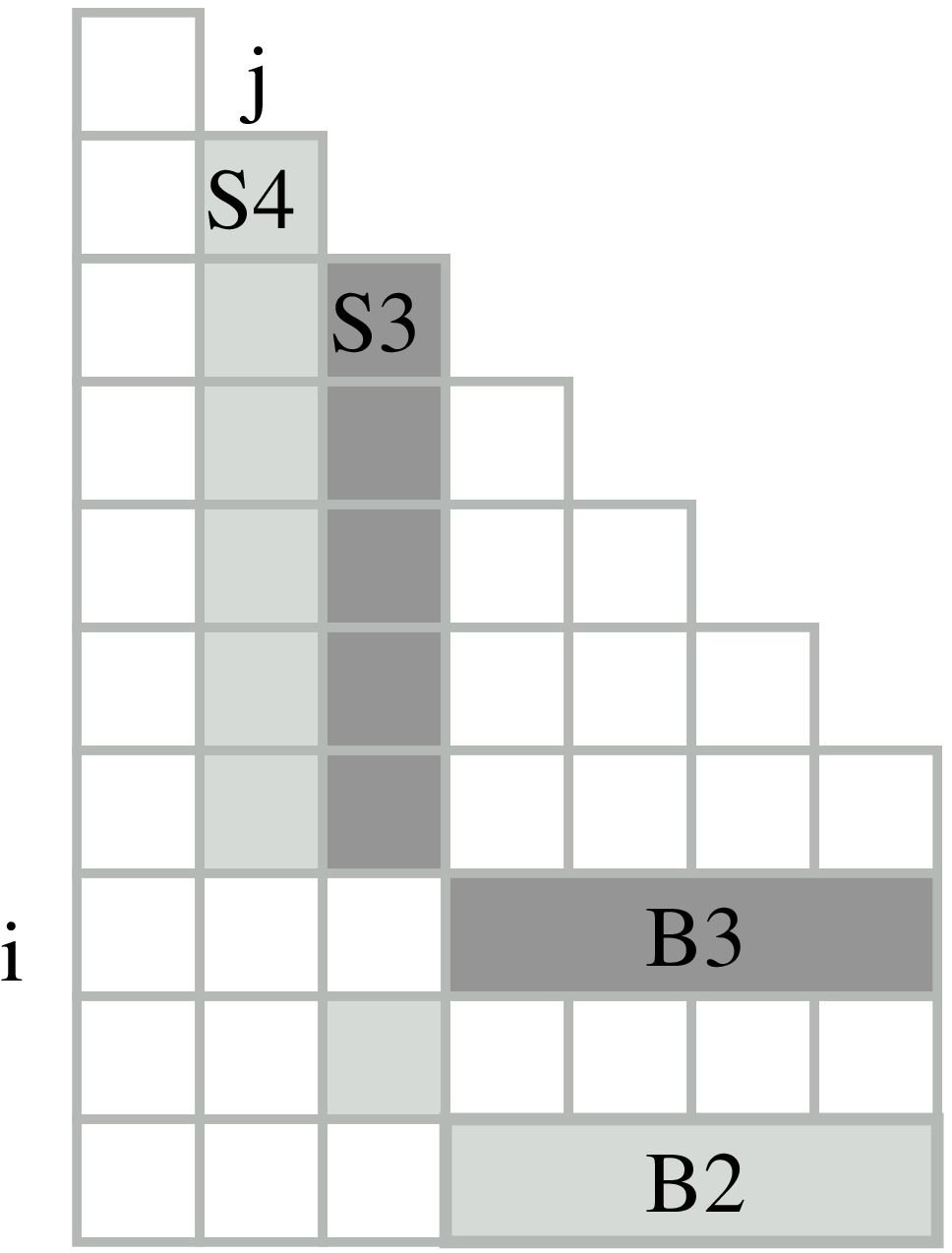}\label{sfl:ensure3}}
    \endgroup
    \caption[]{Illustration of the proof of
      Claim~\ref{claim:BasicStructure},
      parts~\eqref{reduced:firstRowEntry} to~\eqref{reduced:containLeader}.}
    \label{fig:basicStructure}
  \end{figure}

  \begin{proof}
    \begin{partslist}
    \partsitem Assume there exists an SCI in~$\basis$ with shifted
      column~$S$ and bar~$B$ that contains the first nonzero entry of a
      row~$k$, i.e., there is $(k,\ell) \in S$ with $x^\star_{k\ell} > 0$
      and $x^\star_{ks} = 0$ for all $s < \ell$. Let $S' := S \cap
      \orbipartinds{k-1}{q}$ be the entries of~$S$ above row~$k$. Let $C =
      \{(k,1), (k,2), \dots, (k,\ell-1)\}$ and $B' = \row{k} \setminus
      (C+(k,\ell))$. See
      Figure~\ref{fig:basicStructure}~\subref{sfl:ensure1} for an
      illustration.

      Because~$S'$ is a shifting of $\colop(k-1,\ell)$, $x(B') - x(S') \leq
      0$ is an SCI and hence satisfied by~$x^{\star}$. Since we have
      $\card{S'} < \card{S}$ (thus, $\weightop(S') < \weightop(S)$), it
      suffices to show that replacing the original SCI $x(B) - x(S) \leq 0$
      by $x(B') - x(S') \leq 0$ gives another basis~$\basis'$
      of~$x^{\star}$ (which also is reduced), contradicting the minimality
      of the weight of~$\basis$.

      Due to $x^{\star}(\row{k}) = 1$, $x^{\star}(C) = 0$, $x^{\star}(B') -
      x^{\star}(S') \leq 0$, and $S' + (k,\ell) \subseteq S$ we have
      \begin{equation}\label{eq:basicStructureChain1}
        1 = x^\star_{k\ell} + x^\star(B') \leq x^\star_{k\ell} + x^\star(S')
        \leq x^\star(S) =x^{\star}(B) \leq 1.
      \end{equation}
      Therefore, equality must hold throughout this chain. In particular,
      this shows $x^{\star}(B') - x^{\star}(S') = 0$. Thus, its suffices to
      show that the coefficient matrix of the equation system obtained
      from~$\basis'$ is non-singular, which can be seen as follows.

      Since $x^{\star}(S' + (k,\ell)) = 1 = x^{\star}(S)$
      (see~\eqref{eq:basicStructureChain1}), we know that all nonnegativity
      constraints $x_{rs} \geq 0$ with $(r,s) \in S \setminus (S' +
      (k,\ell))$ are contained in~$\basis$ and~$\basis'$. The same holds
      for $x_{ks} \geq 0$ with $(k,s) \in C$ and for $x_{is} \geq 0$ with
      $(i,s) \in \row{i} \setminus B$, where row~$i$ contains bar~$B$
      (since $x^{\star}(B) = 1$ by~\eqref{eq:basicStructureChain1}). Thus,
      we can linearly combine the coefficient vector of $x(B)-x(S)\le 0$
      from the coefficient vectors of the constraints $x(B')-x(S')\le 0$,
      $x(\row{k})=1$, $x(\row{i})=1$, and the nonnegativity constraints
      mentioned above. Since all these constraints are contained
      in~$\basis'$, this shows that the coefficient matrix of~$\basis'$ has
      the same row-span as that of~$\basis$, thus proving that it is
      non-singular as well.

    \partsitem Assume that there exists an SCI in~$\basis$ with leader $(i,j)$,
      bar~$B$, and shifted column~$S$ such that $x^\star_{ij} = 0$. If $S =
      \{\diagcol{1}{c_1}, \diagcol{2}{c_2}, \dots,
      \diagcol{\eta}{c_{\eta}}\}$, then we have $(i,j) =
      \diagcol{\eta}{j}$.  Define $B' := B - (i,j)$, $S' := S -
      \diagcol{\eta}{c_{\eta}}$, and observe that $B' \neq \varnothing$, $S'
      \neq \varnothing$, i.e., $\card{B} > 1$ and $\card{S} > 1$, because a
      reduced basis does not contain trivial SCIs by
      Claim~\ref{claim:TrivialSCIs}; see
      Figure~\ref{fig:basicStructure}~\subref{sfl:ensure2}. Hence, $x(B') - x(S') \leq 0$ is an
      SCI. We therefore have:
      \begin{equation}\label{eq:basicStructureChain2}
        0 = x^{\star}(B) - x^{\star}(S) = x^\star(B') - x^\star(S) \leq x^{\star}(B') - x^{\star}(S') \leq 0,
      \end{equation}
      where the first equation holds because $x(B) - x(S) \leq 0$ is
      satisfied with equality by~$x^{\star}$ and the second equation
      follows from $x^{\star}_{ij} = 0$. Hence, we know that $x^{\star}(B')
      - x^{\star}(S') = 0$. Since we have $\card{S'} < \card{S}$ (and
      consequently $\weightop(S')<\weightop(S)$), again it remains to show
      that the coefficient vector of $x(B) - x(S) \leq 0$ can be linearly
      combined from the coefficient vector of $x(B') - x(S') \leq 0$ and
      some coefficient vectors of nonnegativity constraints in~$\basis$
      and~$\basis'$. But this is clear, as we have $x^{\star}_{ij} = 0$ and
      $x^{\star}_{\diagcol{\eta}{c_{\eta}}} = 0$, where the latter follows
      from~\eqref{eq:basicStructureChain2}.

    \partsitem Assume that in~$\basis$ there exists an SCI
      \begin{equation}\label{eq:claim2:SCI:1}
        x(B_1) - x(S_1) \leq 0
      \end{equation}
      with leader $(i,j) = \diagcol{\eta}{j}$, bar $B_1$, and shifted column
      \[
      S_1 = \{\diagcol{1}{c_1}, \diagcol{2}{c_2}, \dots, \diagcol{\eta}{c_{\eta}}\}
      \]
      (in particular: $c_{\eta} < j$) and another SCI
      \begin{equation}\label{eq:claim2:SCI:2}
        x(B_2) - x(S_2) \leq 0
      \end{equation}
      with bar $B_2$ and shifted column
      \[
      S_2 = \{\diagcol{1}{d_1}, \diagcol{2}{d_2}, \dots, \diagcol{\eta}{j},
      \diagcol{\eta+1}{d_{\eta+1}}, \dots, \diagcol{\tau}{d_{\tau}}\}.
      \]
      Hence, we have $(i,j) = \diagcol{\eta}{j} \in S_2$. Define
      \[
      S_3 := \{\diagcol{1}{d_1}, \diagcol{2}{d_2}, \dots, \diagcol{\eta-1}{d_{\eta-1}}\}
      \]
      (i.e, the part of $S_2$ lying strictly above row~$i$) and
      \[
      S_4 := \{\diagcol{1}{c_1}, \dots, \diagcol{\eta}{c_{\eta}},
      \diagcol{\eta+1}{d_{\eta+1}}, \dots, \diagcol{\tau}{d_{\tau}}\}
      \]
      (i.e, $S_1$ together with the part of $S_2$ strictly below row~$i$).
      Clearly, $S_3$ is a shifting of $\colop\diagcol{\eta-1}{j} =
      \colop(i-1,j)$, and $S_4$ is a shifted column as well (due to
      $c_{\eta} < j \leq d_{\eta+1}$). Thus, with $B_3 = B_1 - (i,j)$, we
      obtain the SCIs
      \begin{equation}\label{eq:claim2:SCI:3}
        x(B_3) - x(S_3) \leq 0
      \end{equation}
      \begin{equation}\label{eq:claim2:SCI:4}
        x(B_2) - x(S_4) \leq 0
      \end{equation}
      (see
      Figure~\ref{fig:basicStructure}~\subref{sfl:ensure3}).

      Since~\eqref{eq:claim2:SCI:1} and~\eqref{eq:claim2:SCI:2} are
      contained in~$\basis$, we have $x^{\star}(B_1) - x^{\star}(S_1) = 0$
      and $x^{\star}(B_2) - x^{\star}(S_2) = 0$. Adding these two equations
      yields
      \begin{equation}\label{eq:claim2:sumSCIs}
        \big(x^{\star}(B_3) - x^{\star}(S_3)\big) +
        \big(x^{\star}(B_2) - x^{\star}(S_4)\big) = 0,
      \end{equation}
      because $x^\star_{ij}$ cancels due to $(i,j) \in B_1 \cap S_2$.
      Since~$x^{\star}$ satisfies the SCIs~\eqref{eq:claim2:SCI:3}
      and~\eqref{eq:claim2:SCI:4}, Equation~\eqref{eq:claim2:sumSCIs} shows that in
      fact we have $x^{\star}(B_3) - x^{\star}(S_3) = 0$ and
      $x^{\star}(B_2) - x^{\star}(S_4) = 0$.

      It is not clear, however, that we can simply
      replace~\eqref{eq:claim2:SCI:1} and~\eqref{eq:claim2:SCI:2}
      by~\eqref{eq:claim2:SCI:3} and~\eqref{eq:claim2:SCI:4} in order to
      obtain a new basis of~$x^{\star}$. Nevertheless, if $v_1, v_2, v_3$,
      and~$v_4$ are the coefficient vectors of~\eqref{eq:claim2:SCI:1},
      \eqref{eq:claim2:SCI:2}, \eqref{eq:claim2:SCI:3},
      and~\eqref{eq:claim2:SCI:4}, respectively, we have $v_1 + v_2 = v_3 +
      v_4$, which implies
      \begin{equation}\label{eq:claim2:v2}
        v_2 = v_3 + v_4 - v_1.
      \end{equation}
      Let $V \subset \R^{\orbipartinds{p}{q}}$ be the subspace of
      $\R^{\orbipartinds{p}{q}}$ that is spanned by the coefficient vectors
      of the constraints different from~\eqref{eq:claim2:SCI:2}
      in~$\basis$. Thus, the linear span of $V \cup \{v_2\}$ is the whole
      space $\R^{\orbipartinds{p}{q}}$. Due to~\eqref{eq:claim2:v2}, the
      same holds for $V \cup \{v_3,v_4\}$ (since $v_1\in V$).  Therefore,
      there is $\alpha\in \{3,4\}$ such that $V \cup \{v_\alpha\}$ spans
      $\R^{\orbipartinds{p}{q}}$.  Let~$(a)$ be the corresponding SCI from
      $\{\eqref{eq:claim2:SCI:3}, \eqref{eq:claim2:SCI:4}\}$. Hence,
      $\basis' := \basis \setminus \{\eqref{eq:claim2:SCI:2}\} \cup
      \{(a)\}$ is a (reduced) basis of~$x^{\star}$ as well.

      Since we have $\card{S_3} < \card{S_2}$ and $\weightop(S_4) <
      \weightop(S_2)$ (due to $c_{\eta} < j$), the weight of~$\basis'$ is
      smaller than that of~$\basis$, contradicting the minimality of the
      weight of~$\basis$.
    \end{partslist}
  \end{proof}

  Before we finish the proof of the proposition by establishing
  Claim~\ref{claim:minrednoSCI}, we need one more structural result on the
  SCIs in a reduced basis of~$x^{\star}$.  Let $S = \{\diagcol{1}{c_1},
  \diagcol{2}{c_2}, \dots, \diagcol{\eta}{c_{\eta}}\}$ be any shifted
  column with $x^{\star}_{\diagcol{\gamma}{c_{\gamma}}} > 0$ for some
  $\gamma \in \ints{\eta}$. We call $\diagcol{\gamma}{c_{\gamma}}$ the
  \emph{first nonzero element} of~$S$ if
  \[
  x^{\star}_{\diagcol{1}{c_1}} = \dots =
  x^{\star}_{\diagcol{\gamma-1}{c_{\gamma-1}}} = 0
  \]
  holds. Similarly, $\diagcol{\gamma}{c_{\gamma}}$ is called the \emph{last
    nonzero element} of~$S$ if we have
  \[
  x^{\star}_{\diagcol{\gamma+1}{c_{\gamma+1}}} = \dots =
  x^{\star}_{\diagcol{\eta}{c_{\eta}}} = 0.
  \]

  \begin{claim}\label{claim:firstlast}
    Let~$\basis$ be a reduced basis of~$x^{\star}$, and let $S_1, S_2$ be
    the shifted columns of some SCIs in~$\basis$ ($S_1 = S_2$ is allowed).
    \begin{myenumerate}
    \item\label{claim:firstlast:firstnonzero}%
      If $(i,j)$ is the first nonzero element of~$S_1$ and $(i,j) \in S_2$,
      then $(i,j)$ is also the first nonzero element of~$S_2$.
    \item\label{claim:firstlast:lastnonzero}%
      If $(i,j)$ is the last nonzero element of~$S_1$ with $x^{\star}(S_1)
      = 1$ and $(i,j) \in S_2$, then $(i,j)$ is also the last nonzero
      element of~$S_2$ and $x^{\star}(S_2) = 1$.
    \item\label{claim:firstlast:both}%
      If $(i,j)$ is the last nonzero element of~$S_1$ with $x^{\star}(S_1)
      = 1$, then $(i,j)$ is not the first nonzero element of $S_2$.
    \end{myenumerate}
  \end{claim}

  \begin{figure}
    \centering
    \newcommand{\Mystyle}[1]{\footnotesize {#1}}
    \newcommand{\mystyle}[1]{\tiny {#1}}
    \psfrag{i}{\Mystyle{$i$}}
    \psfrag{j}{\Mystyle{$j$}}
    \psfrag{S1'}{\mystyle{$S_1'$}}
    \psfrag{S2'}{\textcolor{white}{\mystyle{$S_2'$}}}
    \psfrag{S2b'}{\mystyle{$\overline{S}_2'$}}
    \psfrag{B1}{\textcolor{white}{\mystyle{$B_1$}}}
    \psfrag{B2}{\mystyle{$B_2$}}
    \includegraphics[height=.17\textheight]{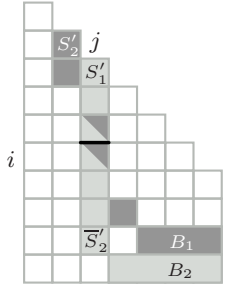}
    \caption[]{Illustration of sets used in the proof of Claim~\ref{claim:firstlast}.}
    \label{fig:firstlast}
  \end{figure}

  \begin{proof}
    Let
    \[
    S_1 = \{\diagcol{1}{c_1}, \diagcol{2}{c_2}, \dots,
    \diagcol{\eta}{c_{\eta}}\} \quad\text{and}\quad
    S_2 =
    \{\diagcol{1}{d_1},\diagcol{2}{d_2},\dots,\diagcol{\tau}{d_{\tau}}\}
    \]
    be two shifted columns of SCIs with bars~$B_1$ and~$B_2$, respectively,
    in the reduced basis~$\basis$ of~$x^{\star}$. Suppose that $(i,j) =
    \diagcol{\gamma}{j} \in S_1 \cap S_2$, i.e., $c_{\gamma} = j =
    d_{\gamma}$ holds. Define
    \begin{align*}
    S_1' & :=\{\diagcol{1}{c_1},\diagcol{2}{c_2},\dots,\diagcol{\gamma-1}{c_{\gamma-1}}\},\\
    S_2' & :=\{\diagcol{1}{d_1},\diagcol{2}{d_2},\dots,\diagcol{\gamma-1}{d_{\gamma-1}}\},
    \end{align*}
    and $\overline{S}_2' := S_2 \setminus S_2'$, see
    Figure~\ref{fig:firstlast}. Since $\diagcol{\gamma}{j} \in S_1 \cap
    S_2$ holds, $S'_1 \cup \overline{S'_2}$ is a shifted column and $x(B_2)
    - x(S'_1 \cup \overline{S'_2}) \leq 0$ is an SCI. Thus, we obtain
    \begin{equation}\label{eq:firstlast:1}
      x^\star(B_2) - x^\star(S_1') - x^{\star}(\overline{S}_2') \leq 0.
    \end{equation}
    Furthermore, since $x(B_2) - x(S_2) \leq 0$ is contained in the
    basis~$\basis$ of~$x^{\star}$, we have
    \begin{equation}\label{eq:firstlast:2}
      x^\star(B_2) - x^\star(S_2') - x^{\star}(\overline{S}_2') = 0.
    \end{equation}
    Subtracting~\eqref{eq:firstlast:2} from~\eqref{eq:firstlast:1} yields
    $x^\star(S_2') - x^\star(S_1') \leq 0$. We thus conclude
    \begin{equation}\label{eq:firstlast:3}
      x^\star(S_2') \leq x^\star(S_1')
      \quad\text{and}\quad
      x^\star(S_1') \leq x^\star(S_2')
    \end{equation}
    (where the second inequality follows by exchanging the roles of~$S_1$
    and~$S_2$ in the argument).
    \smallskip

    \begin{partslist}
    \partsitem If $(i,j)$ is the first nonzero element of~$S_1$, then we have
      $x^{\star}(S'_1) = 0$. Thus, the first inequality
      of~\eqref{eq:firstlast:3} implies $x^{\star}(S'_2) = 0$, showing that
      $(i,j)$ is the first nonzero element of~$S_2$.

    \partsitem If $(i,j)$ is the last nonzero element of~$S_1$ and
      $x^{\star}(S_1) = 1$ holds, then we have $x^{\star}(S'_1 + (i,j)) =
      1$. With the second inequality of~\eqref{eq:firstlast:3} we obtain:
      \[
      1 = x^{\star}(S'_1 + (i,j)) \leq x^{\star}(S'_2 + (i,j)) \leq x^{\star}(S_2) = x^{\star}(B_2) \leq 1,
      \]
      where the last equation holds because $x(B_2) - x(S_2) \leq 0$ is
      contained in~$\basis$. It follows that $x^{\star}(S_2) = 1$ and
      $(i,j)$ is the last nonzero element of~$S_2$.

    \partsitem This follows from the first two parts of the claim,
      since~$\basis$ does not contain any trivial SCIs by
      Claim~\ref{claim:TrivialSCIs}.
    \end{partslist}
  \end{proof}

  We will now proceed with the proof of Claim~\ref{claim:minrednoSCI}.
  Thus, assume that $\basis^{\star}$ is a reduced basis of~$x^{\star}$ of
  minimal weight and suppose that~$\basis^\star$ contains at least one SCI.
  We are going to construct a point~$\tilde{x} \neq x^{\star}$ that
  satisfies the equation system obtained from $\basis^{\star}$,
  contradicting the fact the~$x^{\star}$ is the unique solution to this
  system of equations.

  \begin{figure}
    \centering
    \newcommand{\Mystyle}[1]{\footnotesize {#1}}
    \newcommand{\mystyle}[1]{\tiny {#1}}
    \psfrag{-e}{\Mystyle{$-\lambda$}}
    \psfrag{+e}{\Mystyle{$+\lambda$}}
    \psfrag{*}{\Mystyle{$\star$}}
    \psfrag{0}{\Mystyle{$0$}}
    \psfrag{i}{\Mystyle{$i$}}
    \psfrag{j}{\Mystyle{$j$}}
    \psfrag{S1'}{\mystyle{$S_1'$}}
    \psfrag{S2'}{\mystyle{$S_2'$}}
    \psfrag{S1b'}{\mystyle{$\overline{S}_1'$}}
    \psfrag{S2b'}{\mystyle{$\overline{S}_2'$}}
    \psfrag{B1}{\mystyle{$B_1$}}
    \psfrag{B2}{\textcolor{white}{\mystyle{$B_2$}}}
    \subfloat[][]{\includegraphics[height=.17\textheight]{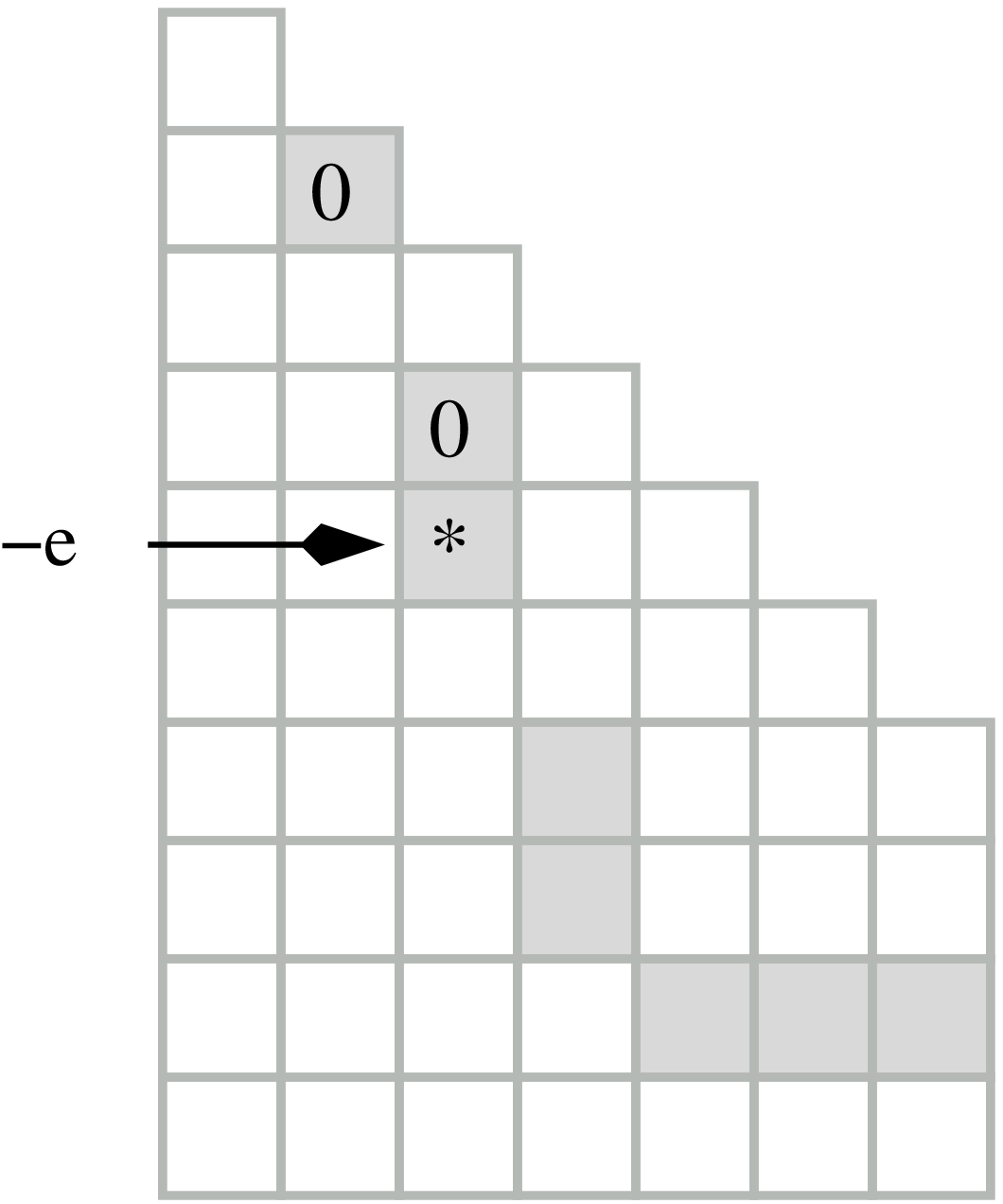}\label{sfl:finish1}}\qquad
    \subfloat[][]{\includegraphics[height=.17\textheight]{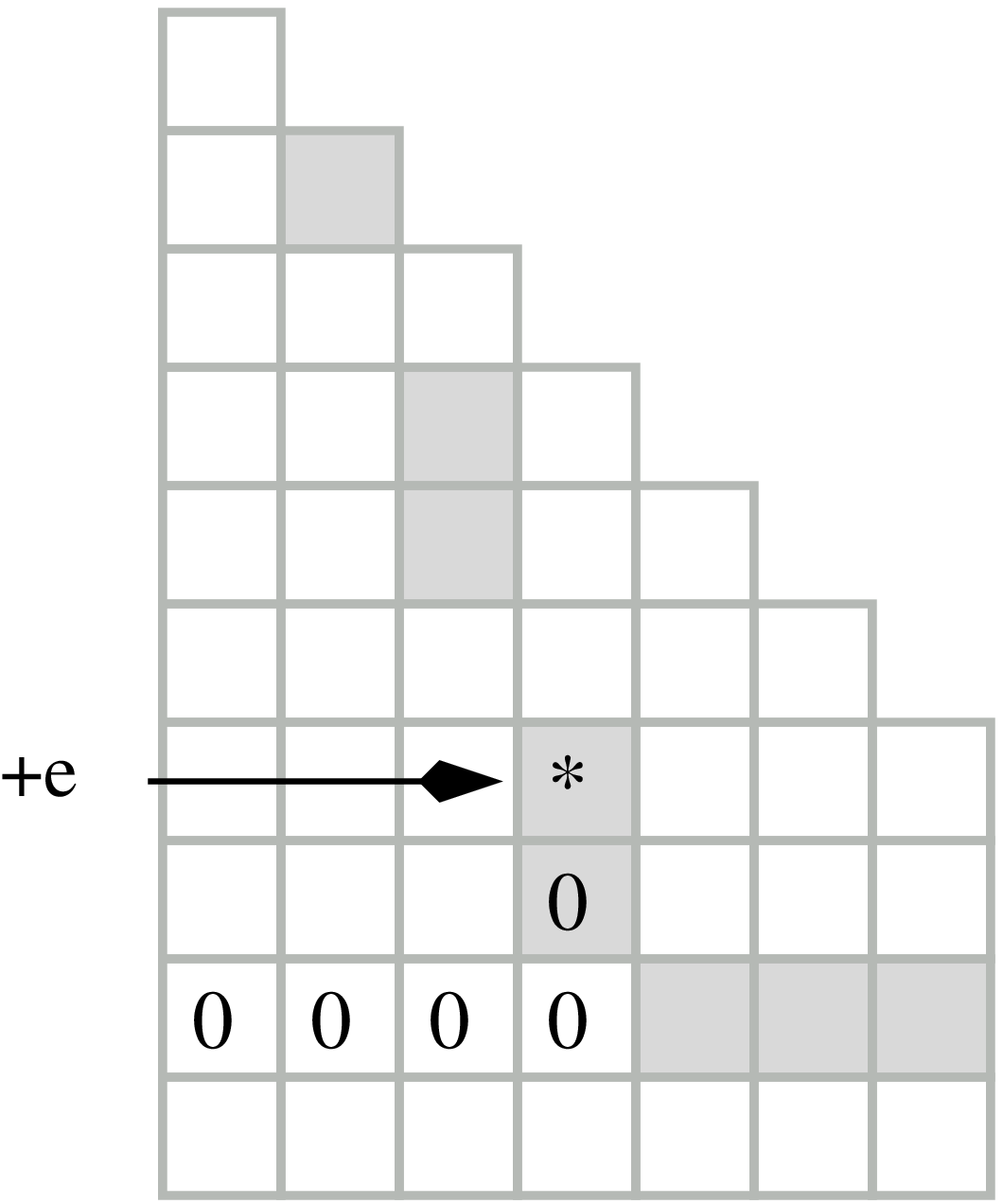}\label{sfl:finish2}}\qquad
    \subfloat[][]{\includegraphics[height=.17\textheight]{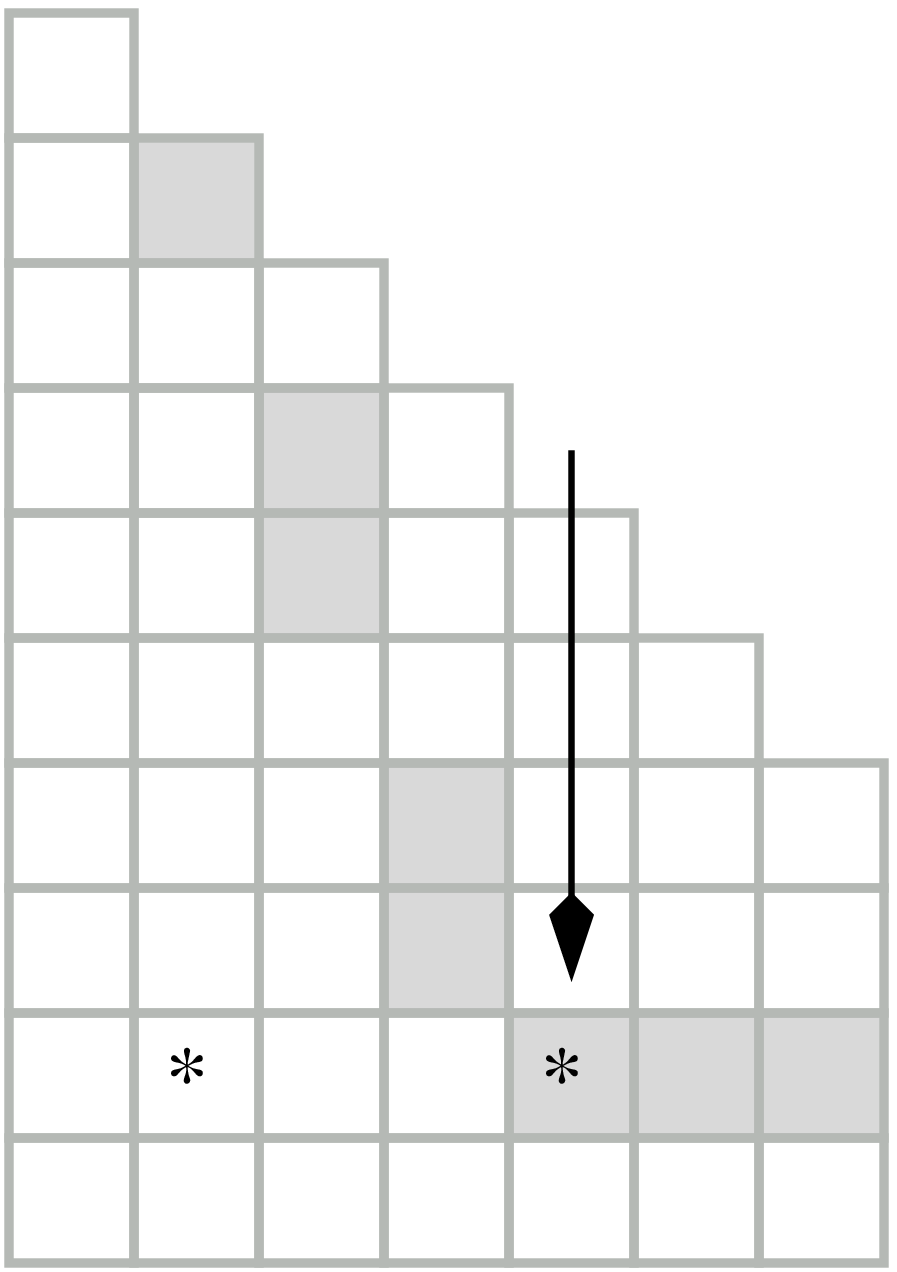}\label{sfl:finish3}}
    \caption[]{Illustration of the construction of $\tilde{x}$,
      Steps~\eqref{step:smallestEntry} to~\eqref{step:adjustBasis}.}
    \label{fig:finish}
  \end{figure}

  At the beginning, we set $\tilde{x} = x^\star$, and let $\lambda > 0$ be
  an arbitrary positive number. Then we perform the following four steps
  (see Figure~\ref{fig:finish} for illustrations of the first three).
  \begin{myenumerate}
  \item\label{step:smallestEntry}%
    For every $(i,j)$ that is the first nonzero element of the shifted
    column of at least one SCI in $\basis^{\star}$, we reduce
    $\tilde{x}_{ij}$ by~$\lambda$.
  \item\label{step:largestEntryValue1}%
    For every $(i,j)$ that is the last nonzero element of the shifted
    column~$S$ of at least one SCI in~$\basis^{\star}$ with $x^{\star}(S) =
    1$, we increase~$\tilde{x}_{ij}$ by $\lambda$.
  \item\label{step:adjustBasis}%
    For each $i \in \ints{p}$ and for all $j = \min\{i,q\}, \min\{i,q\} -
    1, \dots, 1$ (in this order): If $(i,j)$ is the leader of some SCI
    in~$\basis^{\star}$, we adjust~$\tilde{x}_{ij}$ such that, with $B =
    \{(i,j), (i,j+1), \dots, (i,\min\{i,q\})\}$,
    \[
    \tilde{x}(B) =
    \begin{cases}
      1 & \text{if } x^{\star}(B) = 1 \\
      x^{\star}(B) - \lambda & \text{otherwise}
    \end{cases}
    \]
    holds.
  \item\label{step:correctRowSum}%
    For each $i \in \ints{p}$, adjust $\tilde{x}_{ij}$ in order to achieve
    $\tilde{x}(\row{i}) = 1$, where $j =
    \min\setdef{\ell}{x^{\star}_{i\ell} > 0}$.
  \end{myenumerate}

  The reason for treating the case $x^\star(S) = 1$ separately in
  Step~\ref{step:largestEntryValue1} will become evident in the proof of
  Claim~\ref{claim:step4} below.

  The following four claims will yield that~$\tilde{x}$ is a solution of
  the equation system corresponding to~$\basis^\star$.

  \begin{claim}\label{claim:step2}
    After Step~\ref{step:largestEntryValue1}, for each shifted column~$S$
    of some SCI in~$\basis^{\star}$ we have
    \[
    \tilde{x}(S) =
    \begin{cases}
      1 & \text{if } x^{\star}(S) = 1 \\
      x^\star(S) - \lambda & \text{otherwise}.
    \end{cases}
    \]
  \end{claim}

  \begin{proof}
    Let~$S$ be the shifted column of some SCI in~$\basis^\star$. It follows
    from Part~\eqref{claim:firstlast:firstnonzero} of
    Claim~\ref{claim:firstlast} that the first nonzero element $(i,j)$
    of~$S$ is the only element in~$S$ whose $\tilde{x}$-component is
    changed (reduced by~$\lambda$) in Step~\ref{step:smallestEntry}. Thus,
    after Step~\ref{step:smallestEntry} we have $\tilde{x}(S) = x^\star(S)
    - \lambda$.

    If $x^{\star}(S) < 1$, then, by
    Part~\eqref{claim:firstlast:lastnonzero} of
    Claim~\ref{claim:firstlast}, $\tilde{x}(S)$ is not changed in
    Step~\ref{step:largestEntryValue1}.  Otherwise, $x^{\star}(S) = 1$, and
    $\tilde{x}_{k\ell}$ is increased by~$\lambda$ in
    Step~\ref{step:largestEntryValue1}, where~$(k,\ell)$ is the last
    nonzero element of~$S$. According to
    Part~\eqref{claim:firstlast:lastnonzero} of
    Claim~\ref{claim:firstlast}, no other component of~$\tilde{x}$
    belonging to some element in~$S$ is changed in
    Step~\ref{step:largestEntryValue1}. Thus, in both cases the claim
    holds.
  \end{proof}

  \begin{claim}\label{claim:step3a}
    No component of~$\tilde{x}$ belonging to the shifted column of some SCI
    in~$\basis^{\star}$ is changed in Step~\ref{step:adjustBasis}.
  \end{claim}
  \begin{proof}
    Let~$S$ be the shifted column of some SCI in~$\basis^\star$. According
    to Part~\eqref{reduced:containLeader} of
    Claim~\ref{claim:BasicStructure}, $S$ does not contain the leader of
    any SCI in~$\basis^\star$, since~$\basis^{\star}$ is a reduced basis of
    minimal weight.
  \end{proof}

  \begin{claim}\label{claim:step3b}
    After Step~\ref{step:adjustBasis}, for each SCI in~$\basis^{\star}$
    with shifted column~$S$ and bar~$B$ we have $\tilde{x}(S) =
    \tilde{x}(B)$.
  \end{claim}
  \begin{proof}
    For an SCI in~$\basis^\star$ with shifted column~$S$ and bar~$B$, we
    have $x^{\star}(S) = x^{\star}(B)$. Thus, from Claims~\ref{claim:step2}
    and~\ref{claim:step3a} it follows that $\tilde{x}(S) = \tilde{x}(B)$
    holds after Step~\ref{step:adjustBasis}.
  \end{proof}

  \begin{claim}\label{claim:step4}
    Step~\ref{step:correctRowSum} does not change any component
    of~$\tilde{x}$ that belongs to the shifted column or the bar of some
    SCI in~$\basis^{\star}$.
  \end{claim}
  \begin{proof}
    Let $(i,j)$ be such that $x^{\star}_{i\ell} = 0$ for all $\ell < j$ and
    $x^{\star}_{ij} > 0$. By Part~\eqref{reduced:firstRowEntry} of
    Claim~\ref{claim:BasicStructure}, $(i,j)$ is not contained in any
    shifted column of an SCI in~$\basis^{\star}$. If $(i,j)$ is contained
    in the bar~$B$ of some SCI in~$\basis^{\star}$, then clearly
    $x^{\star}(B) = 1$ holds. Thus, after Step~\ref{step:adjustBasis}, we
    have $\tilde{x}(\row{i}) = \tilde{x}(B) = 1$, which shows that
    $\tilde{x}_{ij}$ is not changed in Step~\ref{step:correctRowSum}.
  \end{proof}

  We can now finish the proof of the proposition. Claims~\ref{claim:step3b}
  and~\ref{claim:step4} show that~$\tilde{x}$ satisfies all SCIs contained
  in~$\basis^{\star}$ with equality. Furthermore, in all steps of the
  procedure only components $\tilde{x}_{ij}$ with $x^{\star}_{ij} > 0$ are
  changed (this is clear for Steps~\ref{step:smallestEntry},
  \ref{step:largestEntryValue1}, and~\ref{step:correctRowSum}; for
  Step~\ref{step:adjustBasis} it follows from
  Part~\eqref{reduced:nonzeroLeader} of Claim~\ref{claim:BasicStructure}).
  Since after Step~\ref{step:correctRowSum}, $\tilde{x}$ satisfies all
  row-sum equations, this proves that $\tilde{x}$ is a solution to the
  equation system obtained from~$\basis^{\star}$.

  We assumed that~$\basis^{\star}$ contains at least one SCI. Let~$S$ be
  the shifted column of one of these. We know $x^{\star}(S) > 0$ by
  Claim~\ref{claim:TrivialSCIs}. Thus, let $(i,j)$ be the first nonzero
  element of~$S$. Hence, after Step~\ref{step:smallestEntry}, we have
  $\tilde{x}_{ij} = x^{\star}_{ij} - \lambda$. By
  Part~\eqref{claim:firstlast:both} of Claim~\ref{claim:firstlast}, this
  still holds after Step~\ref{step:largestEntryValue1}. As $\tilde{x}_{ij}$
  is also not changed in Steps~\ref{step:adjustBasis}
  and~\ref{step:correctRowSum} (see Claims~\ref{claim:step3a}
  and~\ref{claim:step4}), we deduce $\tilde{x} \neq x^{\star}$,
  contradicting the fact that~$x^{\star}$ is the unique solution to the
  equation system belonging to~$\basis^{\star}$.

  This concludes the proof of Proposition~\ref{prop:orbipart:descr}.
\end{proof}

We hope that reading this proof was somewhat enjoyable. Anyway, at least it
also gives us a linear description of the packing orbitopes for symmetric
groups almost for free.

\begin{proposition}
  \label{prop:orbipack:descr}
  The packing orbitope $\orbipack{p}{q}$ is completely described by the
  nonnegativity constraints, the row-sum inequalities, and the shifted
  column inequalities:
  \begin{align*}
    \orbipack{p}{q} =
    \{\, & x \in \R^{\orbipartinds{p}{q}} \suchthat
    x \geq \zerovec,\; x(\row{i}) \leq 1 \text{ for } i = 1, \dots, p,\\
    & x(B) - x(S) \leq 0 \text{ for all SCIs with SC } S \text{ and bar }B\, \}.
  \end{align*}
\end{proposition}

\begin{proof}
  Let~$Q \subset \R^{\orbipartinds{p}{q}}$ be the polyhedron on the
  right-hand side of the statement. We define $\mathcal{A} := \setdef{x \in
    \R^{\orbipartinds{p+1}{q+1}}}{x(\row{i}) = 1 \text{ for all }i \in
    \ints{p+1}}$.

  The proof of Proposition~\ref{prop:orbipart:descr} in fact shows that its
  statement remains true if we drop all SCIs with shifted column~$S$ and $S
  \cap \col{1} \neq \varnothing$ from the linear description. This follows
  from the fact that, due to $x^{\star}_{11} = 1$ and
  Claim~\ref{claim:TrivialSCIs}, no such SCI can be contained in any
  reduced basis of~$x^{\star}$ (using the notations from the proof of
  Proposition~\ref{prop:orbipart:descr}). Thus we obtain
  \begin{equation}\label{eq:orbipack:descr:2}
    \orbipart{p+1}{q+1} = \mathcal{A} \cap \tilde{Q},
  \end{equation}
  with
  \begin{align*}
    \tilde{Q} =
    \{ x \in \R^{\orbipartinds{p+1}{q+1}} \suchthat
    & x(B) - x(S) \leq 0 \text{ for all SCIs with bar } B\\
    & \text{and shifted column } S \text{ with } S \cap \col{1} = \varnothing,\\
    & x_{ij} \geq 0 \text{ for all } (i,j) \in \orbipartinds{p+1}{q+1}
    \setminus \col{1},\\
    & x(\row{i} - (i,1)) \leq 1 \text{ for all } i = 2, \dots, p+1 \},
  \end{align*}
  where the last inequalities are equivalent (with respect to
  $\orbipart{p+1}{q+1}$) to the nonnegativity constraints associated with
  the elements of $\col{1}$ by addition of row-sum equations.

  Define $\mathcal{L} := \setdef{x \in \R^{\orbipartinds{p+1}{q+1}}}{x_{i1}
    =0 \text{ for all }i \in \ints{p+1}}$, and denote by $\tilde{\pi} :
  \R^{\orbipartinds{p+1}{q+1}} \rightarrow \mathcal{L}$ the orthogonal
  projection.  Since none of the inequalities defining~$\tilde{Q}$ has a
  nonzero coefficient in~$\col{1}$, we have $\tilde{\pi}^{-1}(\tilde{Q}
  \cap \mathcal{L}) = \tilde{Q}$, hence $\tilde{Q} \cap \mathcal{L} =
  \tilde{\pi}(\tilde{Q})$. This yields $\tilde{\pi}(\mathcal{A} \cap
  \tilde{Q}) = \tilde{\pi}(\mathcal{A}) \cap \tilde{\pi}(\tilde{Q})$,
  which, due to $\tilde{\pi}(\mathcal{A}) = \mathcal{L}$, implies
  $\tilde{\pi}(\mathcal{A} \cap \tilde{Q}) = \tilde{Q} \cap \mathcal{L}$.
  Thus, we obtain
  \[
  \orbipack{p}{q}=\tilde{\pi}(\orbipart{p+1}{q+1})=\tilde{\pi}(\mathcal{A} \cap \tilde{Q})
  =\tilde{Q} \cap \mathcal{L}=Q,
  \]
  where the first equation is due to Proposition~\ref{prop:projpartpack},
  the second equation follows from~\eqref{eq:orbipack:descr:2}, and the final arises
  from identifying~$\mathcal{L}$ with~$\R^{\orbipartinds{p}{q}}$.
\end{proof}

% ------------------------------------------------------
% Facets
% ------------------------------------------------------

\subsection{Facets}
\label{sec:Facets}

In this section, we investigate which of the constraints from the linear
descriptions of $\orbipart{p}{q}$ and $\orbipack{p}{q}$ given in
Propositions~\ref{prop:orbipart:descr} and~\ref{prop:orbipack:descr},
respectively, define facets. This will also yield non-redundant
descriptions.

It seems to be more convenient to settle the packing case first and then to
carry over the results to the partitioning case. Recall that we assume $2
\leq p \leq q$.

\begin{proposition}\label{prop:orbipack:facets}\
  \begin{myenumerate}
  \item\label{prop:orbipack:facets:fulldimensional}%
    The packing orbitope $\orbipack{p}{q} \subset \R^{\orbipartinds{p}{q}}$
    is full dimensional:
    \[
    \dim(\orbipack{p}{q}) = \card{\orbipartinds{p}{q}} =
    pq - \tfrac{q(q-1)}{2} = \big(p-\tfrac{q-1}{2}\big)q.
    \]
  \item A nonnegativity constraint $x_{ij} \geq 0$, $(i,j) \in
    \orbipartinds{p}{q}$, defines a facet of $\orbipack{p}{q}$, unless $i =
    j < q$ holds. The faces defined by $x_{jj} \geq 0$ with $j < q$ are
    contained in the facet defined by~$x_{qq} \geq 0$.%
    \label{prop:orbipack:facets:nonneg}
  \item Every row-sum constraint $x(\row{i}) \leq 1$ for $i \in \ints{p}$
    defines a facet of $\orbipack{p}{q}$.%
    \label{prop:orbipack:facets:rowsum}
  \item\label{prop:orbipack:facets:sci}%
    A shifted column inequality $x(B)-x(S)\le 0$ with bar~$B$ and shifted
    column~$S = \{\diagcol{1}{c_1}, \diagcol{2}{c_2}, \dots,
    \diagcol{\eta}{c_{\eta}}\}$ defines a facet of $\orbipack{p}{q}$,
    unless $\eta \geq 2$ and $c_1 < c_2$ (exception~I) or $\eta=1$ and
    $B\neq\{\diagcol{1}{c_1+1}\}$ (exception~II) hold.  In case of
    exception~I, the corresponding face is contained in the facet defined
    by the SCI with bar~$B$ and shifted column $\{\diagcol{1}{c_2},
    \diagcol{2}{c_2}, \dots, \diagcol{\eta}{c_{\eta}}\}$. In case of
    exception~II, the face is contained in the facet defined by the SCI
    $x_{\diagcol{1}{c_1+1}}-x_{\diagcol{1}{c_1}}\le 0$.
   \end{myenumerate}
\end{proposition}

\begin{figure}
  \centering
  \newcommand{\mystyle}[1]{\footnotesize {#1}}
  \psfrag{k}{\mystyle{$k$}}
  \psfrag{l}{\mystyle{$\ell$}}
  \psfrag{j}{\mystyle{$j$}}
  \psfrag{i}{\mystyle{$i$}}
  \psfrag{B}{\textcolor{white}{\mystyle{$B$}}}
  \psfrag{c1}{\mystyle{$c_1$}}
  \psfrag{c2}{\mystyle{$c_2$}}
  \subfloat[][matrix $V^{k\ell}$]{\includegraphics[height=.17\textheight]{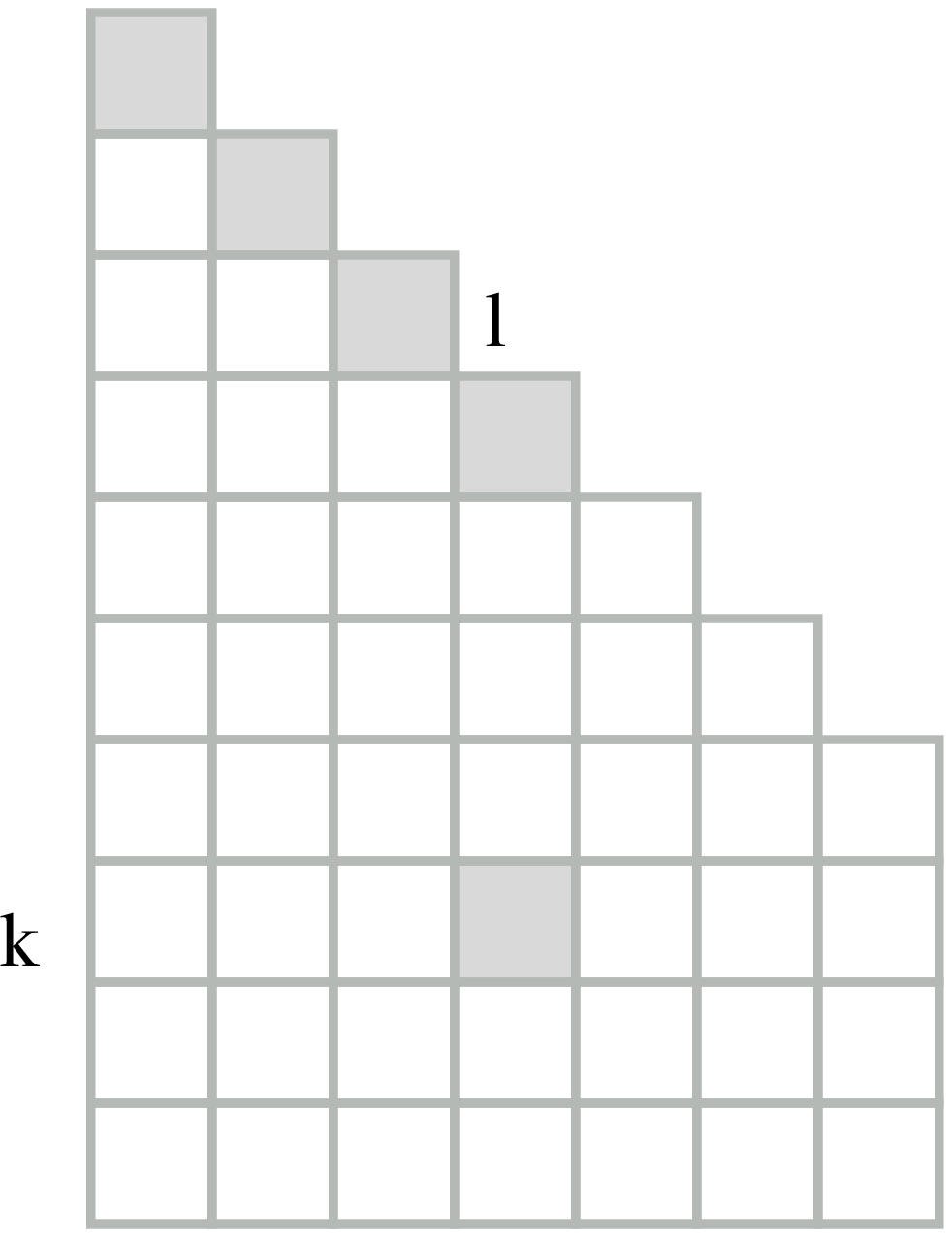}\label{sfl:packnonneg}}\qquad
  \subfloat[][matrix $\hat{V}^{k\ell}$]{\includegraphics[height=.17\textheight]{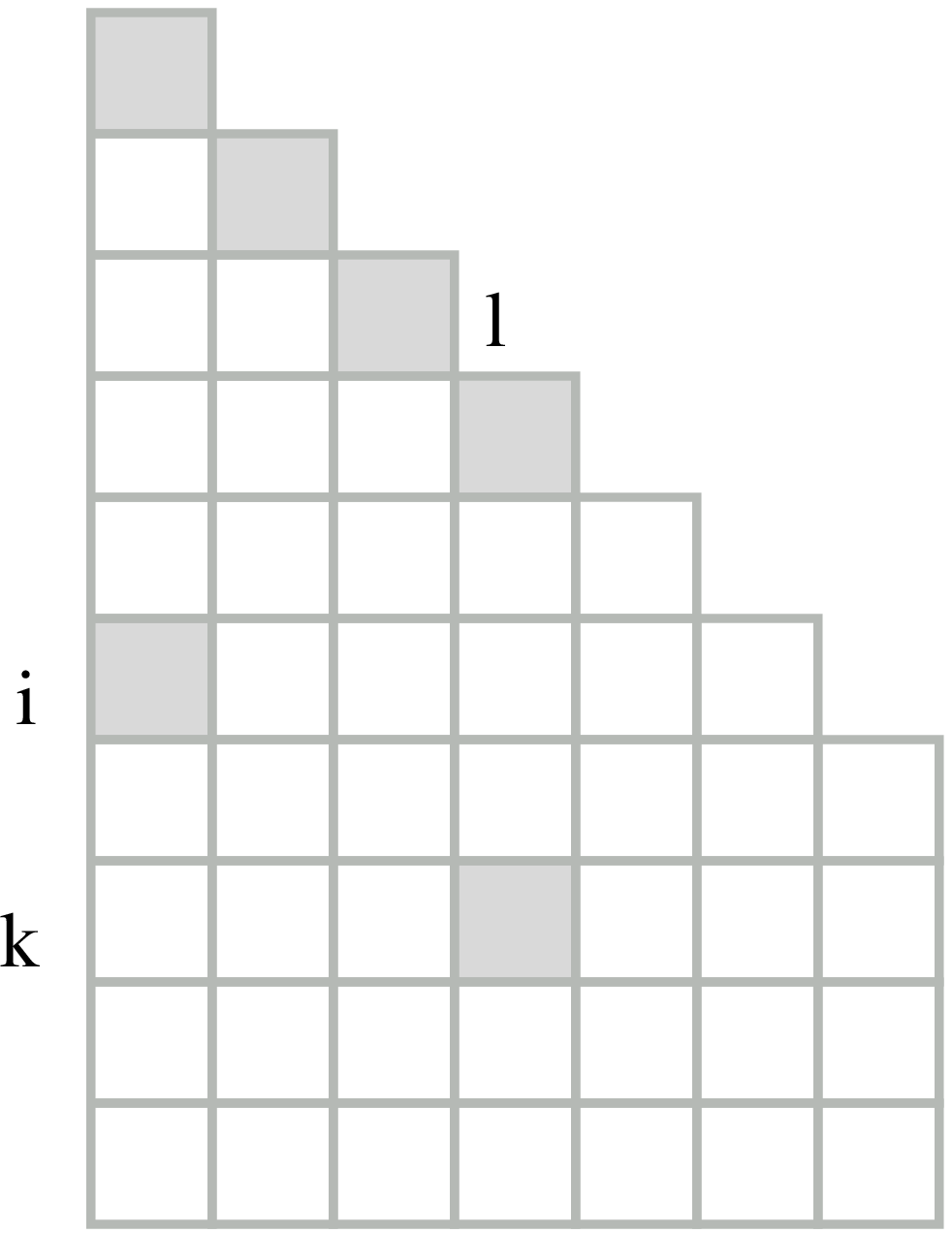}\label{sfl:packrowsum}}\qquad
  \subfloat[][]{\includegraphics[height=.17\textheight]{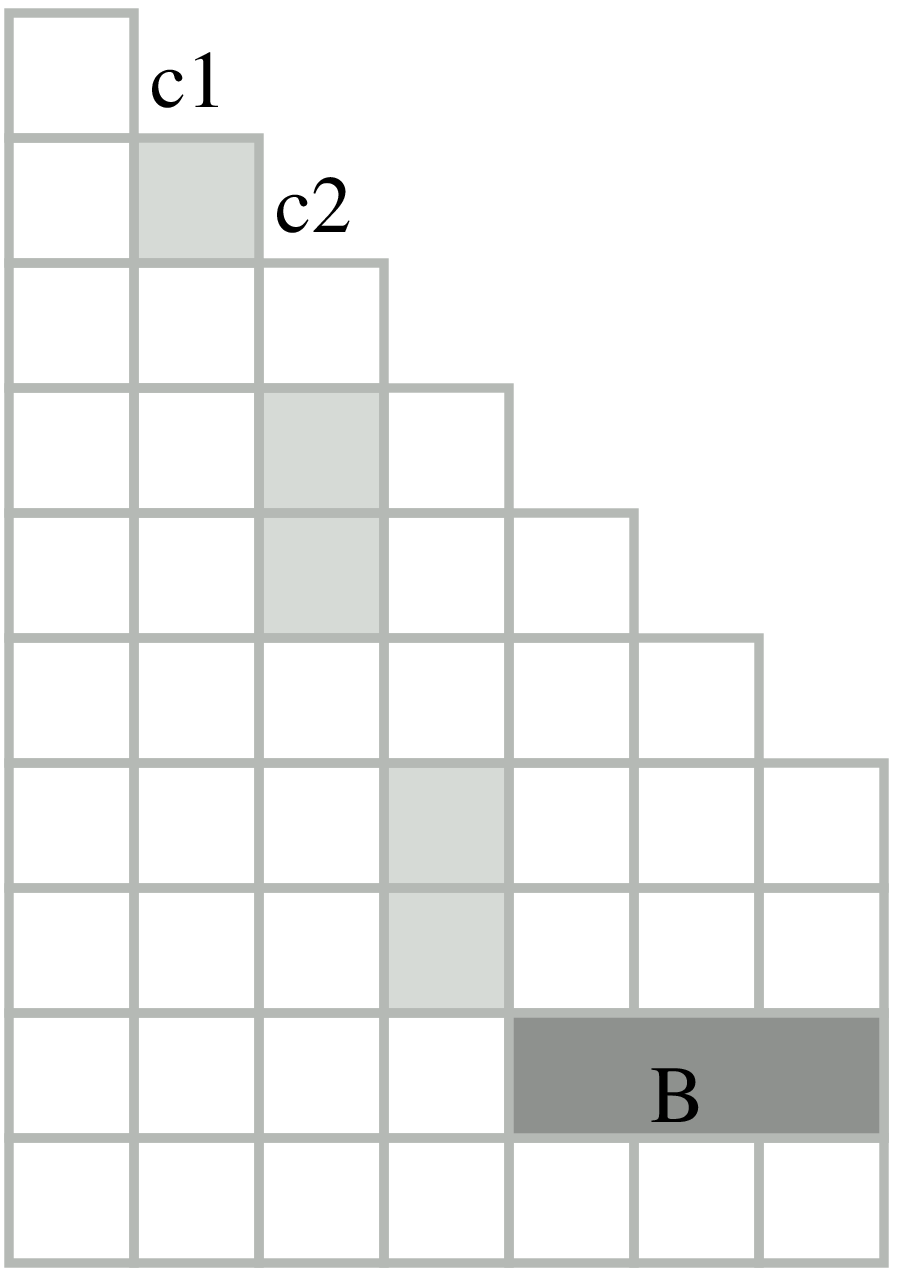}\label{sfl:SCInonFacet}}
  \caption[]{\subref{sfl:packnonneg}--\subref{sfl:packrowsum}: Illustration
    of the matrices used in the proof of
    parts~\eqref{prop:orbipack:facets:fulldimensional}
    and~\eqref{prop:orbipack:facets:rowsum} of
    Proposition~\ref{prop:orbipack:facets}. \subref{sfl:SCInonFacet}:
    Example of an SCI that does not define a facet; see the proof of
    Part~\eqref{prop:orbipack:facets:sci} of
    Proposition~\ref{prop:orbipack:facets}.}
  \label{fig:orbipack:facets}
\end{figure}

\begin{proof}
  \begin{partslist}
  \partsitem For all $(k,\ell) \in \orbipartinds{p}{q}$, we define $V^{k\ell} = (v^{k\ell}_{ij})
    \in \R^{\orbipartinds{p}{q}}$ by
    \[
    v^{k\ell}_{ij} =
    \begin{cases}
      1 & \text{if } \big(i = j \leq \ell \text{ and } j < q\big) \text{ or } (i,j) = (k,\ell)\\
      0 & \text{otherwise}
    \end{cases}
    \quad \text{for }(i,j) \in \orbipartinds{p}{q},
    \]
    that is, $V^{k\ell}$ has $1$-entries at position~$(k,\ell)$ and on the
    main diagonal up to column~$\ell$, except that $v^{k\ell}_{qq} = 0$
    unless $(k,\ell) = (q,q)$; see
    Figure~\ref{fig:orbipack:facets}~\subref{sfl:packnonneg}. The columns
    of each~$V^{k\ell}$ are in non-increasing lexicographic order.  Hence,
    by Part~\eqref{obs:charvert:sym} of Observation~\ref{obs:charvert},
    each $V^{k\ell}$ is a vertex of~$\orbipack{p}{q}$.

    In order to show that these vectors are linearly independent, we fix an
    arbitrary ordering of the $V^{k\ell}$ that starts with $V^{11}, V^{22},
    \dots, V^{q-1,q-1}$. For each $(k,\ell) \in \orbipartinds{p}{q}$, all
    points $V^{rs}$ preceding~$V^{k\ell}$ have a $0$-entry at position
    $(k,\ell)$, while $v^{k\ell}_{k\ell} = 1$. This shows that these
    $\card{\orbipartinds{p}{q}}$ vertices of $\orbipack{p}{q}$ are linearly
    independent.  Together with~$\zerovec$ this
    gives~$\card{\orbipartinds{p}{q}}+1$ affinely independent points
    contained in $\orbipack{p}{q}$, proving that $\orbipack{p}{q}$ is full
    dimensional. The calculations in the statement are straightforward.
  \partsitem For $(i,j) \in \orbipartinds{p}{q} \setminus
    \setdef{(j,j)}{j < q}$ all points $V^{k\ell}$ with $(k,\ell) \neq
    (i,j)$ are contained in the face defined by $x_{ij} \geq 0$. Since this
    is also true for~$\zerovec$, the face defined by $x_{ij} \geq 0$
    contains $\card{\orbipartinds{p}{q}}$ affinely independent points (see
    the proof of Part~\eqref{prop:orbipack:facets:fulldimensional}), i.e.,
    it is a facet of $\orbipack{p}{q}$.

    For every vertex~$x^\star \in \orbipack{p}{q}$ contained in the face
    defined by $x_{jj} \geq 0$ for some $j < q$, we have
    $x^\star_{\ell\ell} = 0$ for all $\ell \geq j$ (because otherwise the
    columns of~$x^{\star}$ would not be in non-increasing lexicographic order).
    This shows that~$x^\star$ is contained in the facet
    defined by~$x_{qq} \geq 0$.

 \partsitem In order to show that $x(\row{i}) \leq 1$ defines a facet of
    $\orbipack{p}{q}$ for $i \in \ints{p}$, we construct points
    $\hat{V}^{k\ell}$ (depending on~$i$) from the points $V^{k\ell}$
    defined in Part~\eqref{prop:orbipack:facets:fulldimensional} by adding
    a $1$ at position $(i,1)$ if $V^{k\ell}(\row{i}) = 0$ (see
    Figure~\ref{fig:orbipack:facets}~\subref{sfl:packrowsum}). The
    $(\card{\orbipartinds{p}{q}}-1)$ points~$\hat{V}^{k\ell}$ for all
    $(k,\ell) \in \orbipartinds{p}{q} - (i,1)$, and the unit
    vector~$E^{i1}$ (with a single~$1$ in position~$(i,1)$) satisfy
    $x(\row{i}) = 1$.  Furthermore, they are affinely independent, since
    subtracting $E^{i1}$ from all vectors~$\hat{V}^{k\ell}$ yields
    vectors~$\tilde{V}^{k\ell}$, which can be shown to be linearly independent similarly to Part~\eqref{prop:orbipack:facets:fulldimensional}; here, we need $(k,\ell)\neq(i,1)$.

  \partsitem Let $x(B) - x(S) \leq 0$ be an SCI with bar~$B$, leader $(i,j) =
    \diagcol{\eta}{j}$, and shifted column $S = \{\diagcol{1}{c_1},
    \diagcol{2}{c_2}, \dots, \diagcol{\eta}{c_{\eta}}\}$.

    If $\eta \geq 2$ and $c_1 < c_2$ hold (exception~I), then the SCI is the sum of the SCI
    \[
    x_{\diagcol{1}{c_1+1}} - x_{\diagcol{1}{c_1}} \leq 0
    \]
    and the SCI with bar~$B$ and shifted column $\{\diagcol{1}{c_1+1},
    \diagcol{2}{c_2}, \dots, \diagcol{\eta}{c_{\eta}}\}$; see
    Figure~\ref{fig:orbipack:facets}~\subref{sfl:SCInonFacet}. Repeating
    this argument $(c_2 - c_1 - 1)$ times proves the second statement of
    Part~\eqref{prop:orbipack:facets:sci} for exception~I.

    If $\eta=1$ and $B=\{\diagcol{1}{j}\}$ with $j> c_1+1$ hold
    (exception~II), then the SCI is the sum of the SCIs
    $x_{\diagcol{1}{c_1+1}}-x_{\diagcol{1}{c_1}}\le 0$, \dots,
    $x_{\diagcol{1}{j}}-x_{\diagcol{1}{j-1}}\le 0$. This proves the second
    statement of Part~\eqref{prop:orbipack:facets:sci} for exception~II.

    Otherwise, let $\mathcal{V}$ be the set of vertices
    of~$\orbipack{p}{q}$ that satisfy the SCI with equality, and let
    $\mathcal{L} = \lin{\mathcal{V} \cup \{E^{ij}\}}$ be the linear span
    of~$\mathcal{V}$ and the unit vector~$E^{ij}$. We will show that
    $\mathcal{L} = \R^{\orbipartinds{p}{q}}$, which proves
    $\dim(\aff{\mathcal{V}}) = \card{\orbipartinds{p}{q}} - 1$ (since
    $\zerovec \in \mathcal{V}$).  Hence, the SCI defines a facet
    of~$\orbipack{p}{q}$.

    \begin{figure}
      \centering
      \newcommand{\mystyle}[1]{\footnotesize {#1}}
      \psfrag{i}{\mystyle{$i$}}
      \psfrag{j}{\mystyle{$j$}}
      \psfrag{s}{\mystyle{$s$}}
      \psfrag{r}[r][r]{\mystyle{$r$}}
      \psfrag{A}{\mystyle{$A$}}
      \psfrag{B}{\textcolor{white}{\mystyle{$B$}}}
      \psfrag{C}{\mystyle{$C$}}
      \psfrag{D}{\mystyle{$D$}}
      \psfrag{S}{\mystyle{$S$}}
      \subfloat[][All cases]{\includegraphics[height=.17\textheight]{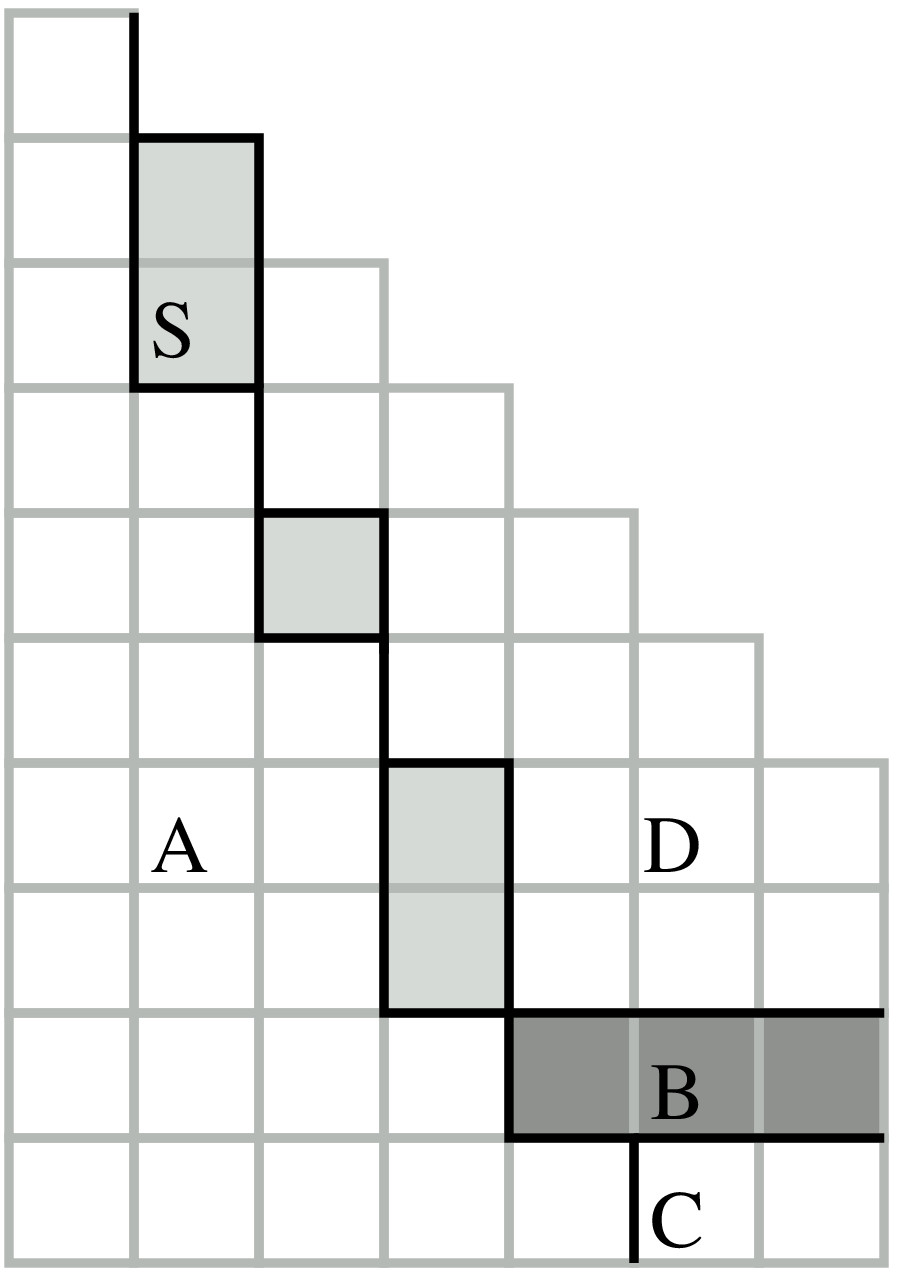}\label{sfl:packsci}}\hfill
      \subfloat[][Case A, $W^{rs}$]%
      {\includegraphics[height=.17\textheight]{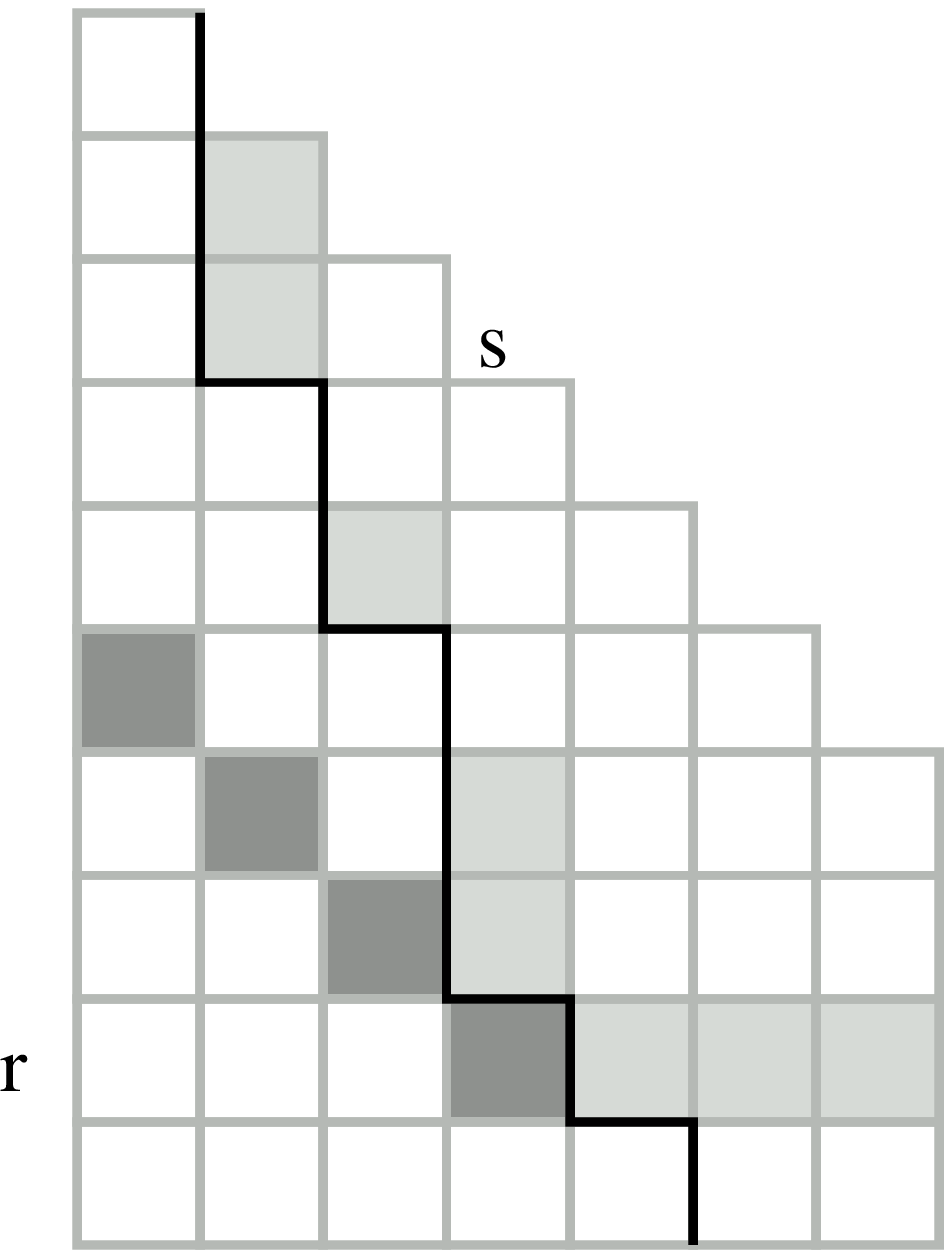}\label{sfl:packsci2}}\hfill
      \subfloat[][Case D,\,$U^{rs}$]%
      {\includegraphics[height=.17\textheight]{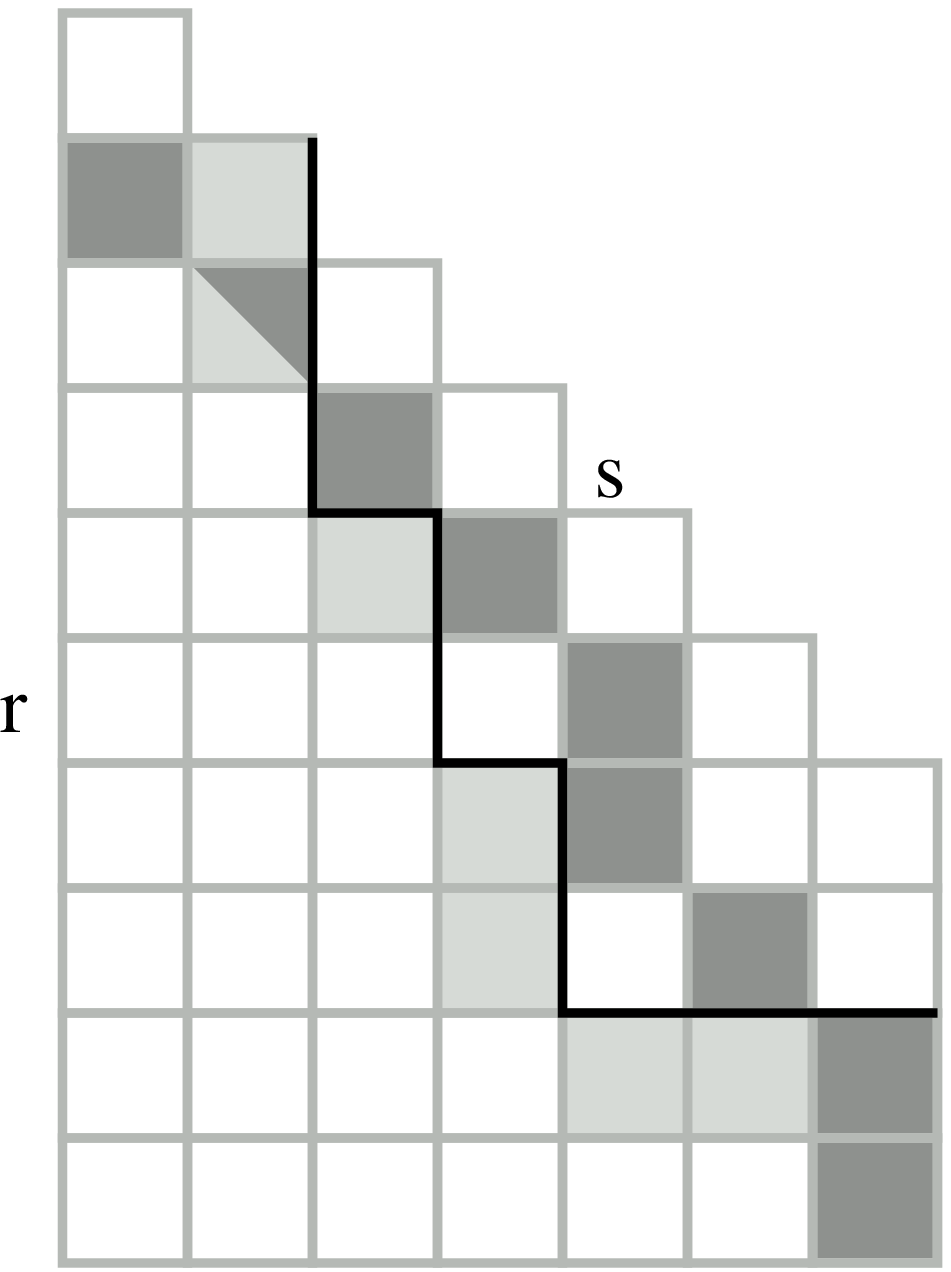}\label{sfl:packsci3}}\hfill
      \subfloat[][Case S]{\includegraphics[height=.17\textheight]{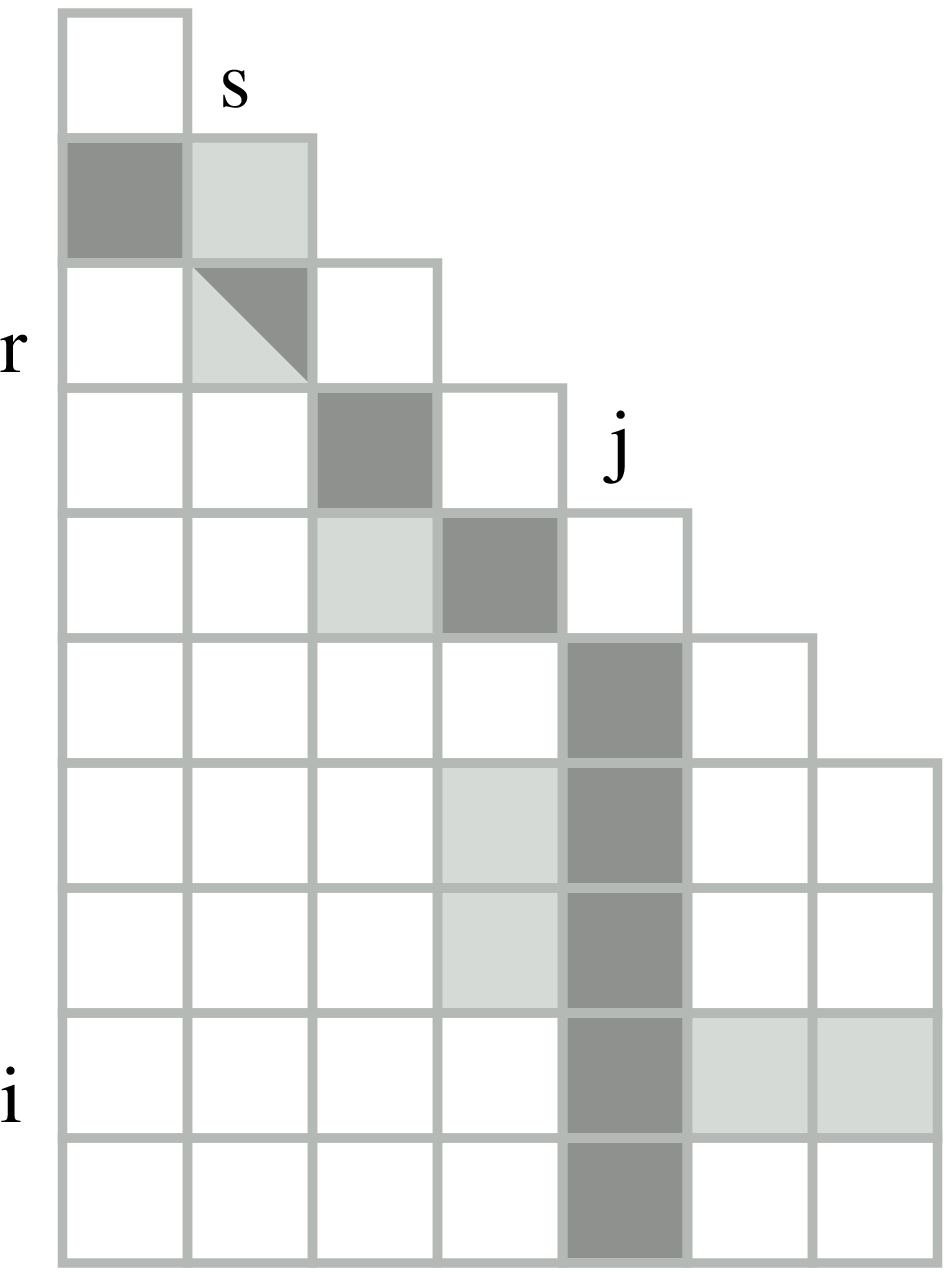}\label{sfl:packsci4}}
      \caption[]{Illustration of the constructions in the proof of
        Part~\eqref{prop:orbipack:facets:sci} of
        Proposition~\ref{prop:orbipack:facets}.}
      \label{fig:orbipack:facets:sci}
    \end{figure}

    To show that $\mathcal{L} = \R^{\orbipartinds{p}{q}}$, we prove that
    $E^{rs} \in \mathcal{L}$ for all $(r,s) \in \orbipartinds{p}{q}$.  We
    partition the set $\orbipartinds{p}{q}\setminus(B\cup S)$ into three
    parts (see Figure~\ref{fig:orbipack:facets:sci}~\subref{sfl:packsci}):
    \begin{align*}
      A := & \setdef{\diagcol{\rho}{s} \in \orbipartinds{p}{q}}%
      {(\rho \leq \eta \text{ and } s < c_{\rho})\text{ or } \rho > \eta},\\
      C := & \setdef{\diagcol{\rho}{s} = (r,s) \in \orbipartinds{p}{q}}%
      {\rho \leq \eta \text{ and } r > i},\text{ and}\\
      D := & \setdef{\diagcol{\rho}{s} = (r,s) \in \orbipartinds{p}{q}}%
      {\rho < \eta,\; s > c_{\rho}, \text{ and } r < i}.
    \end{align*}

    For $(r,s) = \diagcol{\rho}{s}$, denote by $\diagleq{r}{s} =
    \{\diagcol{\rho}{1}, \diagcol{\rho}{2}, \dots, \diagcol{\rho}{s}\}$ the
    diagonal starting at $\diagcol{\rho}{1} = (r-s+1,1)$ and ending at
    $\diagcol{\rho}{s} = (r,s)$.  Similarly, denote by $\diaggeq{r}{s} =
    \{\diagcol{\rho}{s}, \diagcol{\rho}{s+1}, \dots\} \cap
    \orbipartinds{p}{q}$ the diagonal starting at $(r,s)$ and ending in
    $\col{q}$ or in $\row{p}$.

    \begin{claim}\label{claim:facets:AC}
      For all $(r,s) = \diagcol{\rho}{s} \in A \cup C$ we have $E^{rs} \in
      \mathcal{L}$.
    \end{claim}
    \begin{proof}
      Denote the incidence vector of $\diagleq{r}{s}$ by $W^{rs} =
      \pointof{\diagleq{r}{s}}$ (see
      Figure~\ref{fig:orbipack:facets:sci}~\subref{sfl:packsci2}).  Both
      $W^{rs}$ and~$W^{rs}-E^{rs}$ are vertices of $\orbipack{p}{q}$. We
      have $\diagleq{r}{s} \cap (B \cup S) = \varnothing$ for $(r,s) \in
      A$.  Furthermore
      \[
      \card{\diagleq{r}{s} \cap B} = 1 = \card{\diagleq{r}{s} \cap S}
      \]
      for $(r,s) \in C$. Hence, these two vertices satisfy the SCI with
      equality and we obtain $E^{rs} =W^{rs} - (W^{rs} - E^{rs}) \in
      \mathcal{L}$.
    \end{proof}

    \begin{claim}\label{claim:facets:D}
      For all $(r,s) = \diagcol{\rho}{s} \in D$ we have $E^{rs} \in
      \mathcal{L}$.
    \end{claim}
    \begin{proof}
      Define the set
      \[
      U(r,s) := \diagleq{r}{s} \cup \diaggeq{r+1}{s} \cup
      \big(\{\diagcol{\rho+1}{q}, \diagcol{\rho+2}{q},\dots\} \cap \orbipartinds{p}{q}\big),
      \]
      see Figure~\ref{fig:orbipack:facets:sci}~\subref{sfl:packsci3}. Let
      $U^{rs} := \pointof{U(r,s)}$. By construction, the three points
      $U^{rs}$, $U^{rs} - E^{rs}$, and $U^{rs} - E^{r+1,s}$ are vertices of
      $\orbipack{p}{q}$.

      If $\rho = 1$, we have $\card{U(r,s) \cap B} = 1$ and $\card{U(r,s)
        \cap S} = 1$, where we need $c_1 = c_2$ in case of $s = c_1 + 1$
      (notice that in case of $\eta=1$ we have $D=\varnothing$).  Due to
      $(r,s) \notin B \cup S$, both $U^{rs}$ and $U^{rs} - E^{rs}$ satisfy
      the SCI with equality. This yields $E^{rs} = U^{rs} - (U^{rs} -
      E^{rs}) \in \mathcal{L}$.

      If $\rho > 1$, then $\card{U(r,s) \cap S} = 1$ does not hold in all
      cases (e.g., if $s = c_{\rho+1}$, we have $(r+1,s) \in S$). However,
      since $\rho > 1$, $U(r-1,s)$ is well-defined and
      \[
      \card{U(r-1,s) \cap B} = 1
      \qquad\text{and}\qquad
      \card{U(r-1,s) \cap S} = 1
      \]
      hold. Hence the vertices $U^{r-1,s}$ and $U^{r-1,s} - E^{rs}$ satisfy
      the SCI with equality, giving $E^{rs} = U^{r-1,s} - (U^{r-1,s} -
      E^{rs}) \in \mathcal{L}$.
    \end{proof}

    \begin{claim}\label{claim:facets:S}
      For all $(r,s) = \diagcol{\rho}{s} \in S$ we have $E^{rs} \in
      \mathcal{L}$.
    \end{claim}
    \begin{proof}
      Define the set
      \[
      T(r,s) := \diagleq{r+j-s}{j} \cup \big(\{\diagcol{\rho+1}{j},
      \diagcol{\rho+2}{j}, \dots\} \cap \orbipartinds{p}{q}\big),
      \]
      see Figure~\ref{fig:orbipack:facets:sci}~\subref{sfl:packsci4}. The
      incidence vector $T^{rs} := \pointof{T(r,s)}$ is a vertex of
      $\orbipack{p}{q}$, which, due to $T(r,s) \cap S = \{(r,s)\}$ and
      $T(r,s) \cap B = \{(i,j)\}$ satisfies the SCI with equality. Thus,
      from
      \[
      E^{rs} = T^{rs} - E^{ij}
      - \sum_{(k,\ell) \in T(r,s) \cap A} E^{k\ell}
      - \sum_{(k,\ell) \in T(r,s) \cap C} E^{k\ell}
      - \sum_{(k,\ell) \in T(r,s) \cap D} E^{k\ell}
      \]
      we conclude $E^{rs} \in \mathcal{L}$, since $E^{ij} \in \mathcal{L}$ by
      definition of~$\mathcal{L}$, and $E^{k\ell} \in \mathcal{L}$ for all
      $(k,\ell) \in A \cup C \cup D$ by Claims~\ref{claim:facets:AC}
      and~\ref{claim:facets:D}.
    \end{proof}

    \begin{claim}\label{claim:facets:B}
      For all $(i,s) = \diagcol{\rho}{s} \in B$ we have $E^{rs} \in
      \mathcal{L}$.
    \end{claim}
    \begin{proof}
      The vector $W^{is} := \pointof{\diagleq{i}{s}}$ is a vertex of
      $\orbipack{p}{q}$ that satisfies the SCI with equality. Furthermore, we
      have
      \[
      E^{is} = W^{is} - E^{rc_{\rho}}
      - \sum_{(k,\ell) \in \diagleq{i}{s} \cap A} E^{k\ell}
      - \sum_{(k,\ell) \in \diagleq{i}{s} \cap D} E^{k\ell},
      \]
      where $(r,c_{\rho}) := \diagcol{\rho}{c_{\rho}} \in S$. Thus, we
      conclude $E^{is}\in\mathcal{L}$, since $E^{k\ell}\in\mathcal{L}$ for all
      $(k,\ell)\in A\cup D \cup S$ by Claims~\ref{claim:facets:AC},
      \ref{claim:facets:D}, and \ref{claim:facets:S}.
    \end{proof}

    Claims~\ref{claim:facets:AC} to~\ref{claim:facets:B} show $E^{rs} \in
    \mathcal{L}$ for all $(r,s) \in \orbipartinds{p}{q}$. This proves that
    the SCI defines a facet of $\orbipack{p}{q}$ (unless exception~I or~II
    hold).
  \end{partslist}
\end{proof}

Finally, we carry the results of Proposition~\ref{prop:orbipack:facets}
over to partitioning orbitopes.

\begin{proposition}\label{prop:orbipart:facets}\
  \begin{myenumerate}
  \item The partitioning orbitope $\orbipart{p}{q} \subset
    \R^{\orbipartinds{p}{q}}$ has dimension
    \[
    \dim(\orbipart{p}{q}) = \card{\orbipartinds{p-1}{q-1}} =
    \card{\orbipartinds{p}{q}} - p =
    \big( p - \tfrac{q}{2}\big)(q-1).
    \]
    The constraints $x(\row{i}) = 1$ form a complete and non-redundant
    linear description of $\aff{\orbipart{p}{q}}$.
  \item\label{prop:orbipart:facets:nonnegative}%
    A nonnegativity constraint $x_{ij} \geq 0$, $(i,j) \in
    \orbipartinds{p}{q}$, defines a facet of $\orbipart{p}{q}$, unless $i =
    j < q$ holds. The faces defined by $x_{jj} \geq 0$ with $j < q$ are
    contained in the facet defined by~$x_{qq} \geq 0$.
  \item\label{prop:orbipart:facets:sci}%
    A shifted column inequality $x(B) - x(S) \leq 0$ with bar~$B$ and
    shifted column~$S = \{\diagcol{1}{c_1}, \diagcol{2}{c_2}, \dots,
    \diagcol{\eta}{c_{\eta}}\}$ defines a facet of $\orbipart{p}{q}$,
    unless $c_1=1$ (Exception~I) or $\eta \geq 2$ and $c_1 < c_2$
    (Exception~II) or $\eta=1$ and $B\neq\{\diagcol{1}{c_1+1}\}$
    (Exception~III).  In case of Exception~I, the corresponding face is
    contained in the facet defined by $x_{i1}\ge 0$, where~$i$ is the index
    of the row containing~$B$.  In case of Exception~II, the face is
    contained in the facet defined by the SCI with bar~$B$ and shifted
    column $\{\diagcol{1}{c_2}, \diagcol{2}{c_2}, \dots,
    \diagcol{\eta}{c_{\eta}}\}$. In case of Exception~III, the face is
    contained in the facet defined by the SCI
    $x_{\diagcol{1}{c_1+1}}-x_{\diagcol{1}{c_1}}\le 0$.
 \end{myenumerate}
\end{proposition}

\begin{proof}
  According to Proposition~\ref{prop:projpartpack}, $\orbipack{p-1}{q-1}$
  is isomorphic to $\orbipart{p}{q}$ via the orthogonal projection of the
  latter polytope to the space
  \[
  \mathcal{L} := \setdef{x \in \R^{\orbipartinds{p}{q}}}{x_{i1} = 0 \text{
      for all } i \in \ints{p}}
  \]
  (and via the canonical identification of~$\mathcal{L}$ and
  $\R^{\orbipartinds{p-1}{q-1}}$). This shows the statement on the
  dimension of $\orbipart{p}{q}$; the calculations and the claim on the
  non-redundancy of the equation system are straightforward.

  Furthermore, this projection (which is one-to-one on
  $\aff{\orbipart{p}{q}}$) maps every face of $\orbipart{p}{q}$ that is
  defined by some inequality
  \[
  \scal{a}{x} := \sum_{(i,j) \in \orbipartinds{p}{q}} a_{ij}\,x_{ij} \leq
  a_0,
  \]
  with $a \in \R^{\orbipartinds{p}{q}}$, $a_0 \in \R$, and $a_{i1} = 0$
  for all $i \in \ints{p}$ to a face of $\orbipack{p-1}{q-1}$ of the same
  dimension defined by
  \[
  \sum_{(i,j) \in \orbipartinds{p-1}{q-1}} a_{i+1,j+1}\, x_{ij} \leq a_0.
  \]
  Conversely, if $\scal{\tilde{a}}{x} \leq \tilde{a}_0$ defines a face of
  $\orbipack{p-1}{q-1}$ for $\tilde{a}\in\R^{\orbipartinds{p-1}{q-1}}$ and
  $\tilde{a}_0\in\R$, then the inequality \[ \sum_{(i,j) \in
    \orbipartinds{p}{q}} \tilde{a}_{ij}\, x_{i+1,j+1} \leq \tilde{a}_0
  \]
  defines a face of $\orbipart{p}{q}$ of the same dimension.

  Due to parts~\eqref{prop:orbipack:facets:nonneg}
  and~\eqref{prop:orbipack:facets:rowsum} of
  Proposition~\ref{prop:orbipack:facets}, this proves
  Part~\eqref{prop:orbipart:facets:nonnegative} of the proposition, where
  we use the fact that the inequalities $x_{i1} \geq 0$ are equivalent to
  $x\big(\row{i} - (i,1)\big) \leq 1$ with respect to $\orbipart{p}{q}$.

  Furthermore, due to Part~\eqref{prop:orbipack:facets:sci} of
  Proposition~\ref{prop:orbipack:facets}, the above arguments also imply
  the statements of Part~\eqref{prop:orbipart:facets:sci} for $c_1\geq 2$
  (including Exception~II and~III).  Finally, we consider the case $c_1=1$
  (Exception~I). Since we have $x_{1,1}=1$ for all $x\in\orbipart{p}{q}$,
  the equation $x(B)-x(S)=0$ implies
  \[
  1\geq x(B)=x(S)\ge x_{1,1}= 1,
  \]
  and hence $x_{i,1}=0$ (using the row-sum equation for row~$i$ containing~$B$).
   This concludes the proof.
\end{proof}

% ------------------------------------------------------
% Summary of Results on the Symmetric Group
% ------------------------------------------------------

\subsection{Summary of Results on the Symmetric Group}
\label{SummaryOfResultsSymmetricGroup}

We collect the results on the packing- and partitioning orbitopes for
symmetric groups.

\begin{theorem}\label{thm:resultsPart}
  The partitioning orbitope $\orbipartG{p}{q}{\symgr{q}}$ $($for $p \geq q
  \geq 2)$ with respect to the symmetric group $\symgr{q}$ equals the set
  of all $x \in \R^{\orbipartinds{p}{q}}$ that satisfy the following linear
  constraints:
  \begin{myitemize}
  \item the row-sum equations $x(\row{i}) = 1$ for all $i \in \ints{p}$,
  \item the nonnegativity constraints $x_{ij} \geq 0$ for all $(i,j) \in
    \orbipartinds{p}{q} \setminus \setdef{(j,j)}{j < q}$,
  \item the shifted column inequalities $x(B) - x(S) \leq 0$ for all bars
    \[
    B = \{(i,j), (i,j+1), \dots, (i,\min\{i,q\})\}
    \]
    with $(i,j) = \diagcol{\eta}{j} \in \orbipartinds{p}{q}$, $j \geq 2$, and shifted
    columns
    \[
    S = \{\diagcol{1}{c_1}, \diagcol{2}{c_2}, \dots, \diagcol{\eta}{c_{\eta}}\}
    \text{ with }2 \leq c_1 = c_2 \leq \dots \leq c_{\eta} \leq j-1,
    \]
    where in case of $\eta=1$ the last condition reduces to $2\leq
    c_1$ and we additionally require $j=c_1+1$.
  \end{myitemize}
  This system of constraints is non-redundant. The corresponding separation
  problem can be solved in time $\bigo{pq}$.
\end{theorem}

For the result on the completeness of the description, see
Proposition~\ref{prop:orbipart:descr}, for the question of redundancy see
Proposition~\ref{prop:orbipart:facets}, and for the separation algorithm
see Corollary~\ref{cor:sepSCI}. Note that the SCI with shifted column
$\{(1,1)\}$ and bar $\{(2,2)\}$ defines the same facet of $\orbipart{p}{q}$
as the nonnegativity constraint $x_{2,1}\geq 0$.

\begin{theorem}\label{thm:resultsPack}
  The packing orbitope $\orbipackG{p}{q}{\symgr{q}}$ $($for $p \geq q \geq
  2)$ with respect to the symmetric group $\symgr{q}$ equals the set of all
  $x \in \R^{\orbipartinds{p}{q}}$ that satisfy the following linear
  constraints:
  \begin{myitemize}
  \item the row-sum inequalities $x(\row{i}) \leq 1$ for all $i \in
    \ints{p}$,
  \item the nonnegativity constraints $x_{ij} \geq 0$ for all
    $(i,j) \in \orbipartinds{p}{q} \setminus \setdef{(j,j)}{j < q}$,
  \item the shifted column inequalities $x(B)-x(S) \leq 0$ for all bars
    \[
    B = \{(i,j), (i,j+1), \dots, (i,\min\{i,q\})\}
    \]
    with $(i,j) = \diagcol{\eta}{j} \in \orbipartinds{p}{q}$, $j \geq 2$,
    and shifted columns
    \[
    S = \{\diagcol{1}{c_1}, \diagcol{2}{c_2}, \dots,\diagcol{\eta}{c_{\eta}}\}
    \text{ with } c_1 = c_2 \leq \dots \leq c_{\eta} \leq j-1,
    \]
    where in case of $\eta=1$  we additionally require $j=c_1+1$.
  \end{myitemize}
  This system of constraints is non-redundant. The corresponding separation
  problem can be solved in time $\bigo{pq}$.
\end{theorem}

For the result on the completeness of the description, see
Proposition~\ref{prop:orbipack:descr}, for the question of redundancy see
Proposition~\ref{prop:orbipack:facets}, and for the separation algorithm
see Corollary~\ref{cor:sepSCI}.

% +++++++++++++++++++++++++++++++++++++++++++
% Remarks
% +++++++++++++++++++++++++++++++++++++++++++

\section{Concluding Remarks}
\label{sec:ClosingRemarks}

We close with some remarks on the technique used in the proof of
Proposition~\ref{prop:orbipart:descr}, on the combination of SCIs and
clique-inequalities for the graph-coloring problem, and on full and covering
orbitopes.

\subsubsection*{The Proof Technique.}

Our technique to prove Proposition~\ref{prop:orbipart:descr} can be
summarized as follows. Assume a polytope~$Q \subset \R^n$ is described by
some (finite) system~$\mathcal{Q}$ of linear equations and inequalities.
Suppose that~$\mathcal{Q}'$ is a subsystem of~$\mathcal{Q}$ for which it is
known that~$\mathcal{Q}'$ defines an integral polytope $Q'\supseteq Q$. One
can prove that~$Q$ is integral by showing that every vertex~$x^{\star}$
of~$Q$ is a vertex of~$Q'$ in the following way. Here we call a basis (with
respect to~$\mathcal{Q}$) of~$x^{\star}$ \emph{reduced} if it contains as
many constraints from~$\mathcal{Q}'$ as possible:
\begin{myenumerate}
\item\label{item:remarks:weight} Starting from an arbitrary reduced
  basis~$\mathcal{B}$ of~$x^{\star}$, construct iteratively a reduced
  basis~$\mathcal{B}^\star$ of~$x^{\star}$ that satisfies some properties
  that are useful for the second step.
\item Under the assumption that~$\mathcal{B}^\star \not\subseteq
  \mathcal{Q'}$, modify~$x^{\star}$ to some $\Tilde{x}\not= x^{\star}$ that
  also satisfies the equation system corresponding to~$\mathcal{B}^\star$
  (contradicting the fact that~$\mathcal{B}^\star$ is a basis).
\end{myenumerate}
(In our proof of Proposition~\ref{prop:orbipart:descr},
Step~\eqref{item:remarks:weight} was done by showing that a reduced basis
of ``minimal weight'' has the desired properties.)

Such a proof is conceivable for every 0/1-polytope~$Q$ by choosing~$Q' =
[0,1]^n$ as the whole $0/1$-cube and~$\mathcal{Q}'$ as the set of the $2n$
trivial inequalities $0 \leq x_i \leq 1$, for $i = 1, \dots, n$ (if
necessary, modifying~$\mathcal{Q}$ in order to contain them all).

We do not know whether this kind of integrality proof has been used in the
literature. It may well be that one can interpret some of the classical
integrality proofs in this setting. Anyway, it seems to us that the
technique might be useful for other polytopes as well.

\subsubsection*{The Graph-Coloring Problem.}

As mentioned in the introduction, for concrete applications like the graph
coloring problem one can (and probably has to) combine the polyhedral
knowledge on orbitopes with the knowledge on problem specific polyhedra.
We illustrate this by the example of clique inequalities for the graph
coloring model~\eqref{eq:intro:model} described in the introduction.

\begin{figure}
  \centering
  \newcommand{\mystyle}[1]{\footnotesize {#1}}
  \psfrag{j}{\mystyle{$j$}}
  \includegraphics[height=.17\textheight]{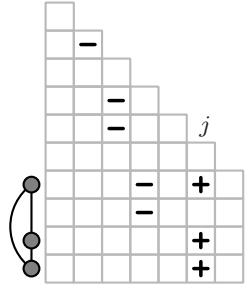}
  \caption{Combination of a clique inequality and an SCI.}
  \label{fig:cliquesci}
\end{figure}

Fix a color index $j \in \ints{C}$. If $W \subseteq V$ is a clique in the
graph $G = (V,E)$, then clearly the inequality $\sum_{i\in W}x_{ij} \leq 1$
is valid.  In fact, the strengthened inequalities $\sum_{i\in W}x_{ij} \leq
y_j$ are known to be facet-defining for the convex hull of the solutions
to~\eqref{eq:intro:model}, see~\cite{CollMDZ02}. Suppose that $S \subset
\orbipartinds{\card{V}}{C}$ is a shifted column and that we have $\eta \leq
\card{S}$ for all $\diagcol{\eta}{j} = (i,j)$ with $i \in W$. Then the
inequality
\[
\sum_{i\in W}x_{ij} - x(S) \leq 0
\]
is valid for all solutions to the model obtained
from~\eqref{eq:intro:model} by adding inequalities~\eqref{eq:SymmetryBreak}
(which are all ``column inequalities'' in terms of orbitopes), see
Figure~\ref{fig:cliquesci}. The details and a computational study will be
the subject of a follow-up paper.

\subsubsection*{Full and Covering Orbitopes.}

As soon as one starts to consider 0/1-matrices that may have more than one
$1$-entry per row, things seem to become more complicated.

With respect to cyclic group actions, we loose the simplicity of the
characterizations in Observation~\ref{obs:charvert}. The reason is that the
matrices under investigation may have several equal nonzero columns. In
particular, the lexicographically maximal column may not be unique.

With respect to the action of the symmetric group, we still have the
characterization of the representatives as the matrices whose columns are
in non-increasing lexicographic order (see Part~\ref{obs:charvert:sym} of
Observation~\ref{obs:charvert}). The structures of the respective full and
covering orbitopes, however, become much more complicated. In particular,
we know from computer experiments that several powers of two arise as
coefficients in the facet-defining inequalities. This increase in
complexity is reflected by the fact that optimization of linear functionals
over these orbitopes seems to be more difficult than over packing and
partitioning orbitopes (see the remarks at the end of
Section~\ref{sec:OptimizingOverOrbitopes}).

Let us close with a comment on our choice of the set of representatives as
the maximal elements with respect to a lexicographic ordering (referring to
the row-wise ordering of the components of the matrices). It might be that
the difficulties for full and covering orbitopes mentioned in the previous
paragraph can be overcome by the choice of a different system of
representatives. The choice of representatives considered in this paper,
however, seems to be appropriate for the packing and partitioning cases.
\smallskip

Whether the results presented in this paper are useful in practice will
turn out in the future. In any case, we hope that the reader shares our
view that orbitopes are neat mathematical objects. It seems that symmetry
strikes back by its own beauty, even when mathematicians start to fight it.

\subsubsection*{Acknowledgment.}

We thank the referees for their work. In particular, we are indepted to one
of them for several insightful and constructive remarks, including the
proof of Theorem~\ref{thm:orbipackcycl} that we have in this final version.
We furthermore thank Yuri Faenza, Andreas Loos, and Matthias Peinhardt for
helpful comments.

%\bibliographystyle{mod_siam}
%\bibliography{orbitopes}

\end{document}